# A model structure for locally constant factorization algebras

Victor Carmona[a], Ramon Flores[b] and Fernando Muro[c*]

January 2026

[a]Max-Planck Institut for Mathematics in the Sciences, Leipzig, Germany
*Email address:* victor.carmona@mis.mpg.de *(Corresponding author)*

[b]Universidad de Sevilla, Departamento de Geometría y Topología, Sevilla, Spain
*Email address:* ramonjflores@us.es

[c]Universidad de Sevilla, Departamento de Álgebra, Sevilla, Spain
*Email address:* fmuro@us.es

**Abstract**

Several model structures related to the homotopy theory of locally constant factorization algebras are constructed. This answers a question raised by D. Calaque in his habilitation thesis about the existence of such a model category. Our methods also solve a problem related to cosheafification and factorization algebras identified by O. Gwilliam - K. Rejzner in the locally constant case.

**Keywords:** *Factorization algebras, factorization homology, model categories, homotopy theory, Bousfield localization, operad.*

*The authors were partially supported by the Spanish Ministry of Economy grant MTM2016-76453-C2-1-P (AEI/FEDER, UE); the Andalusian Ministry of Economy and Knowledge and the Operational Program FEDER 2014-2020 grant US-1263032; the Spanish Ministry of Science and Innovation grant PID2020-117971GB-C21 and Junta de Andalucía grant FQM-213. V.C. was also partly supported by Spanish Ministry of Science, Innovation and Universities grant FPU17/01871.

# Contents

# 1 Introduction

*Factorization algebras* are a relatively new algebraic structure arising in a variety of fields: $(i)$ they are the topological analogues of Beilinson-Drinfeld's chiral algebras [8] (appearing in their work on Conformal Field Theories, also denoted CFTs); $(ii)$ they were used by Gaitsgory and Lurie to prove Weil's Conjecture over finite fields; $(iii)$ they provide strong manifold invariants, gathered under the name *factorization homology* [1, 2, 40], reminiscent of the classical configuration spaces of summable labels of McDuff, Salvatore and Segal [46]; $(iv)$ they are central to Costello-Gwilliam's work on Perturbative Quantum Field Theory (pQFT) [21, 22]; $(v)$ among factorization algebras, the locally constant ones are tightly connected to Topological Field Theories (TFTs) and iterated loop spaces [14, 27, 30, 40]...

Being more precise, factorization algebras are algebraic devices which satisfy (homotopical) codescent conditions. Their underlying algebraic structure is simple to describe and receives the name of prefactorization algebra structure.

**Definition 1.1.** A *prefactorization algebra* $\mathcal{F}$ over a space $X$ valued in chain complexes is given by the following data:

- a chain complex $\mathcal{F}(U)$ for each open subset $U$ in $X$,
- a chain map $\mathcal{F}(U) \to \mathcal{F}(V)$ for any inclusion of open subsets $U \subset V$ in $X$,
- a product map $\mathcal{F}(U_1) \otimes \cdots \otimes \mathcal{F}(U_m) \to \mathcal{F}(Z)$ for each inclusion of disjoint open subsets $U_1 \sqcup \cdots \sqcup U_m \subset Z$ in $X$, subject to associativity, unitality and equivariancy conditions (see §2 for more details).

The codescent conditions that a prefactorization algebra has to satisfy in order to be a factorization algebra (Definition 2.6) are designed to provide a local-to-global principle that accounts for multilocal data, i.e. simultaneous information around finite families of points. Manifold calculus [11] showed how important this principle is, since global objects like embedding spaces are multilocal in this sense under appropriate assumptions.

Observe that prefactorization algebras admit an operadic description, as they are algebras over a discrete operad $\mathsf{Disj}(X)$ (see Definition 2.3). Hence, it is possible, under suitable conditions, to endow the category of prefactorization algebras with a Quillen model structure (e.g. [43, §9.7]) which allows one to access its homotopy theory. However, the local-to-global property of factorization algebras obstructs the passage from this Quillen model structure to another one that presents the homotopy theory of factorization algebras. Also, the lack of a model category presentation seriously complicates computations. In fact, it was a question raised by Calaque in [13, Remark 2.1.15] if such a model category exists and to the best of our knowledge no satisfactory answer has been provided so far.

Assuming local constancy (Definition 2.10)–informally, insensitivity to the size of discs–, we find a model category presenting the homotopy theory of such factorization algebras, thereby answering Calaque's question.

**Theorem A.** *[Theorem 4.25] The category of prefactorization algebras over a smooth manifold valued in a suitable symmetric monoidal model category $\mathcal{V}$ admits a model structure such that the bifibrant objects are the (projectively bifibrant) locally constant factorization algebras and the weak equivalences between them are just the objectwise equivalences. Namely, this holds if $\mathcal{V} = \mathcal{C}\mathsf{h}(R)$ is the projective model structure on chain complexes over a commutative ring $R$ containing $\mathbb{Q}$.*

Mostly for simplicity in the exposition, we focus on smooth manifolds without boundary, although our methods do apply to more general settings. Remarkably, they work for smooth manifolds with boundary, conically smooth stratified manifolds, treated in [5], or topological manifolds (see §5 for these generalizations). Moreover, we also explore what happens with less hypotheses on the base symmetric monoidal model category. In such cases we obtain left semimodel categories (see [6, §1] or [16, §2]) instead of complete model categories. This subtlety comes from applying left Bousfield localizations at a set of maps in the absence of left properness (see [7, Theorem A] and [16, Theorem B]).

Along the way, we have found a variety of different Quillen equivalent model structures. The more remarkable one presents the homotopy theory of what we have called enriched factorization algebras (Definition 3.15 and Theorem 4.10). They should be seen as factorization algebras whose underlying algebraic structure is sensitive to the topology of embedding spaces. This idea is implemented by replacing the discrete operad $\mathsf{Disj}(X)$ by an operad $\mathcal{M}_X$ whose spaces of operations involve embedding spaces. The use of enriched factorization algebras is a fundamental step to obtain Theorem **A**. In fact, they are used to show that, in the locally constant case, the homotopical codescent conditions that a factorization algebra satisfies are equivalent to being completely characterized by a variant of factorization homology (see §3.2). This idea can be motivated by the following principle: *factorization homology produces fully extended topological quantum field theories* (see [47]). Furthermore, locally constant factorization algebras are expected to be tightly connected to topological field theories beyond this principle. For instance, one may consult the surveys [3, 23], and the references therein, for more on this connection.

It is also important to point out that, by dual analogy with sheaves, one should recognize two main sources of problems to impose codescent conditions on prefactorization algebras. On the one hand, cosheafification, as a formal machine to force precosheaves to satisfy a local-to-global principle, is much more elusive than sheafification. Delving into [45, Theorem 2.1], one can find that usual categories such as sets or Grothendieck abelian categories do permit cosheafification, but we only know its existence by means of abstract adjoint functor theorems, and hence there is no manageable expression for it. On the other hand, one must ensure that, when existing, cosheafification respects the additional algebraic structure that a prefactorization algebra has. This problem was actually identified by Gwilliam-Rejzner in [29, Remark 2.33]. To the best of our knowledge, there is no solution to these issues in general and Theorem **A** is the only known answer to Gwilliam-Rejzner's problem.

**Outline:** Let us summarize the content of this paper.

- The document begins with two mostly preliminary sections, §2 and §3. They contain a reminder of notions related to factorization algebras as well as fundamental lemmas about codescent properties of embedding spaces. Factorization homology (with context) is presented in §3.2, whose main goal is Theorem 3.19, which is an original contribution. This result gives a proof of [4, Theorem 2.15] in the setting of smooth manifolds only relying on the fundamental codescent property stated in Lemma 3.10 and the existence of good Weiss covers (Lemma 3.14). Such theorem is fundamental for technical reasons and it appears, without proof, in [4]; see Remark 3.17.

- A summary of a model categorical construction from [15, §4] that we will extensively apply can be found in §4. It is based on a generalization of a method due to Bousfield-Friedlander to localize model structures (see [12, Theorem A.7]).

- Model categories presenting the homotopy theory of enriched (lax) factorization algebras are discussed in §4.1, see Theorem 4.10 (resp. Theorem 4.7). These model categories are constructed via the construction in the previous item.

- The analogous model categories for locally constant (lax) factorization algebras, Theorem 4.25 (resp. Theorem 4.21 plus Proposition 4.24), are the main goal of §4.2. They are combinations of a left Bousfield localization, forcing local constancy, and a right Bousfield localization, forcing codescent conditions.

- Section §5 is devoted to the discussion of different settings where our results hold. The exposition highlights what ingredients are required, separating them from general abstract arguments. Remarkably, our methods are suitable for topological manifolds or conically smooth stratified manifolds (developed within [4, 5]).

- Finally, §A contains a discussion about partial monoidal categories which is required for the recognition of bifibrant objects in the (lax) lc-factorization model structure (see Lemma 4.22). On the other hand, §B establishes the left properness of certain model structures on the category of prefactorization algebras under some conditions.

**Conventions:**

- In this document, we make use of semimodel structures. These devices have been prominently applied to the study of $\infty$-categories of operadic algebras; see [50, 51] and references therein. There are several inequivalent notions (see [16, §2]), but we will only consider Spitzweck left semimodel categories [16, Definition 2.5]. For this reason, we will omit the adjective *Spitzweck*.
- We will adopt the following notation: symbols like $\mathtt{V}$, $\mathtt{M}$... denote ordinary categories, while $\mathcal{V}$, $\mathcal{M}$ will be reserved for (semi)model categories.
- Given a category $\mathtt{C}$ with a subcategory $\mathtt{W}$ of weak equivalences, $\mathtt{C}[\![\mathtt{W}^{-1}]\!]$ stands for its $\infty$-localization, provided it exists. When $\mathcal{M}$ is a (semi)model category, $\mathcal{M}_\infty$ refers to the associated $\infty$-localization, or in other words, to its underlying $\infty$-category.
- We will denote by $\mathcal{S}\mathrm{pc}$ the model category of simplicial sets with the Kan-Quillen model structure. The notation comes from the fact that its associated $\infty$-category, $\mathcal{S}\mathrm{pc}_\infty$, is the $\infty$-category of spaces/$\infty$-groupoids.
- By homotopy cosmos we mean a closed symmetric monoidal model category $\mathcal{V}$ equipped with a lax symmetric monoidal left Quillen functor $\mathcal{S}\mathrm{pc} \to \mathcal{V}$. We also ask the unit in $\mathcal{V}$ to be cofibrant to simplify the exposition, although one can work around this assumption, e.g., as in [17]. We fix a homotopy cosmos $\mathcal{V}$ and consider that everything is $\mathcal{V}$-enriched. Otherwise, it will be explicitly specified.
- To denote categories of enriched functors we use $[-,-]$, and for categories of lax monoidal enriched functors and monoidal natural transformations $\mathtt{LaxMon}[-,-]$.
- $\mathtt{Alg}_{\mathcal{O}}(\mathtt{V})$ (or simply $\mathtt{Alg}_{\mathcal{O}}$) denotes the category of $\mathcal{O}$-algebras in $\mathtt{V}$. When required and when it exists, $\mathtt{Alg}_{\mathcal{O}}(\mathtt{V})$ carries the projective model structure, denoted $\mathcal{A}\mathrm{lg}_{\mathcal{O}}(\mathcal{V})$. In particular, proj-cofibrant algebras makes reference to cofibrant objects in this model category.
- Let $\mathcal{O}$ be an operad. We will denote by $\overline{\mathcal{O}}$ the underlying category of $\mathcal{O}$, and by $\mathcal{O}\left[\begin{smallmatrix}\underline{c}\\ d\end{smallmatrix}\right]$ the object of multimorphisms from $\underline{c}$ to $d$. In other words, $\overline{\mathcal{O}} \hookrightarrow \mathcal{O}$ is the suboperad on unary operations. The induced restriction functor between algebras $\mathtt{Alg}_{\mathcal{O}}(\mathtt{V}) \to \mathtt{Alg}_{\overline{\mathcal{O}}}(\mathtt{V}) = [\overline{\mathcal{O}}, \mathtt{V}]$ is denoted by $A \mapsto \overline{A}$ and its left adjoint by ext.
- $\mathbb{0}$ and $\mathbb{1}$ will refer to the initial and the terminal objects, always interpreted in the corresponding context. We denote by $\mathbb{I}$ the monoidal unit of $\mathcal{V}$.
- Given a cospan of categories $B \xrightarrow{f} C \xleftarrow{g} A$, we denote $f \downarrow g$ the ordinary slice category associated to it.

**Further assumptions on $\mathcal{V}$:** In order to ensure that the two ways to obtain $\infty$-categories of (pre)factorization algebras in the sense of [21] with values in $\mathcal{V}_\infty$ match, i.e. as a full subcategory of $\infty$-algebras over the operad of disjoint open subsets in $\mathcal{V}_\infty$ or as the $\infty$-localization of a homotopically replete full subcategory $\mathtt{Fact}_X(\mathtt{V})$ of $\mathtt{Alg}_{\mathsf{Disj}(X)}(\mathtt{V})$ (see Definitions 2.3 and 2.6), we will impose the following restrictive assumption on $\mathcal{V}$ in the sequel:

**Hypothesis 1.2.** *The closed symmetric monoidal model category $\mathcal{V}$ satisfies the following list of properties (see [43, 44]):*

- *$\mathcal{V}$ is tractable, i.e. it is cofibrantly generated and the generating sets of (trivial) cofibrations have cofibrant domains.*
- *$\mathcal{V}$ is combinatorial and weak equivalences are closed under filtered colimits.*
- *$\mathcal{V}$ is symmetric* h*-monoidal, symmetroidal and symmetric flat.*

This strong restriction on $\mathcal{V}$ implies that any operad is admissible [43, Theorem 5.11], admits rectification [43, Theorem 7.5] and has the correct $\infty$-category of algebras by [33, Theorem 1.1] and [43, Theorem 7.11]. The most basic examples are chain complexes over a commutative ring $R$ containing $\mathbb{Q}$ with the projective model structure or symmetric spectra with the positive flat model structure. For further examples, consult [43, 44].

It might be possible to consider more general homotopy cosmoi $\mathcal{V}$ satisfying analogues of the previous results which suffice for our purposes. For example, if we restrict ourselves to use $\Sigma$-cofibrant operads, the previous results hold true with virtually no assumption on $\mathcal{V}$ at the expense of working with semimodel categories; see [17, 50] and [51, Theorem 7.3.1]. The main obstruction to follow this alternative is that $\mathrm{Disj}(X)$, and other relevant operads for our constructions, such as $\mathrm{Discs}(X)$ or $\mathcal{D}\mathrm{iscs}_n$, are not $\Sigma$-free because of the presence of the empty set $\varnothing$ as a color. Notice that the value $\mathcal{F}(\varnothing)$ of a prefactorization algebra $\mathcal{F}$, as considered in [21], carries a unital commutative algebra structure and those behave badly in general symmetric monoidal model categories (cf. [49]). We believe that our arguments can be carried out in greater generality regarding $\mathcal{V}$ if we do not consider the empty set as a color in all the operads we use. This would not be an important restriction in applications since, by weak monadicity, the unit map $\mathbb{I}_V \to \mathcal{F}(\varnothing)$ is an equivalence for any factorization algebra $\mathcal{F}$. In other words, as long as weak monadicity is imposed, the technicalities which arise from dealing with the commutative algebra $\mathcal{F}(\varnothing)$ should be avoidable.

## 2 A factorization overview

**Little $n$-discs algebras and factorization homology.** Algebras over the little $n$-discs operad $\mathbb{E}_n$ appear in many branches of mathematics. In fact, $\mathbb{E}_n$-algebras interpolate between (homotopy coherent unital[i]) associative algebras and (homotopy coherent unital) commutative algebras. Our interest in $\mathbb{E}_n$-algebras here relies on the fact that they may be used as coefficient systems for "homology theories of $n$-manifolds", also called *factorization homology* or *chiral homology*, as envisioned by [1]. We will go deeper into the lands of *factorization homology* in §3.2, but for the moment, let us just present an overview of some of the central ideas about it.

*Remark* 2.1. As discussed in [2] and [40, §5.5], factorization homology can be studied for general topological manifolds. For simplicity, we will assume that all manifolds are *smooth* unless otherwise specified.

Given a framed $n$-manifold $Z$ (i.e. a smooth manifold with a trivialization of its tangent bundle), one can define a right $\mathbb{E}_n$-module (right module over the monoid $\mathbb{E}_n$ for the $\circ$-product on symmetric collections of topological spaces $\Sigma\mathtt{Coll}(\mathtt{Top})$) using framed-embeddings of finite disjoint unions of the open $n$-disc $\mathbb{D}$ into $Z$. Recalling [36, Definition 1.8, Remark 1.9] and denoting by $\mathsf{Env}_{\otimes}(\mathbb{E}_n)$ the symmetric monoidal envelope of $\mathbb{E}_n$, this is just the datum of the $\mathtt{Top}$-enriched functor

$$\begin{aligned} \mathsf{Emb}^{\mathsf{fr}}(\star, Z) \colon \mathsf{Env}_{\otimes}(\mathbb{E}_n)^{\mathsf{op}} &\longrightarrow \mathtt{Top} \\ \{1, \dots, k\} &\longmapsto \mathsf{Emb}^{\mathsf{fr}}(\mathbb{D}^{\sqcup k}, Z) \end{aligned} .$$

[i] In this document, we consider only unital versions of $\mathbb{E}_n$ operads and their variations, i.e. they are considered to have one and only one nullary operation.

Then, given an $\mathbb{E}_n$-algebra $A$ in $\mathcal{V}$, *factorization homology of* $Z$ *with coefficients in* $A$ is simply the derived relative $\circ$-product of the right $\mathbb{E}_n$-module associated to $Z$ and the left-module $A$:

$$\int_Z A = \mathsf{Emb}^{\mathsf{fr}}(\star, Z) \overset{\mathbb{L}}{\underset{\mathbb{E}_n}{\circ}} A\,.$$ [ii]

This definition yields a pairing

$$\int : \begin{Bmatrix} \text{framed} \\ n\text{-manifolds} \end{Bmatrix} \times \mathcal{A}\mathsf{lg}_{\mathbb{E}_n}(\mathcal{V})_\infty \longrightarrow \mathcal{V}_\infty, \qquad (Z, A) \longmapsto \int_Z A$$

fulfilling a list of axioms close in spirit to Eilenberg-Steenrod axioms for ordinary homology theory (see [2, 28] for a fully detailed account):

(i) (*dimension*) there is a natural equivalence $\int_{\mathbb{R}^n} A \simeq \mathsf{fgt}(A)$, where $\mathsf{fgt}$ denotes the forgetful functor $\mathcal{A}\mathsf{lg}_{\mathbb{E}_n}(\mathcal{V})_\infty \to \mathcal{V}_\infty$;

(ii) (*weak monoidality*) the functor $Z \mapsto \int_Z A$ is lax monoidal and the structure maps $\bigotimes_i \int_{Z_i} A \to \int_{\bigsqcup_i Z_i} A$ are equivalences for any set of framed $n$-manifolds $\{Z_i\}_i$;[iii]

(iii) ($\otimes$-*excision*) the functor $Z \mapsto \int_Z A$ sends collar-gluings of framed $n$-manifolds to relative tensor products, i.e.

$$Z \cong Z_l \underset{Z_c \times \mathbb{R}}{\cup} Z_r \xrightarrow{\int_\star A} \left(\int_{Z_l} A\right) \underset{\left(\int_{Z_c \times \mathbb{R}} A\right)}{\overset{\mathbb{L}}{\otimes}} \left(\int_{Z_r} A\right) \simeq \int_Z A\ .$$

In fact, using handlebody decompositions, Ayala-Francis in [2] showed that these properties uniquely characterize *factorization homology*, just as Eilenberg-Steenrod axioms characterize ordinary homology, by constructing $\int_Z A$ over any $Z$ by combining (i)-(iii).

*Remark* 2.2. As mentioned before, this scheme can be pushed further and works for much more general kinds of manifolds if one modifies the algebraic input accordingly. For instance, one can define factorization homology for general smooth $n$-manifolds by coupling them with "unoriented $\mathbb{E}_n$-algebras", i.e. algebras over the operad $\mathcal{E}_n$ from Notation 3.4.

[ii] Recall that we assume $\mathcal{V}$ to come equipped with a lax symmetric monoidal functor from simplicial sets.

[iii] In (ii), it is written "any *set* of framed manifolds" instead of *finite set*. Adding the adjective finite to this axiom would restrict the uniqueness of factorization homology to *finitary framed manifolds* (see [3]). However, requiring (ii) just on finite sets of manifolds and an additional axiom for filtered colimits yields the same concept as the one spelled out here.

One might interpret this uniqueness result by saying that there is an unique way to extend an $\mathbb{E}_n$-algebra A to a lax-monoidal functor from the symmetric monoidal (or sm-)Top-category of framed $n$-manifolds and framed embeddings if one requires axioms (i)-(iii). So, in a sense, by prescribing an algebraic structure over $\mathbb{R}^n$ or equivalently $\mathbb{D}$, one can create a kind of continuous multiplicative precosheaf on all framed $n$-manifolds. Moreover, this extension controls important algebraic and manifold invariants, e.g. if A is an unital associative algebra, $\int_{\mathbb{S}^1} A$ computes its Hochschild homology [2, Theorem 3.19].

Now, fixing a framed $n$-manifold X and restricting the functor $\int_\star A$ to the poset Open(X) of open subsets of X, one obtains an actual precosheaf over X

$$\mathcal{F}_A \colon \mathrm{Open}(X) \longrightarrow \mathcal{V}.$$

The role of (pre)factorization algebras is to replace these precosheaves associated to $\mathbb{E}_n$-algebras by more general coefficient systems. The idea is to produce generalized "homology theories of manifolds" by replacing "local systems" as coefficients by "(co)sheaves".

**(Pre)factorization algebras.** In this subsection, we only require that X is a topological space since the definitions make sense in this level of generality.

The notion of prefactorization algebra is quite elementary and it does not even require operadic machinery to specify what it is. In fact, a *prefactorization algebra* $\mathcal{F}$ on X (with values in V) is given by:

(a) an object $\mathcal{F}(Z) \in V$ for each open subset $Z \subseteq X$, and

(b) a product map $\mu\left[{(Z_i)_i \atop Z}\right] \colon \bigotimes_i \mathcal{F}(Z_i) \to \mathcal{F}(Z)$ for any finite disjoint union of open subsets $Z_i \subseteq Z$, "independent of the order" of $Z_i$'s,[iv] subject to unitality conditions imposing that $\mu\left[{(Z) \atop Z}\right] = \mathrm{id}_{\mathcal{F}(Z)}$ and a simple associativity condition stating that compositions of product maps yield product maps [27, Definition 12]. See Figure 2.1.

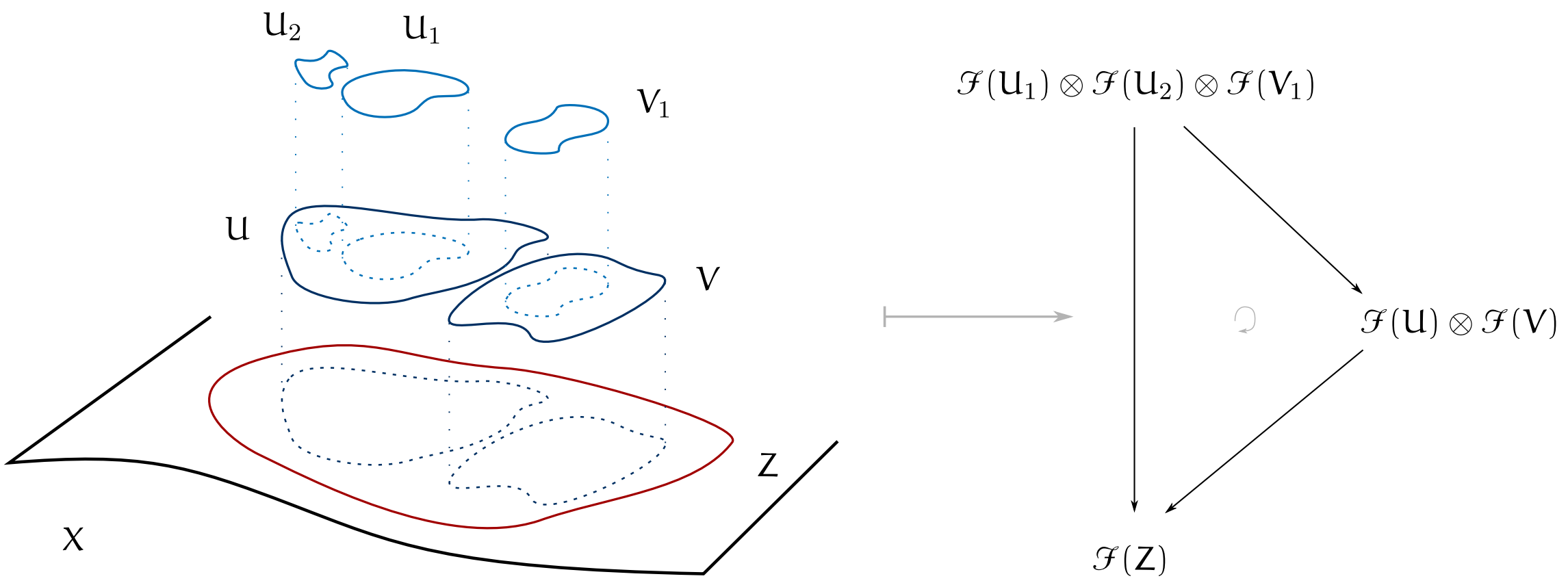

Figure 2.1: Associativity condition for prefactorization algebras.

Despite the simplicity of the previous definition, it is useful to bring operadic language into the picture.

[iv] In particular, there are maps $\mu\left[{\emptyset \atop Z}\right] \colon \mathbb{1}_V \to \mathcal{F}(Z)$ and $\mu\left[{(Z') \atop Z}\right] \colon \mathcal{F}(Z') \to \mathcal{F}(Z)$.

**Definition 2.3.** The operad in sets of *disjoint open subsets* in $X$, $\mathsf{Disj}(X)$, is the operad with colors $\mathrm{col}(\mathsf{Disj}(X)) = \mathrm{ob}(\mathsf{Open}(X))$ the open subsets of $X$ and operations

$$\mathsf{Disj}(X)\begin{bmatrix}(Z_i)_i\\ Z\end{bmatrix} = \begin{cases} \mathbb{1} & \text{if } \bigsqcup_i Z_i \subseteq Z, \\ \mathbb{0} & \text{otherwise.}\end{cases}$$

**Definition 2.4.** A *prefactorization algebra* on $X$ is a $\mathsf{Disj}(X)$-algebra. We denote by $\mathtt{pFact}_X(\mathtt{V})$ the category of prefactorization algebras on $X$.

*Remark* 2.5. Since $\mathsf{Disj}(X)$ is the operad associated to the partial monoidal category $(\mathsf{Open}(X), \sqcup)$ (see §A),

$$\mathtt{pFact}_X(\mathtt{V}) = \mathtt{Alg}_{\mathsf{Disj}(X)}(\mathtt{V}) = \mathtt{LaxMon}\left[\mathsf{Open}(X), \mathtt{V}\right].$$

This fact illustrates how prefactorization algebras are a sort of multiplicative precosheaves.

Factorization algebras are to prefactorization algebras what sheaves are to presheaves; that is, they are prefactorization algebras satisfying "gluing properties" that account for multilocal data in $X$, i.e. local information around any finite family of points simultaneously. Motivated by manifold calculus (see for example [11]), one may encode multilocal information using the *Weiss topology*.

**Definition 2.6.** Let $\mathcal{F}$ be a prefactorization algebra on $X$.

(a) A *Weiss cover* of an open subset $U \subseteq X$ is a family of open subsets $(U_i)_{i\in I}$ of $U$ such that any non-empty finite subset $T \subset U$ is contained in one of them, $T \subset U_i$.

(b) $\mathcal{F}$ is a *factorization algebra* on $X$ if it satisfies:

(i) (*Weiss homotopy codescent*) for any Weiss cover $(U_i)_{i\in I}$ of an open subset $U \subseteq X$, the map

$$\operatorname*{hocolim}_{S\subseteq I} \mathcal{F}(U_S) \longrightarrow \mathcal{F}(U)$$

is an equivalence, where the homotopy colimit runs over non-empty finite subsets $S$ of $I$ and we have adopted the standard notation $U_S = \bigcap_{i\in S} U_i$;

(ii) (*weak monadicity*) for any finite collection of disjoint open subsets $\{V_j\}_j$ of $X$, the structure map

$$\mu\colon \bigotimes_j \mathcal{F}(V_j) \longrightarrow \mathcal{F}\Big(\bigsqcup_j V_j\Big)$$

is an equivalence. In particular, $\mathbb{1}_{\mathtt{V}} \to \mathcal{F}(\varnothing)$ is an equivalence.

If $\mathcal{F}$ only satisfies (i) (Weiss homotopy codescent), it is said to be a *lax factorization algebra* on $X$.

The full subcategory of $\mathtt{pFact}_X(\mathtt{V})$ spanned by factorization algebras (resp. lax factorization algebras) on $X$ is denoted $\mathtt{Fact}_X(\mathtt{V})$ (resp. $\mathtt{laxFact}_X(\mathtt{V})$). Both are seen as relative categories by equipping them with the class of colorwise equivalences.

*Remarks* 2.7.

1. Given a Weiss cover $\mathscr{U} = \{U_i\}_{i\in I}$ of $U$, it is common to present the Weiss homotopy codescent for $\mathscr{U}$ by constructing the Čech-type $\Delta^{op}$-diagram

$$C_\bullet(\mathscr{U},\mathscr{F}): \quad \cdots \quad \coprod_{i,j,k\in I} \mathscr{F}(U_{ijk}) \substack{\longrightarrow\\ \leftarrow\text{-}\text{-}\text{-}\\ \longrightarrow\\ \leftarrow\text{-}\text{-}\text{-}\\ \longrightarrow} \coprod_{i,j\in I} \mathscr{F}(U_{ij}) \substack{\longrightarrow\\ \leftarrow\text{-}\text{-}\text{-}\\ \longrightarrow} \coprod_i \mathscr{F}(U_i)$$

and asking for the canonical map $\operatorname{hocolim} C_\bullet(\mathscr{U},\mathscr{F}) \to \mathscr{F}(U)$ to be an equivalence (see [21, Definition 6.1.4] or [37, Remark 2.57] for a careful analysis).

2. Some authors prefer to combine Weiss homotopy codescent and weak monadicity into a single gluing property: *homotopy codescent for factorizing covers*. Recall that an open cover $\mathscr{U}$ of $U \subseteq X$ is *factorizing* if for any finite subset $S \subseteq U$, there is a finite list of pairwise disjoint open subsets in the cover $\{U_j\}_j \subseteq \mathscr{U}$ (possibly indexed by a finite subset different from $S$) so that $S \subseteq \bigsqcup_j U_j$. Note that any factorizing cover $\mathscr{U}$ generates a Weiss cover by considering

$$\mathscr{U}' = \{U_1 \sqcup \cdots \sqcup U_k : U_i \in \mathscr{U}\}.$$

Combining this fact with the Čech-type presentation of homotopy codescent for Weiss covers, it is possible to define a $\Delta^{op}$-diagram $C^{\otimes}_\bullet(\mathscr{U},\mathscr{F})$ accounting for homotopy codescent for the factorizing cover $\mathscr{U}$ (see [13, §2.1.2] and [27, §4.1]). Roughly speaking, one decomposes the finite intersections of elements in $\mathscr{U}'$ as disjoint unions of finite intersections of elements in $\mathscr{U}$ and takes the resulting disjoint unions outside of $\mathscr{F}$.

An important property of factorization algebras is that they can be defined on a nice basis of the topology of $X$ and then extended to all open subsets, just as sheaves can be defined by descent. First, note that given a basis $\mathscr{B}$ of the topology of $X$ and considering the full suboperad of $\mathsf{Disj}(X)$ whose objects are elements of $\mathscr{B}$, one can define the category $\mathtt{pFact}_{\mathscr{B}}(\mathtt{V})$ of prefactorization algebras on $\mathscr{B}$. Assuming that $\mathscr{B}$ is a *factorizing basis* of $X$, i.e. a basis for the topology of $X$ so that for any $U \in \mathsf{Open}(X)$, $\mathscr{B}_U = \{Z \in \mathscr{B} : Z \subseteq U\}$ is a factorizing cover of $U$ (see [13, §2.1.2]), it makes sense to define $\mathtt{Fact}_{\mathscr{B}}(\mathtt{V})$ and in fact, one has:

**Proposition 2.8.** *If $\mathscr{B}$ is a factorizing basis of $X$, the restriction functor*

$$\mathtt{Fact}_X(\mathtt{V})[\![\mathtt{Eq}^{-1}]\!] \longrightarrow \mathtt{Fact}_{\mathscr{B}}(\mathtt{V})[\![\mathtt{Eq}^{-1}]\!]$$

*is an equivalence of $\infty$-categories.*

*Proof idea.* See [13, Theorem 2.1.9] or [37, §4.3]. Restriction along $\mathsf{Disj}(\mathscr{B}) \hookrightarrow \mathsf{Disj}(X)$ admits a homotopy inverse construction by setting $\mathrm{ext}(\mathscr{F}_{|})(U) = \operatorname{hocolim} C^{\otimes}_\bullet(\mathscr{B}_U, \mathscr{F}_{|})$ for any $\mathscr{B}$-factorization algebra $\mathscr{F}_{|}$. □

Given a continuous map $f: X \to Y$, there is an associated morphism of set-valued operads $f^{-1}: \mathsf{Disj}(Y) \to \mathsf{Disj}(X)$, $U \mapsto f^{-1}(U)$ and restriction along $f^{-1}$ defines the *pushforward of (pre/lax-)factorization algebras along* $f$

$$\begin{array}{ccc} \mathtt{pFact}_X(\mathtt{V}) & \xrightarrow{f_*} & \mathtt{pFact}_Y(\mathtt{V}) \\ \uparrow & & \uparrow \\ (\mathtt{lax})\mathtt{Fact}_X(\mathtt{V}) & \xrightarrow[f_*]{} & (\mathtt{lax})\mathtt{Fact}_Y(\mathtt{V}) \end{array},$$

where $f_*\mathcal{F}(U) = \mathcal{F}(f^{-1}(U))$. In fact, $f_*$ preserves objectwise equivalences and so, it induces functors between the associated $\infty$-categories of (pre/lax-)factorization algebras. Pushforward along $X \to \{*\}$ evaluates the (pre/lax-)factorization algebra on X. This object is alternatively referred to as *global cosections* of $\mathcal{F}$ and denoted $\Gamma(\mathcal{F}) \equiv \Gamma(X, \mathcal{F})$.

**Locally constant factorization algebras.** When X is a (smooth) $n$-manifold, the collection of open subsets of X diffeomorphic/homeomorphic to $\mathbb{R}^n$ (*open discs* in X) is enough to control factorization algebras in X. In fact, the collection of *bounded geodesically convex open subsets* of X (for some choice of riemannian metric)[v] is comprised of open discs and yields a factorizing basis of X (see [11, 28]). Hence, by Proposition 2.8, one can study factorization algebras by restricting them to open discs in X. This fact motivates the following:

**Definition 2.9.** The full suboperad Discs(X) of Disj(X) is the one whose colors are the open subsets of X diffeomorphic to finite disjoint unions of copies of $\mathbb{R}^n$. $E_X$ denotes the full suboperad of Disj(X) spanned by open subsets diffeomorphic to $\mathbb{R}^n$.

Recall that given an unoriented $\mathbb{E}_n$-algebra A valued in $\mathcal{V}$, i.e. an algebra over the operad $\mathcal{E}_n$ from Notation 3.4, we saw that restricting the factorization homology functor $\int_\star A$ to Open(X) gives rise to a multiplicative precosheaf $\mathcal{F}_A\colon \mathrm{Open}(X) \to \mathcal{V}$. By definition, $\mathcal{F}_A$ is actually a prefactorization algebra on X. While it might be complicated to explicitly determine $\mathcal{F}_A(U) = \int_U A$ for a general open subset $U \subseteq X$ in terms of A, the following two observations are, as we will see, enough to characterize $\mathcal{F}_A$:

- Any inclusion of open discs $D \subseteq D'$ is sent to an equivalence $\mathcal{F}_A(D) \simeq \mathcal{F}_A(D')$ since both objects are just $\int_{\mathbb{R}^n} A \simeq A$.
- By the dimension axiom and weak monadicity,
$$\mathcal{F}_A(D_1 \sqcup \cdots \sqcup D_k) \simeq A^{\otimes k}$$
for any finite disjoint union of open discs $D_i \subseteq X$, so the restriction of $\mathcal{F}_A$ to Discs(X) and $E_X$ is readily determined by A.

The first property can be seen as a "scale independence" of $\mathcal{F}_A$. It is important to identify the subcategory of (pre)factorization algebras satisfying it.

**Definition 2.10.** A (pre)factorization algebra $\mathcal{F}$ on X is called *locally constant* if for any inclusion $D \subset D'$ of open subsets diffeomorphic to finite disjoint unions of copies of $\mathbb{R}^n$ which is a bijection on connected components, the associated map $\mathcal{F}(D) \to \mathcal{F}(D')$ is an equivalence. We denote by $\mathtt{Fact}_X^{\mathtt{lc}}(\mathtt{V})$ (resp. $\mathtt{pFact}_X^{\mathtt{lc}}(\mathtt{V})$) the category of locally constant (pre)factorization algebras on X.

*Remark* 2.11. Note that pushforward along arbitrary smooth maps $f\colon X \to Y$ might not preserve local constancy. See [27, §6.1].

[v]Recall that a geodesically convex subset of a riemannian manifold is a subset so that any pair of distinct points in the subset can be connected by a unique minimal geodesic segment and this segment belongs completely to the subset.

**Notation 2.12.** We will use the superscript lc to refer to any kind of locally constant algebras, not just (pre)factorization algebras. For example, $\mathtt{Alg}^{\mathtt{lc}}_{E_X}(\mathtt{V})$ makes reference to the full subcategory of $E_X$-algebras spanned by those algebras sending unary morphisms in $E_X$ to equivalences.

One important consequence of local constancy is:

**Proposition 2.13.** *Locally constant factorization algebras on* $X$ *are completely characterized by their restriction to* $E_X$*, i.e. as locally constant* $E_X$*-algebras. More precisely,*

$$\mathtt{Fact}^{\mathtt{lc}}_X(\mathtt{V})[\![\mathtt{Eq}^{-1}]\!] \xrightarrow{\sim} \mathtt{Alg}^{\mathtt{lc}}_{E_X}(\mathtt{V})[\![\mathtt{Eq}^{-1}]\!].$$

*is an equivalence of* $\infty$*-categories.*[vi]

*Proof idea.* Let $\mathscr{B}$ be the factorizing basis of $X$ given by bounded geodesically convex open subsets (for some riemannian metric). By Proposition 2.8, it suffices to produce an equivalence replacing $\mathtt{Fact}^{\mathtt{lc}}_X(\mathtt{V})$ by $\mathtt{Fact}^{\mathtt{lc}}_{\mathscr{B}}(\mathtt{V})$. The result follows if we check that $A \in \mathtt{Alg}^{\mathtt{lc}}_{E_X}(\mathtt{V})$ (assume that $A$ is colorwise-cofibrant to be precise), seen as a locally constant $\mathscr{B}$-prefactorization algebra, satisfies factorizing homotopy codescent. This is follows from the following assertion: for any factorizing cover $\mathscr{U}$ of $U \in \mathscr{B}$ made of elements in $\mathscr{B}$, the canonical cocone $C^{\otimes}_{\bullet}(\mathscr{U}, A) \to A(U)$ is comprised of equivalences by local constancy. See [37, §5.3] for a complete proof. □

On the other hand, Lurie ([40, Theorem 5.4.5.9]) showed:

**Theorem 2.14.** *The construction* $\mathtt{Alg}_{\mathbb{E}_n}(\mathtt{V}) \to \mathtt{pFact}_{\mathbb{R}^n}(\mathtt{V})$, $A \mapsto \mathscr{F}_A$*, produces locally constant factorization algebras over* $\mathbb{R}^n$*. Moreover,* $A \mapsto \mathscr{F}_A$ *induces an equivalence of* $\infty$*-categories*

$$\mathtt{Alg}_{\mathbb{E}_n}(\mathtt{V})[\![\mathtt{Eq}^{-1}]\!] \xrightarrow{\sim} \mathtt{Fact}^{\mathtt{lc}}_{\mathbb{R}^n}(\mathtt{V})[\![\mathtt{Eq}^{-1}]\!].$$

*Remark* 2.15. Lurie's result can be split in two observations: (1) the prefactorization algebra given by factorization homology $\mathscr{F}_A$ satisfies Weiss homotopy codescent, and (2) $\mathbb{E}_n$ is the $\infty$-localization of $E_{\mathbb{R}^n}$ (as an $\infty$-operad) at all unary maps. Hence, (2) is saying that giving an $\mathbb{E}_n$-algebra is essentially the same thing as giving a locally constant $E_{\mathbb{R}^n}$-algebra. Then, Proposition 2.13 implies that such a datum produces a locally constant factorization algebra on $\mathbb{R}^n$ by forcing Weiss homotopy codescent and weak monadicity to extend the $E_{\mathbb{R}^n}$-algebra. In this sense, Theorem 2.14 may be seen as the coincidence of two lists of instructions to define invariants over arbitrary open subsets of $\mathbb{R}^n$ from an $\mathbb{E}_n$-algebra.

In §4 and §5, similar ideas will be applied to define model structures presenting the homotopy theory of locally constant (constructible) factorization algebras. In particular, we will provide a generalization of Theorem 2.14 together with an alternative proof to the one given by Lurie, which we hope is more conceptual.

# 3 Enriched factorization algebras

Motivated by manifold calculus as discussed in [11], we will define a variant of factorization algebras sensitive to the topology of embedding spaces, and we will see how they can be used to model locally constant factorization algebras.

[vi]For an alternative proof of this fact, see Theorem 4.25.

## 3.1 Embedding spaces and categories of manifolds

In the rest of this document, except for §5, we will consider that all manifolds are smooth. This choice simplifies the exposition at some points, but it is not essential as we will see in the cited section.

Let us start by defining spaces of embeddings between manifolds. Our exposition is based on [1, 36]. Recall that the compact-open topology on $\mathrm{Map}(X, Y)$ is designed to give, in the best scenario, an adjunction $(\star \times X) \dashv \mathrm{Map}(X, \star)$. If compact subsets of $X$ are closed (which is the case for all manifolds), this topology is determined by the basis with elements

$$\mathcal{V}_{(K,U)} = \left\{ f \in \mathrm{Map}(X, Y) : \ \mathrm{graph}(f_{|K}) \subseteq U \right\}$$

indexed by pairs $(K, U)$ with $K \subseteq X$ compact and $U \subseteq X \times Y$ open.

When smooth structures are considered, i.e. $X, Y$ are smooth manifolds, one equips the set of smooth maps $\mathrm{Map}_{\mathcal{C}^\infty}(X, Y)$ with a variation of this topology which takes into account proximity of derivatives as well. Fixing $1 \leqslant r < \infty$, one defines the *weak* $\mathcal{C}^r$*-topology* on $\mathrm{Map}_{\mathcal{C}^\infty}(X, Y)$ as the subspace topology inherited from the $r^{\mathrm{th}}$-jet prolongation map

$$j^r \colon \mathrm{Map}_{\mathcal{C}^\infty}(X, Y) \hookrightarrow \Gamma\big(X, J^r(X, Y)\big),$$

where the right hand side is the space of continuous sections of the jet projection $J^r(X, Y) \to X$, considered with the compact-open topology. For $r = \infty$, take the union of all these topologies. The weak $\mathcal{C}^r$-topology on $\mathrm{Map}_{\mathcal{C}^\infty}(X, Y)$ is not only interesting because it accounts for derivatives of maps, it makes the composition product

$$\circ \colon \mathrm{Map}_{\mathcal{C}^\infty}(Y, Z) \times \mathrm{Map}_{\mathcal{C}^\infty}(X, Y) \longrightarrow \mathrm{Map}_{\mathcal{C}^\infty}(X, Z)$$

continuous (see [20, 3.2.23]).

Inside $\mathrm{Map}_{\mathcal{C}^\infty}(X, Y)$, we find the *space of embeddings*

$$\mathrm{Emb}(X, Y) = \left\{ \begin{array}{c} f \in \mathrm{Map}_{\mathcal{C}^\infty}(X, Y) : \ f \text{ is a homeo onto its image} \\ \text{and its derivative is injective at any point} \end{array} \right\}.$$

*Remark* 3.1. We will be interested just on the homotopy type of embedding spaces and, as mentioned in [36, Remark 5.1][i], the weak $\mathcal{C}^r$ topologies for $1 \leqslant r \leqslant \infty$ on $\mathrm{Emb}(X, Y)$ are all (weak) homotopy equivalent. For that reason, we equip $\mathrm{Emb}(X, Y)$ with the weak $\mathcal{C}^1$-topology in the sequel.

Now, we want to add a manifold $X$ as context to the embedding spaces in a homotopically sensible way.

**Definition 3.2.** Given embeddings $\nu_Y \colon Y \hookrightarrow X \hookleftarrow Z : \nu_Z$, define the *space of embeddings with context* $\mathrm{Emb}_{/X}(\nu_Y, \nu_Z)$ or simply $\mathrm{Emb}_{/X}(Y, Z)$ via the homotopy pullback

$$\begin{array}{ccc} \mathrm{Emb}_{/X}(\nu_Y, \nu_Z) & \longrightarrow & \mathrm{Emb}(Y, Z) \\ \big\downarrow & {}^{h}\lrcorner & \big\downarrow {\scriptstyle \nu_{Z,*}} \\ \mathbb{1} & \xrightarrow[\nu_Y]{} & \mathrm{Emb}(Y, X) \end{array}$$

in $\mathtt{Top}$ (or equiv. in $\mathtt{Top}{\downarrow}\mathrm{Map}_{\mathcal{C}^\infty}(Y, Z)$; see [36, Remark 5.4]).

[i] In fact, consult [36, §6] for a detailed study of the homotopy type of embedding spaces.

*Remarks* 3.3. We have not specified the model structure considered on $\mathtt{Top}$ in the previous definition since we can use Hurewicz, Quillen or mixed structures indifferently. We ultimately care about the weak homotopy type of embedding spaces.

Note that spaces of embeddings with context inherit strict composition laws and units from those of plain embedding spaces when taking an appropriate model for the homotopy pullback (see [1, Definition V.10.1] and [36, §6] for details).

These embedding spaces may be organized to give topological categories and operads.

**Notation 3.4.** Let $\mathcal{M}\mathrm{fld}_n$ denote the sm-$\mathtt{Top}$-category of smooth $n$-manifolds and smooth embeddings with disjoint union as monoidal structure. In order to avoid set theoretic difficulties later on, we consider every manifold embedded in $\mathbb{R}^\infty$.

We denote by $\mathcal{D}\mathrm{iscs}_n$ the full sm-$\mathtt{Top}$-subcategory of $\mathcal{M}\mathrm{fld}_n$ spanned by finite disjoint unions of $n$-dimensional discs, i.e. $\mathbb{D}^{\sqcup k}$ for $0 \leqslant k < \infty$. We also denote by $\mathcal{E}_n$ the full $\mathtt{Top}$-suboperad of $\mathcal{D}\mathrm{iscs}_n$ on the color $\{\mathbb{D}\}$.

*Remark* 3.5. We are interested in the algebraic structures that the above operads (and variants of them) encode. For example, one should notice that $\mathcal{E}_n$-algebras are just unoriented $\mathbb{E}_n$-algebras as a consequence of Kister-Mazur theorem ([3, Proposition 2.6]).

Fixing a context $X \in \mathcal{M}\mathrm{fld}_n$, we have:

**Definition 3.6.** The $\mathtt{Top}$-category $\mathcal{M}\mathrm{fld}_{n/X}$ of $n$-manifolds with context $X$ has: $(i)$ as objects, pairs $(N, \nu_N)$ with $N \in \mathcal{M}\mathrm{fld}_n$ and $\nu_N \in \mathrm{Emb}(N, X)$; and $(ii)$ mapping spaces

$$\mathcal{M}\mathrm{fld}_{n/X}\big((T, \nu_T), (N, \nu_N)\big) = \mathrm{Emb}_{/X}(\nu_T, \nu_N).$$

Composition and units are induced from those of $\mathcal{M}\mathrm{fld}_n$. We often abuse notation and refer to an object in $\mathcal{M}\mathrm{fld}_{n/X}$ by its underlying manifold.

*Remark* 3.7. Note that the category $\mathcal{M}\mathrm{fld}_{n/X}$ has to be seen as a homotopical version of the ordinary slice construction.

The main drawback of $\mathcal{M}\mathrm{fld}_{n/X}$ is that it does not inherit the symmetric monoidal structure of $\mathcal{M}\mathrm{fld}_n$, nor a suitable *partial monoidal structure* in contrast with the discrete case (meaning $\mathrm{Open}(X)$). However, $\mathcal{M}\mathrm{fld}_{n/X}$ is still the underlying category of a natural operad $\mathcal{M}_X$.

**Definition 3.8.** The $\mathtt{Top}$-operad $\mathcal{M}_X$ is determined by:

**Obj:** its set of objects/colors is $\mathrm{col}(\mathcal{M}_X) = \mathrm{ob}(\mathcal{M}\mathrm{fld}_{n/X})$;

**Mor:** its spaces of multimorphisms fit into homotopy pullback squares

$$\begin{array}{ccc} \mathcal{M}_X\left[\begin{smallmatrix}(T_i)_i\\ N\end{smallmatrix}\right] & \longrightarrow & \prod_i \mathrm{Emb}_{/X}(T_i, N) \\ \downarrow & h & \downarrow{\scriptstyle \text{forget}} \\ \mathrm{Emb}\left(\bigsqcup_i T_i, N\right) & \xrightarrow[\text{rest}]{} & \prod_i \mathrm{Emb}(T_i, N) \end{array}.$$

Composition and units are induced from those of $\mathcal{M}\mathrm{fld}_n$ and $\mathcal{M}\mathrm{fld}_{n/X}$.

The full suboperad $\mathcal{D}_X$ of $\mathcal{M}_X$ is the one spanned by colors $(Z, \nu_Z)$ such that $Z \cong (\mathbb{R}^n)^{\sqcup k}$ for some $0 \leqslant k < \infty$. $\mathcal{E}_X$ denotes the full suboperad of $\mathcal{M}_X$ spanned by colors $(Z, \nu_Z)$ so that $Z \cong \mathbb{R}^n$.

Summarizing,

| Without context | $\mathcal{M}\mathrm{fld}_n$ sm category | $\mathcal{D}\mathrm{iscs}_n$ sm category | $\mathcal{E}_n$ operad |
|---|---|---|---|
| With context | $\mathcal{M}_X$ operad | $\mathcal{D}_X$ operad | $\mathcal{E}_X$ operad |

For later reference, let us introduce a bit of notation.

**Notation 3.9.** Motivated by the analogy with manifolds with context just presented, we use the following notation:

$$M_X = \mathrm{Disj}(X), \qquad\qquad D_X = \mathrm{Discs}(X).$$

Hence, there is a chain of inclusions of (set valued) operads $E_X \hookrightarrow D_X \hookrightarrow M_X$ and a canonical commutative diagram of $\mathtt{Top}$-operads

$$\begin{array}{ccccc} E_X & \hookrightarrow & D_X & \hookrightarrow & M_X \\ \downarrow & & \downarrow & & \downarrow \\ \mathcal{E}_X & \hookrightarrow & \mathcal{D}_X & \hookrightarrow & \mathcal{M}_X \end{array}.$$

**All we need about embedding spaces.** We collect the (co)descent properties of embedding spaces with context which are going to play a role in §3.2, §4.1 and §4.2.

**Lemma 3.10.** *Let* $T$ *be a finite disjoint union of discs in* $\mathcal{M}\mathrm{fld}_{n/X}$*, i.e. a color of* $\mathcal{D}_X$*. Then, the* $\mathtt{Top}$*-functor* $\mathrm{Emb}_{/X}(T; \star)\colon \mathcal{M}\mathrm{fld}_{n/X} \to \mathtt{Top}$ *is a Weiss homotopy cosheaf. That is, for any Weiss cover* $(U_i)_{i\in I}$ *of* $Z \in \mathcal{M}\mathrm{fld}_{n/X}$*, the canonical map*

$$\operatorname*{hocolim}_{S\subseteq I} \mathrm{Emb}_{/X}(T; U_S) \longrightarrow \mathrm{Emb}_{/X}(T; Z)$$

*is a weak homotopy equivalence.*

*Proof.* The result holds if we drop the context $X$ by [11, Equation (10)], i.e.

$$\operatorname*{hocolim}_{S\subseteq I} \mathrm{Emb}(T; U_S) \longrightarrow \mathrm{Emb}(T; Z)$$

is a weak homotopy equivalence. Alternatively, one can deduce this result by relating this map with its analogue replacing embeddings by configuration spaces of $|\pi_0(T)|$-points. The claim for embeddings with context $X$ follows by the homotopy pullback definition of $\mathrm{Emb}_{/X}(T; Z)$ since the $\infty$-category of spaces $\mathtt{Top}[\![\mathtt{whe}^{-1}]\!] \simeq \mathcal{S}\mathrm{pc}_\infty$ is an $\infty$-topos. □

*Remark* 3.11. Lemma 3.10 can be also be deduced from Lurie's version of the Seifert-van Kampen theorem [40, Proposition A.3.1], which is a source of broad generalizations for this result, such as the following lemma.

**Lemma 3.12.** *Let* $\mathrm{T}$ *be a color in* $\mathcal{D}_X$*. Then, for any* $Z \in \mathcal{M}\mathrm{fld}_{n/X}$*, there is a canonical weak homotopy equivalence*

$$\operatorname*{hocolim}_{V \in \overline{D}_Z} \mathrm{Emb}_{/X}(T; V) \longrightarrow \mathrm{Emb}_{/X}(T; Z).$$

*Proof.* Note that $\mathrm{Emb}(T, Y)$ is homotopy equivalent to the total space of a $\mathrm{Gl}(n)^{\times t}$-principal bundle over $\mathrm{Conf}_t(Y)$, the configuration space of $t$-points in $Y$ with $t$ being the number of discs in $T$ (see [1, Proposition V.14.4]). Hence, [40, Proposition A.3.1] applied to the subdiagram of open subsets of $Z$ which are diffeomorphic to finite disjoint unions of discs, i.e. $\overline{D}_Z$, gives rise to a weak homotopy equivalence

$$\operatorname*{hocolim}_{V \in \overline{D}_Z} \mathrm{Conf}_t(V) \longrightarrow \mathrm{Conf}_t(Z),$$

that, in turn, yields the corresponding weak homotopy equivalence for $\mathrm{Emb}(T, Z)$. Passing to $\mathrm{Emb}_{/X}(T, Z)$ only requires a homotopy pullback, which respects homotopy colimits, like in Lemma 3.10. □

**Lemma 3.13.** *Let* $\{(D_i, \nu_i)\}_{i \in I}$ *be a finite family of colors in* $\mathcal{E}_X$*. Then, for any pair of* $n$*-manifolds with context* $Z_1, Z_2 \in \mathcal{M}\mathrm{fld}_{n/X}$*, the map induced by considering disjoint unions of embeddings*

$$\coprod_{I = I_1 \amalg I_2} \mathrm{Emb}_{/X}\Big(\bigsqcup_{p \in I_1} D_p, Z_1\Big) \times \mathrm{Emb}_{/X}\Big(\bigsqcup_{q \in I_2} D_q, Z_2\Big) \xrightarrow{\ \sqcup\ } \mathrm{Emb}_{/X}\Big(\bigsqcup_{i \in I} D_i, Z_1 \sqcup Z_2\Big)$$

*is a weak homotopy equivalence.*

*The analogous result holds for any finite collection of* $Z_i$*'s.*

*Proof.* By connectivity and continuity, $\mathrm{Emb}(D_i, Z_1 \sqcup Z_2) \cong \mathrm{Emb}(D_i, Z_1) \amalg \mathrm{Emb}(D_i, Z_2)$. From here, a decomposition similar to that of the statement for $\mathrm{Emb}\left(\bigsqcup_{i \in I} D_i, Z_1 \sqcup Z_2\right)$ follows by tracking each index $i \in I$ and noting that indices sent to the same component of $Z_1 \sqcup Z_2$ should form a joint embedding. To reintroduce the context $X$, substitute this decomposition into the homotopy pullback definition of $\mathrm{Emb}_{/X}\left(\bigsqcup_{i \in I} D_i, Z_1 \sqcup Z_2\right)$, along with the canonical weak homotopy equivalence of spaces

$$\mathbb{1} \overset{h}{\underset{U}{\times}} \Big(\coprod_t V_t\Big) \simeq \coprod_t \Big(\mathbb{1} \overset{h}{\underset{U}{\times}} V_t\Big).$$ □

Additionally, we will need a covering result for the Weiss topology.

**Lemma 3.14.** *[11, Proposition 2.10] Every (smooth) manifold* $Z$ *admits a Weiss cover whose elements are finite disjoint unions of open discs in* $Z$*, i.e. objects of* $D_Z = \mathrm{Discs}(Z)$*, and their finite intersections remain so.*

**Enriched factorization algebras.** Using the topological version $\mathcal{M}_X$ of $M_X$, we propose an enriched analogue of (lax) factorization algebras.

**Definition 3.15.** An *enriched lax factorization algebra* on $X$ is an $\mathcal{M}_X$-algebra satisfying homotopical codescent with respect to Weiss covers (see Definition 2.6). An *enriched factorization algebra* is an enriched lax factorization algebra which satisfies weak monadicity.

The weak monadicity axiom for $\mathcal{M}_X$-algebras is subtle and it will enjoy a distinguished treatment in the main results. Recall that the operad $\mathcal{M}_X$ does not come from a partial monoidal (or pm-)$\mathtt{Top}$-category structure unlike $\mathsf{M}_X = \mathsf{Disj}(X)$ (by Remark A.3, this fact implies $\mathtt{pFact}_X(\mathtt{V}) = \mathtt{LaxMon}[\mathsf{Open}(X), \mathtt{V}]$). However, it is possible to define weak monadicity for $\mathcal{M}_X$-algebras in the following sense: for any finite family of embeddings $\{\nu_i\colon Z_i \hookrightarrow X\}_i$ with disjoint images in X, the map $(\nu_i)_i\colon \bigsqcup_i Z_i \hookrightarrow X$ is also an embedding. Thus, $\big(\bigsqcup_i Z_i, (\nu_i)_i\big)$ is canonically a color of $\mathcal{M}_X$ and it comes with a preferred operation $m[(\nu_i)_i] \in \mathcal{M}_X\left[\begin{smallmatrix}(Z_i)_i\\ \bigsqcup_i Z_i\end{smallmatrix}\right]$. An $\mathcal{M}_X$-algebra F satisfies *weak monadicity* if for any such a family of embeddings, the associated morphism

$$F(m[(\nu_i)_i])\colon F(Z_1)\otimes\cdots\otimes F(Z_k) \longrightarrow F(Z_1\sqcup\cdots\sqcup Z_k)$$

is an equivalence.

The advantage of enriched factorization algebras over ordinary factorization algebras is that the chain of inclusions of $\mathtt{Top}$-operads

$$\mathcal{E}_X \hookrightarrow \mathcal{D}_X \hookrightarrow \mathcal{M}_X$$

carries enough information to construct model categories presenting enriched factorization algebras on X. Our approach to undertake this task is based on the reduction to questions about embedding and configuration spaces.

## 3.2 Factorization homology with context

In this subsection, we discuss a variant of factorization homology (see [1, 2, 36]) which introduces a manifold as context. We also explain by simple categorical arguments the relation between different perspectives about factorization homology.

**On factorization homology.** As we mentioned in §2, factorization homology is a pairing which takes a n-manifold Z and an unoriented $\mathbb{E}_n$-algebra (equiv. an $\mathcal{E}_n$-algebra) A and returns

$$\int_Z A = \mathsf{Emb}(\star, Z) \underset{\mathcal{E}_n}{\overset{\mathbb{L}}{\circ}} A,$$

where $\mathsf{Emb}(\star, Z)$ is considered as a right $\mathcal{E}_n$-module. Taking into account functoriality on the manifold variable, $\int_\star A\colon \mathcal{M}\mathsf{fld}_n \to \mathcal{V}$ is really the derived operadic Lan of $A\colon \mathcal{E}_n \simeq \mathbb{E}_n^{\mathsf{un}} \to \mathcal{V}$ into $\mathcal{M}\mathsf{fld}_n$. Denoting by $j\colon \mathcal{E}_n \hookrightarrow (\mathcal{M}\mathsf{fld}_n)_{\mathsf{opd}}$ the obvious inclusion of operads, we may refer to this observation as

$$\int_\star A \equiv \mathbb{L}j_\sharp A \equiv (\mathcal{M}\mathsf{fld}_n)_{\mathsf{opd}} \underset{\mathcal{E}_n}{\overset{\mathbb{L}}{\circ}} A\,.$$

Notice that since $\mathcal{M}\mathsf{fld}_n$ is a sm-$\mathtt{Top}$-category, $\mathtt{Alg}_{(\mathcal{M}\mathsf{fld}_n)_{\mathsf{opd}}}(\mathtt{V}) \simeq \mathtt{LaxMon}[\mathcal{M}\mathsf{fld}_n, \mathtt{V}]$ (see [32]).

Summarizing, factorization homology computes the derived left adjoint of the restriction functor

$$j^*\colon \mathcal{A}\mathsf{lg}_{(\mathcal{M}\mathsf{fld}_n)_{\mathsf{opd}}}(\mathcal{V}) \longrightarrow \mathcal{A}\mathsf{lg}_{\mathcal{E}_n}(\mathcal{V})$$

and it actually lands in weak monoidal functors inside $\mathtt{LaxMon}[\mathcal{M}\mathsf{fld}_n, \mathtt{V}]$.

On the other side, there is another formula in the literature for factorization homology (e.g. [3, Definition 3.2]), namely

$$\int_Z A = \underset{U\in\mathcal{D}\mathsf{iscs}_{n/Z}}{\mathsf{hocolim}} A^{\otimes \pi_0(U)}\,.$$

The coincidence of the two perspectives is simple provided one has the right ideas at hand. First, note that for the second formula to make sense, we extend the $\mathcal{E}_n$-algebra A to a $\mathcal{D}iscs_n$-algebra $A^{\otimes}$ by setting: $A^{\otimes}(D_1 \sqcup \cdots \sqcup D_k) = A(D_1) \otimes \cdots \otimes A(D_k)$. Then, the (hocolim)-formula can be interpreted as

$$\int_Z A = \operatorname*{hocolim}_{U \in \mathcal{D}iscs_{n/Z}} A^{\otimes \pi_0(U)} = \mathrm{Emb}(\star, Z) \overset{\mathbb{L}}{\underset{\mathcal{D}iscs_n}{\otimes}} A^{\otimes}$$

or alternatively, as

$$\int_\star A = \mathbb{L}\, i_!\, A^{\otimes} \equiv \mathcal{M}fld_n \overset{\mathbb{L}}{\underset{\mathcal{D}iscs_n}{\otimes}} A^{\otimes}.$$

In other words, $\int_\star A$ computes the derived left adjoint of $i^*\colon [\mathcal{M}fld_n, \mathcal{V}] \to [\mathcal{D}iscs_n, \mathcal{V}]$ when evaluated on $A^{\otimes}\colon \mathcal{D}iscs_n \to \mathcal{V}$.

The coincidence of both formulas/perspectives comes from two simple facts:

1. Since $\mathcal{M}fld_n$ and $\mathcal{D}iscs_n$ are sm-Top-categories, we see in §A that the following diagram of ∞-categories commutes

$$\begin{array}{ccc}
[\mathcal{D}iscs_n, \mathcal{V}]_\infty & \xrightarrow{\quad \mathbb{L}\, i_! \quad} & [\mathcal{M}fld_n, \mathcal{V}]_\infty \\
{\scriptstyle \mathrm{fgt}}\uparrow & & \uparrow{\scriptstyle \mathrm{fgt}} \\
\mathrm{LaxMon}\,[\mathcal{D}iscs_n, \mathcal{V}]_\infty & \xrightarrow{\mathbb{L}\, i_!} & \mathrm{LaxMon}\,[\mathcal{M}fld_n, \mathcal{V}]_\infty \\
\wr\uparrow & & \uparrow\wr \\
\mathcal{A}lg_{(\mathcal{D}iscs_n)_{\mathrm{opd}}}(\mathcal{V})_\infty & \xrightarrow{\quad \mathbb{L}\, i_\sharp \quad} & \mathcal{A}lg_{(\mathcal{M}fld_n)_{\mathrm{opd}}}(\mathcal{V})_\infty
\end{array}.$$

   Equivalently, $\mathbb{L}\, i_!\, \overline{B} \simeq \overline{\mathbb{L}\, i_\sharp\, B}$ or simply, $\mathbb{L}\, i_!$ preserves monoidal structures.

2. The map of operads $\mathcal{E}_n \hookrightarrow (\mathcal{D}iscs_n)_{\mathrm{opd}}$ induces an equivalence of sm-Top-categories $\mathrm{Env}_{\otimes}(\mathcal{E}_n) \simeq \mathcal{D}iscs_n$. This is just a simple consequence of continuity and connectivity.

More concretely, the equivalence

$$\mathrm{Emb}(\star, Z) \overset{\mathbb{L}}{\underset{\mathcal{D}iscs_n}{\otimes}} A^{\otimes} \simeq \mathrm{Emb}(\star, Z) \overset{\mathbb{L}}{\underset{\mathcal{E}_n}{\circ}} A$$

between the two ways of computing $\int_Z A$ for $A \in \mathcal{A}lg_{\mathcal{E}_n}(\mathcal{V})$ is the natural equivalence exhibiting commutativity of the following diagram of ∞-categories (when evaluated on $(A, Z)$)

$$\begin{array}{ccccccc}
 & & \mathrm{Fun}(\mathcal{D}iscs_n, \mathcal{V}_\infty) & \xrightarrow{\quad \mathbb{L}\, i_! \quad} & \mathrm{Fun}(\mathcal{M}fld_n, \mathcal{V}_\infty) & & \\
 & & {\scriptstyle \mathrm{fgt}}\uparrow & & \uparrow{\scriptstyle \mathrm{fgt}} & & \\
 & & \mathrm{Fun}^{\mathrm{lax}\otimes}(\mathcal{D}iscs_n, \mathcal{V}_\infty) & \xrightarrow{\mathbb{L}\, i_!} & \mathrm{Fun}^{\mathrm{lax}\otimes}(\mathcal{M}fld_n, \mathcal{V}_\infty) & = & \\
 & & \uparrow & & \uparrow & & \| \\
\mathcal{A}lg_{\mathcal{E}_n}(\mathcal{V})_\infty & \xrightarrow{\sim} & \mathrm{Fun}^{\otimes}(\mathcal{D}iscs_n, \mathcal{V}_\infty) & \xrightarrow[\mathbb{L}\, i_!]{\qquad} & \mathrm{Fun}^{\otimes}(\mathcal{M}fld_n, \mathcal{V}_\infty) & \hookrightarrow & \mathcal{A}lg_{(\mathcal{M}fld_n)_{\mathrm{opd}}}(\mathcal{V})_\infty \\
 & & & \xrightarrow[\mathbb{L}\, i_\sharp]{\qquad\qquad\qquad\qquad} & & &
\end{array},$$

where $\mathrm{Fun}^{\otimes}$ (resp. $\mathrm{Fun}^{\mathrm{lax}\otimes}$) denotes the $\infty$-category of (resp. lax) symmetric monoidal functors between symmetric monoidal $\infty$-categories.

In conclusion, factorization homology yields the derived operadic Lan along a (homotopy) fully-faithful inclusion of operads $\mathcal{B} \hookrightarrow (\mathcal{M}\mathrm{fld}_n)_{\mathrm{opd}}$, where $\mathcal{B}$ is some $\mathtt{Top}$-operad of building blocks, such as discs or finite disjoint unions of discs.

*Remark* 3.16. The (hocolim-)formula computing $\int_Z A$ admits simplifications, usually obtained by (homotopy) cofinality statements; e.g. see [2, §3.2], [36, §7] and [40, §5.5.3].

**Factorization homology with context.** Now, we want to introduce a manifold $X$ as context. To do so, we should replace:

$$\mathcal{M}\mathrm{fld}_n \rightsquigarrow \mathcal{M}_X, \qquad \mathcal{D}\mathrm{iscs}_n \rightsquigarrow \mathcal{D}_X, \qquad \mathcal{E}_n \rightsquigarrow \mathcal{E}_X\,.$$

Then, *factorization homology with context* $X$ has inputs: $(i)$ an $\mathcal{E}_X$-algebra $A$ (resp. a $\mathcal{D}_X$-algebra $B$) and $(ii)$ a manifold $Z \in \mathcal{M}\mathrm{fld}_{n/X}$; and it gives back:

$$\int_Z A = \mathrm{Emb}_{/X}(\star, Z) \overset{\mathbb{L}}{\underset{\mathcal{E}_X}{\circ}} A \qquad \left(\text{resp. } \int_Z B = \mathrm{Emb}_{/X}(\star, Z) \overset{\mathbb{L}}{\underset{\mathcal{D}_X}{\circ}} B\right).$$

In other words, $\int_\star A = \mathbb{L}\, j_\sharp A$ (resp. $\int_\star B = \mathbb{L}\, i_\sharp B$), where

$$\begin{array}{ccc} \mathcal{E}_X & \longrightarrow & \mathcal{D}_X \\ & {\scriptstyle j}\searrow \quad \swarrow{\scriptstyle i} & \\ & \mathcal{M}_X & \end{array} \quad\rightsquigarrow\quad \begin{array}{ccc} \mathcal{A}\mathrm{lg}_{\mathcal{E}_X}(\mathcal{V}) & \longrightarrow & \mathcal{A}\mathrm{lg}_{\mathcal{D}_X}(\mathcal{V}) \\ & {\scriptstyle \mathbb{L} j_\sharp}\searrow \quad \swarrow{\scriptstyle \mathbb{L} i_\sharp} & \\ & \mathcal{A}\mathrm{lg}_{\mathcal{M}_X}(\mathcal{V}) & \end{array}\,.$$

However, it is not clear at all whether or not there is a "simplified" formula to compute factorization homology with context corresponding to derived Lan along $\overline{\mathcal{D}}_X \hookrightarrow \overline{\mathcal{M}}_X$. The obstruction is that the operads $\mathcal{D}_X$ and $\mathcal{M}_X$ do not readily come from symmetric or partial monoidal categories unlike for $\mathcal{D}\mathrm{iscs}_n$ or $\mathcal{M}\mathrm{fld}_n$.

*Remark* 3.17. It is claimed in [4, Lemma 2.16], using different terminology and just as a particular case, that $\mathrm{Emb}_{/X}(\star, Z) \otimes^{\mathbb{L}}_{\overline{\mathcal{D}}_X} \overline{B}$ computes $\int_Z B$ for any $B \in \mathcal{A}\mathrm{lg}_{\mathcal{D}_X}(\mathcal{V})$, or equiv. $\overline{\mathbb{L}\, i_\sharp B} \simeq \mathbb{L}\, i_!\, \overline{B}$. They also assert that this claim follows from an argument similar to that given in the proof of [4, Lemma 2.15]. To be precise, [4, Lemma 2.15] is an $\infty$-categorical version of the fact that ordinary Kan extension preserves monadicity. We feel that such argument is not fully complete as it stands. In fact, we doubt that as stated this result holds due to various reasons. For instance, $(1)$ $\mathcal{M}_X$ and $\mathcal{D}_X$ are operads not coming from partial monoidal categories, and so it is not possible to identify algebras over them by prescribing actions over universal operations; $(2)$ the finality assumption in [4, Lemma 2.16 (5)] has an issue. Using their notation, that assumption says that given objects $Z$ and $Z'$ in $\mathcal{M}_{/X}$, one may construct $Z \otimes Z'$ in $\mathcal{M}$, but one cannot combine their structural maps into $X$ in a canonical way to obtain an object in $\mathcal{M}_{/X}$. Hence, the tensor product functor appearing in $(5)$ exists, but $Z \otimes Z'$ is no longer an object in $\mathcal{M}_{/X}$ so we cannot consider the slice $\infty$-category $(\mathcal{M}_{/X})_{/Z\otimes Z'}$ which is claimed to be equivalent to $\mathcal{M}_{/Z\otimes Z'}$.

We propose a different approach to see that the "simplified" formula computes factorization homology with context. Our task consists of checking if $\mathbb{L}\, i_!$ computes the derived operadic Lan along $i\colon \mathcal{D}_X \hookrightarrow \mathcal{M}_X$ in the sense that the following square commutes

$$\begin{array}{ccccccc} [\overline{\mathcal{D}}_X, \mathcal{V}] & = & \mathcal{A}lg_{\overline{\mathcal{D}}_X}(\mathcal{V}) & \xrightarrow{\mathbb{L}\, i_!} & \mathcal{A}lg_{\overline{\mathcal{M}}_X}(\mathcal{V}) & = & [\overline{\mathcal{M}}_X, \mathcal{V}] \\ & & \uparrow{\scriptstyle \overline{(\star)}} & & \uparrow{\scriptstyle \overline{(\star)}} & & \\ & & \mathcal{A}lg_{\mathcal{D}_X}(\mathcal{V}) & \xrightarrow[\mathbb{L}\, i_\sharp]{} & \mathcal{A}lg_{\mathcal{M}_X}(\mathcal{V}) & & \end{array} . \tag{3.1}$$

Proposition A.6 provides a criterion to see if this is the case.

To apply the cited result, we use the following lemma.

**Lemma 3.18.** *The canonical map of* Top*-valued precosheaves over* $\overline{\mathcal{M}}_X$

$$\mathbb{L}\, i_!\, \mathcal{D}_X \begin{bmatrix} \underline{T} \\ \star \end{bmatrix} \longrightarrow \mathcal{M}_X \begin{bmatrix} i(\underline{T}) \\ \star \end{bmatrix}$$

*is a weak homotopy equivalence for any* $\underline{T} = (T_1, \dots, T_m)$ *finite collection of colors in* $\mathrm{col}(\mathcal{D}_X)$.

*Proof.* The claim can be checked by evaluating the natural transformation on each $Z \in \mathcal{M}_X$, i.e. $\mathbb{L}\, i_!\, \mathcal{D}_X \left[\begin{smallmatrix} \underline{T} \\ \star \end{smallmatrix}\right](Z) \longrightarrow \mathcal{M}_X \left[\begin{smallmatrix} i(\underline{T}) \\ Z \end{smallmatrix}\right]$ is a weak homotopy equivalence for all $Z$. We deduce this equivalence using codescent to reduce the problem to $Z \in \mathcal{D}_X$. Clearly, if $Z \in \mathcal{D}_X$, the claim holds since $i\colon \mathcal{D}_X \hookrightarrow \mathcal{M}_X$ is fully-faithful. For a general $Z \in \mathcal{M}_X$, pick a nice Weiss cover $\{U_i\}_{i \in I}$ of $Z$ by finite disjoint unions of discs (whose existence is ensured by Lemma 3.14) and use the commutative square

$$\begin{array}{ccc} \operatorname*{hocolim}_{S \subseteq I} \mathbb{L}\, i_!\, \mathcal{D}_X \left[\begin{smallmatrix} \underline{T} \\ \star \end{smallmatrix}\right](U_S) & \longrightarrow & \operatorname*{hocolim}_{S \subseteq I} \mathcal{M}_X \left[\begin{smallmatrix} i(\underline{T}) \\ U_S \end{smallmatrix}\right] \\ \downarrow & & \downarrow \\ \mathbb{L}\, i_!\, \mathcal{D}_X \left[\begin{smallmatrix} \underline{T} \\ \star \end{smallmatrix}\right](Z) & \longrightarrow & \mathcal{M}_X \left[\begin{smallmatrix} i(\underline{T}) \\ Z \end{smallmatrix}\right] \end{array} .$$

By 2-out of-3, the lower horizontal map is an equivalence if the rest of the maps forming the square are equivalences. The upper horizontal map is an equivalence since $U_S \in \mathcal{D}_X$ for all $\varnothing \neq S \subseteq I$. Using Lemma 3.10, it is easy to see that the right vertical arrow is an equivalence. The left vertical map is an equivalence also by the cited lemma; just note

$$\begin{aligned} \operatorname*{hocolim}_{S \subseteq I} \mathbb{L}\, i_!\, \mathcal{D}_X \begin{bmatrix} \underline{T} \\ \star \end{bmatrix}(U_S) &= \operatorname*{hocolim}_{S \subseteq I} \left( \overline{\mathcal{M}}_X \begin{bmatrix} i(\bullet) \\ U_S \end{bmatrix} \overset{\mathbb{L}}{\underset{\overline{\mathcal{D}}_X}{\otimes}} \mathcal{D}_X \begin{bmatrix} \underline{T} \\ \bullet \end{bmatrix} \right) \\ &\simeq \left( \operatorname*{hocolim}_{S \subseteq I} \overline{\mathcal{M}}_X \begin{bmatrix} i(\bullet) \\ U_S \end{bmatrix} \right) \overset{\mathbb{L}}{\underset{\overline{\mathcal{D}}_X}{\otimes}} \mathcal{D}_X \begin{bmatrix} \underline{T} \\ \bullet \end{bmatrix} \\ &\simeq \overline{\mathcal{M}}_X \begin{bmatrix} i(\bullet) \\ Z \end{bmatrix} \overset{\mathbb{L}}{\underset{\overline{\mathcal{D}}_X}{\otimes}} \mathcal{D}_X \begin{bmatrix} \underline{T} \\ \bullet \end{bmatrix} \\ &= \mathbb{L}\, i_!\, \mathcal{D}_X \begin{bmatrix} \underline{T} \\ \star \end{bmatrix}(Z). \end{aligned}$$

□

Consequently, we obtain:

**Theorem 3.19.** *The square (3.1) commutes up to natural equivalence, i.e. the mate natural transformation* $\mathbb{L}\,\mathrm{i}_!\,\overline{\mathrm{B}} \to \overline{\mathbb{L}\,\mathrm{i}_\sharp\,\mathrm{B}}$ *is an equivalence. In other words,*

$$\mathrm{Emb}_{/X}(\star, Z) \overset{\mathbb{L}}{\underset{\overline{\mathcal{D}}_X}{\otimes}} \overline{\mathrm{B}} \simeq \int_Z \mathrm{B} \simeq \mathrm{Emb}_{/X}(\star, Z) \overset{\mathbb{L}}{\underset{\mathcal{D}_X}{\circ}} \mathrm{B}$$

*for any* $\mathrm{B} \in \mathcal{A}\mathrm{lg}_{\mathcal{D}_X}(\mathcal{V})$ *and* $Z \in \mathcal{M}\mathrm{fld}_{n/X}$.

*Remarks* 3.20.

1. We see in Proposition A.5 that derived operadic Lan along a strong monoidal functor of partial monoidal categories coincides with derived ordinary/categorical Lan under our assumptions on $\mathcal{V}$. This implies without further work that the discrete version of the discussion above holds; meaning: there is a natural equivalence $\mathbb{L}\iota_!\overline{\mathrm{B}} \simeq \overline{\mathbb{L}\iota_\sharp \mathrm{B}}$, where $\mathrm{B} \in \mathcal{A}\mathrm{lg}_{\mathrm{D}_X}(\mathcal{V}) = \mathrm{LaxMon}\,[\mathrm{D}_X, \mathcal{V}]$ and $\iota\colon \mathrm{D}_X \hookrightarrow \mathrm{M}_X$ is the obvious inclusion of partial monoidal categories.

2. In [36, Proposition 3.9], it is computed in general the derived operadic left Kan extension using enveloping sm-categories. Our result simplifies this presentation since we use the underlying category instead of sm-envelopes. More precisely, without context, Horel's formula yields

$$\int_Z \mathrm{A} = \mathrm{Emb}(\star, Z) \overset{\mathbb{L}}{\underset{\mathcal{E}_n}{\circ}} \mathrm{A} \simeq \mathrm{Emb}(\star, Z) \overset{\mathbb{L}}{\underset{\mathrm{Env}_\otimes(\mathcal{E}_n)}{\otimes}} \mathrm{A}^\otimes \simeq \mathrm{Emb}(\star, Z) \overset{\mathbb{L}}{\underset{\mathcal{D}\mathrm{iscs}_n}{\otimes}} \mathrm{A}^\otimes$$

   since $\mathrm{Env}_\otimes(\mathcal{E}_n) \simeq \mathcal{D}\mathrm{iscs}_n$ as sm-Top-categories. However, introducing a manifold $X$ as context, even for $\mathcal{D}_X$-algebras, it yields

$$\int_Z \mathrm{B} = \mathrm{Emb}_{/X}(\star, Z) \overset{\mathbb{L}}{\underset{\mathcal{D}_X}{\circ}} \mathrm{B} \simeq \mathrm{Emb}_{/X}(\star, Z) \overset{\mathbb{L}}{\underset{\mathrm{Env}_\otimes(\mathcal{D}_X)}{\otimes}} \mathrm{B}^\otimes,$$

   which is no longer a simplification.

3. In the proof of Lemma 3.18 we applied a codescent property of embedding spaces (Lemma 3.10), but one can use any other suitable codescent property such as the one given in Lemma 3.12.

# 4 A bunch of model structures

In this section, we construct the model structures we were looking for in this paper. Let us start with a brief reminder of [15, §4], because we will use the results obtained there extensively in our current constructions.

**Extension model structure.** Local-to-global conditions for sheaves are treated by means of left Bousfield localizations in homotopy theory [25, 41]. By dual analogy, colocality principles for cosheaves may be treated with right Bousfield localizations, although it is not clear what right Bousfield localization could correspond to those colocality principles. Furthermore, classical methods to find right Bousfield localizations are rarely applicable. In [15, §4], we give a general procedure to construct certain right Bousfield localizations for model categories of operadic algebras which certainly capture colocality principles. Because of their applicability in this work, we include a brief discussion of the results obtained in the cited reference.

Let $\iota\colon \mathcal{B} \to \mathcal{N}$ be a map of $\mathcal{V}$-operads. This simple data already determines a local-to-global principle: $\mathcal{N}$-algebras that are completely characterized by their restriction to $\mathcal{B}$. A canonical choice of reconstruction process from the restricted $\mathcal{B}$-algebra to an $\mathcal{N}$-algebra is via operadic left Kan extension. More concretely, $\iota$ induces an adjoint pair $\iota_\sharp\colon \mathtt{Alg}_{\mathcal{B}} \rightleftarrows \mathtt{Alg}_{\mathcal{N}}\colon \iota^*$. An $\mathcal{N}$-algebra $A$ satisfies the above local-to-global principle if $\iota_\sharp\iota^* A \to A$ is an isomorphism, because $\iota^* A$ is the restricted $\mathcal{B}$-algebra and $\iota_\sharp$ the operadic left Kan extension.

In order to introduce homotopy theory into the picture, the simplest approach is to enhance the adjunction $\iota_\sharp \dashv \iota^*$ to a Quillen pair. Fixing the projective model structure for operadic algebras on both sides, one achieves such an enhancement. Then, the derived functor $\mathbb{L}\iota_\sharp$ gives a homotopy meaningful reconstruction from $\mathcal{B}$-algebras to $\mathcal{N}$-algebras. To refer to $\mathcal{N}$-algebras which satisfy the corresponding homotopical local-to-global principle, we introduce the following notion.

**Definition 4.1.** An $\mathcal{N}$-algebra $A$ such that $\mathbb{L}\iota_\sharp\iota^* A \to A$ is an equivalence is said to be $\mathcal{B}$-colocal.

*Remark* 4.2. The above discussion is quite trivial in nature and that is why is so surprising that it has produced so deep mathematics. The formal idea is the germ of chiral homology [40] and recent local-to-global constructions in algebraic quantum field theories [10]. In a different direction, manifold calculus [11] also builds from a similar idea.

In [15, §4], it is discussed in detail how one can produce a right Bousfield localization of the projective model structure on $\mathtt{Alg}_{\mathcal{N}}$ whose colocal objects are $\mathcal{B}$-colocal algebras, i.e. algebras that can be reconstructed from its restriction to $\mathcal{B}$. However, this program requires some assumptions to work. We collect them together with their raison d'être.

**Hypothesis 4.3.** *We assume the following conditions.*

- *$\mathcal{V}$ is cofibrantly generated and $\mathcal{B}$, $\mathcal{N}$ are admissible operads, i.e. $\mathtt{Alg}_{\mathcal{B}}$, $\mathtt{Alg}_{\mathcal{N}}$ admit the projective model structure. We impose this condition to ensure there is a Quillen pair $\iota_\sharp\colon \mathcal{A}\mathrm{lg}_{\mathcal{B}} \rightleftarrows \mathcal{A}\mathrm{lg}_{\mathcal{N}} : \iota^*$.*
- *$\mathbb{L}\iota_\sharp\colon \mathrm{Ho}\,\mathcal{A}\mathrm{lg}_{\mathcal{B}} \to \mathrm{Ho}\,\mathcal{A}\mathrm{lg}_{\mathcal{N}}$ is fully-faithful. We impose this condition to ensure that by restricting $\mathbb{L}\iota_\sharp\iota^* A$ and $A$ back to $\mathcal{B}$ we get the same result.*

*Remark* 4.4. Notice that the first point above is already implied by the restrictive Hypothesis 1.2 that we consider in this document.

**Theorem 4.5.** *[15, Theorem 4.9 and Proposition 4.10] Assuming Hypothesis 4.3, the projective model structure $\mathcal{A}\mathrm{lg}_{\mathcal{N}}(\mathcal{V})$ admits a right Bousfield localization, called the extension model structure and denoted $\mathcal{A}\mathrm{lg}_{\mathcal{N}}(\mathcal{V})_{\mathcal{B}}$, which enjoys the following properties:*

- *The weak equivalences are those morphisms of $\mathcal{N}$-algebras which are equivalences when restricted to $\mathcal{B}$.*

- *A cofibrant object is a proj-cofibrant algebra which is also $\mathcal{B}$-colocal.*

- *The extension model structure is cofibrantly generated.*

*Furthermore, the Quillen pair $\iota_\sharp \dashv \iota^*$ descends to a Quillen equivalence between the projective model structure $\mathcal{A}\mathrm{lg}_{\mathcal{B}}(\mathcal{V})$ and the extension model structure $\mathcal{A}\mathrm{lg}_{\mathcal{N}}(\mathcal{V})_{\mathcal{B}}$.*

This general construction is applicable to our current search of a model structure for locally constant factorization algebras. The main point is that there are two chains of fully-faithful inclusions of operads

$$\begin{array}{cc} \mathrm{E}_X \hookrightarrow \mathrm{D}_X \hookrightarrow \mathrm{M}_X, & \mathcal{E}_X \hookrightarrow \mathcal{D}_X \hookrightarrow \mathcal{M}_X \\ \text{(discrete version)} & \text{(enriched version)} \end{array}$$

suitable to feed the extension model machine.

## 4.1 Enriched factorization model structure

Our goal in this subsection is constructing model structures on $\mathtt{Alg}_{\mathcal{M}_X}(\mathtt{V})$ presenting the homotopy theory of enriched (lax) factorization algebras.

**Homotopical codescent and enriched lax factorization algebras.** Let us fix the inclusion of operads $i\colon \mathcal{D}_X \hookrightarrow \mathcal{M}_X$ along this subsection.

The main idea here is to describe a relation between the homotopical codescent condition for Weiss covers and $\mathcal{D}_X$-colocality. Our main tool will be the computation of the derived operadic left Kan extension (factorization homology with context)

$$\mathbb{L}\, i_\sharp\colon \mathcal{A}\mathrm{lg}_{\mathcal{D}_X}(\mathcal{V}) \to \mathcal{A}\mathrm{lg}_{\mathcal{M}_X}(\mathcal{V})$$

given in Theorem 3.19. The extension model machine will produce the claimed model structure for enriched lax factorization algebras.

The main deduction from Theorem 3.19 in this subsection is the recognition principle described in the following result.

**Proposition 4.6.** *For any $B \in \mathcal{A}\mathrm{lg}_{\mathcal{D}_X}(\mathcal{V})$, its extension $\mathbb{L}\, i_\sharp B$ is an enriched lax factorization algebra, i.e. it satisfies homotopical codescent with respect to Weiss covers.*

*Moreover, the following conditions are equivalent for an $\mathcal{M}_X$-algebra $F$:*

- *$F$ is an enriched lax factorization algebra (see Definition 3.15).*

- *$F$ is $\mathcal{D}_X$-colocal (see Definition 4.1).*

*Proof.* Let $(U_t)_{t\in I}$ be a Weiss cover of $U \subseteq X$. The first statement asserts that the canonical map

$$\operatorname*{hocolim}_{S\subseteq I} \mathbb{L}\, i_\sharp B(U_S) \longrightarrow \mathbb{L}\, i_\sharp B(U)$$

is an equivalence. We deduce this fact using 2-out of-3 and the commutative diagram

$$\begin{array}{ccc}
\operatorname*{hocolim}_{S\subseteq I} \mathbb{L}\, i_\sharp\, B(U_S) & \xleftarrow{\sim} & \operatorname*{hocolim}_{S\subseteq I} \left( \mathrm{Emb}_{/X}(\star, U_S) \overset{\mathbb{L}}{\underset{\overline{\mathcal{D}}_X}{\otimes}} \overline{B} \right) \\
\Big\downarrow & & \Big\downarrow \wr \\
 & & \left( \operatorname*{hocolim}_{S\subseteq I} \mathrm{Emb}_{/X}(\star, U_S) \right) \overset{\mathbb{L}}{\underset{\overline{\mathcal{D}}_X}{\otimes}} \overline{B} \ , \\
 & & \Big\downarrow \wr \\
\mathbb{L}\, i_\sharp\, B(U) & \xleftarrow{\sim} & \mathrm{Emb}_{/X}(\star, U) \overset{\mathbb{L}}{\underset{\overline{\mathcal{D}}_X}{\otimes}} \overline{B}
\end{array}$$

where the horizontal arrows are equivalences by Theorem 3.19 and the lower right vertical map is so by Lemma 3.10.

Now, in order to prove the equivalence of the rest of the statements, fix an $\mathcal{M}_X$-algebra $F$. Observe that $\mathcal{D}_X$-colocality means that $\mathfrak{Q}\, F \to F$ is an equivalence, where $\mathfrak{Q}\, F = i_\sharp\, Q\, i^*\, F$ and $Q$ denotes proj-cofibrant replacement. Moreover, $\mathfrak{Q}\, F$ satisfies homotopical codescent with respect to Weiss covers by the former argument and hence $F$ as well.

Conversely, we deduce $\mathcal{D}_X$-colocality of $F$ by evaluation on a manifold $Z$ in $\mathrm{col}(\mathcal{M}_X)$, i.e. showing that $\mathfrak{Q}\, F(Z) \to F(Z)$ is an equivalence. Lemma 3.14 provides a particular Weiss cover $(Z_i)_{i\in I}$ of $Z$ which we employ to obtain the above equivalence using the commutative square

$$\begin{array}{ccc}
\operatorname*{hocolim}_{S\subseteq I} \mathfrak{Q}\, F(Z_S) & \xrightarrow{\sim} & \operatorname*{hocolim}_{S\subseteq I} F(Z_S) \\
\wr \Big\downarrow & & \Big\downarrow \wr \\
\mathfrak{Q}\, F(Z) & \longrightarrow & F(Z)
\end{array} ,$$

since both $\mathcal{M}_X$-algebras coincide on colors in $\mathrm{col}(\mathcal{D}_X)$ and satisfy homotopical codescent with respect to Weiss covers. □

The previous characterization of enriched lax factorization algebras is the keystone to the construction of the enriched lax factorization model structure.

**Theorem 4.7.** *The category* $\mathtt{Alg}_{\mathcal{M}_X}(\mathtt{V})$ *supports the* enriched lax factorization model structure $\mathrm{lax}\,\mathcal{F}\mathrm{act}_X^{\mathrm{enr}}(\mathcal{V})$ *characterized by:*

- *The class of weak equivalences consists on those maps that are equivalences when evaluated on colors of* $\mathcal{D}_X$.
- *The class of fibrations is that of proj-fibrations.*
- *Cofibrant objects are proj-cofibrant enriched lax factorization algebras.*

*Proof.* It is an immediate application of Theorem 4.5 together with the recognition principle given in Proposition 4.6. □

**Weak monadicity and enriched factorization algebras.** Now we work instead with the operad inclusion $j\colon \mathcal{E}_X \hookrightarrow \mathcal{M}_X$.

As we did with enriched lax factorization algebras, we want to further relate weak monadicity with $\mathcal{E}_X$-colocality.

**Proposition 4.8.** *For any $A \in \mathcal{A}lg_{\mathcal{E}_X}(\mathcal{V})$, its extension $\mathbb{L}\, j_\sharp\, A$ is an enriched factorization algebra, i.e. it satisfies Weiss homotopical codescent and weak monadicity.*

*Moreover, an $\mathcal{M}_X$-algebra $F$ such that $F(Z) \in \mathcal{V}$ is cofibrant or flat for each $Z$ is an enriched factorization algebra if and only if it is $\mathcal{E}_X$-colocal.*

*Proof.* Noting that $j$ factors through the morphism $i\colon \mathcal{D}_X \hookrightarrow \mathcal{M}_X$, it is clear that $\mathbb{L}\, j_\sharp\, A$ satisfies homotopical codescent with respect to Weiss covers by Proposition 4.6.

It remains to check weak monadicity, that is

$$\mathbb{L}\, j_\sharp\, A\left(m[(\nu_i)_i]\right)\colon \bigotimes_i \mathbb{L}\, j_\sharp\, A(Z_i) \longrightarrow \mathbb{L}\, j_\sharp\, A\left(\bigsqcup_i Z_i\right)$$

is an equivalence for any finite family $\{(Z_i, \nu_i)\}_i \subset \mathrm{col}(\mathcal{M}_X)$ of disjoint manifolds with context. Let us argue the case $\{Z_i\}_i = \{Z, Y\}$, since the general case is analogous. Unwrapping the definitions, one sees that this operation is induced by taking disjoint union of embeddings

$$\bigsqcup_{J = J_Z \amalg J_Y} \mathcal{M}_X \begin{bmatrix} \underline{T}_Z \\ Z \end{bmatrix} \otimes \mathcal{M}_X \begin{bmatrix} \underline{T}_Y \\ Y \end{bmatrix} \xrightarrow{\sqcup} \mathcal{M}_X \begin{bmatrix} \underline{T} \\ Z \sqcup Y \end{bmatrix} \tag{4.1}$$

where $\underline{T} = (T_j)_{j \in J}$ is a sequence of objects $T_j \in \mathrm{col}(\mathcal{E}_X)$, i.e. each one is a disc embedded in $X$ ($Z$ and $Y$ being disjoint in $X$ ensure that this map exists). Since the elements that conform the collection $\underline{T}$ are all connected, (4.1) is an equivalence as a consequence of Lemma 3.13. We should pursue the constructions involved to get $\mathbb{L}\, j_\sharp\, A\left(m[\nu_Z, \nu_Y]\right)$ from (4.1) in order to deduce that weak monadicity holds. Recall that Hypothesis 1.2 guarantees that proj-cofibrant $\mathcal{M}_X$-algebras are fgt-cofibrant. A somewhat lengthy but straightforward analysis of the homotopy colimit formula that defines the derived operadic left Kan extension $\mathbb{L}\, j_\sharp\, A$ shows that these constructions are just: taking homotopy colimits over finite groupoids, tensor products with cofibrant objects and homotopy geometric realization. Hence, $\mathbb{L}\, j_\sharp\, A\left(m[\nu_Z, \nu_Y]\right)$ is an equivalence.

Now, assume that $F$ is a colorwise flat $\mathcal{M}_X$-algebra. If $F$ is $\mathcal{E}_X$-colocal, i.e. $\epsilon\colon \widehat{\mathfrak{Q}}F \to F$ is a colorwise equivalence for $\widehat{\mathfrak{Q}} = j_\sharp\, Q\, j^*$, we use the commutative square

$$\begin{array}{ccc} \bigotimes_i \widehat{\mathfrak{Q}}F(Z_i) & \xrightarrow{\widehat{\mathfrak{Q}}F(m)} & \widehat{\mathfrak{Q}}F\left(\bigsqcup_i Z_i\right) \\ {\scriptstyle \otimes_i \epsilon}\downarrow & & \downarrow{\scriptstyle \epsilon} \\ \bigotimes_i F(Z_i) & \xrightarrow[F(m)]{} & F\left(\bigsqcup_i Z_i\right) \end{array}$$

to see that $F$ satisfies weak monadicity. The top arrow is an equivalence by the above argument, $\epsilon$ is an equivalence by definition of $\mathcal{E}_X$-colocality and $\otimes_i \epsilon$ is an equivalence since $F$ is colorwise flat and $\widehat{\mathfrak{Q}}F$ is proj-cofibrant.

Conversely, we are going to prove that $\epsilon\colon \widehat{\mathfrak{Q}}F \to F$ is an equivalence by evaluation on each $Z \in \mathrm{col}(\mathcal{M}_X)$ when $F$ is an enriched factorization algebra and colorwise flat. The argument is an adaptation of the one given in Proposition 4.6. By definition, if $Z$ belongs to $\mathrm{col}(\mathcal{E}_X)$, i.e. it is an embedded disc into $X$, the claim holds. In general, consider a Weiss cover $(Z_i)_{i\in I}$ of $Z$ by finite disjoint unions of discs whose intersections remain so, Lemma 3.14. Combining homotopical codescent with weak monadicity one can further decompose each finite intersection $Z_S$ of the Weiss cover into its constituent discs $Z_S = \bigsqcup_{j\in J_S} Z_j$ and provide an equivalence

$$\underset{S\subseteq I}{\mathrm{hocolim}} \bigotimes_{j\in J_S} F(Z_j) \xrightarrow{\sim} F(Z). \tag{4.2}$$

This is also valid for $\widehat{\mathfrak{Q}}F$. Therefore, the commutative square

$$\begin{array}{ccc} \underset{S\subseteq I}{\mathrm{hocolim}} \bigotimes_{j\in J_S} \widehat{\mathfrak{Q}}F(Z_j) & \xrightarrow{\sim} & \underset{S\subseteq I}{\mathrm{hocolim}} \bigotimes_{j\in J_S} F(Z_j) \\ \wr\downarrow & & \downarrow\wr \\ \widehat{\mathfrak{Q}}F(Z) & \longrightarrow & F(Z) \end{array},$$

allows us to conclude the claim, where the top map is an equivalence due to $F$ being colorwise flat and $\widehat{\mathfrak{Q}}F$ proj-cofibrant. □

*Remark* 4.9. The equivalence (4.2) is strongly related to homotopy codescent with respect to factorizing covers.

Analogously to Theorem 4.7, the former characterization permits the construction of the so called enriched factorization model structure.

**Theorem 4.10.** *The category* $\mathtt{Alg}_{\mathcal{M}_X}(\mathtt{V})$ *supports the* enriched factorization model structure $\mathcal{F}\mathrm{act}^{\mathrm{enr}}_X(\mathcal{V})$ *characterized by:*

- *The class of weak equivalences consists on those maps that are equivalences when evaluated on colors of* $\mathcal{E}_X$.
- *The class of fibrations is that of proj-fibrations.*
- *Cofibrant objects are proj-cofibrant enriched factorization algebras.*

*Proof.* Immediate application of Theorem 4.5 together with Proposition 4.8. □

## 4.2 LC-Factorization model structure

The next section is devoted to the discretization of the enriched factorization model structure constructed in Theorem 4.10. This goal will be achieved after two steps. First, using the equivalence [2, §2.4] (or equivalently [40, Lemma 5.4.5.11]), we construct a left Bousfield localization which encodes local constancy; and then, a right Bousfield localization will incorporate the local-to-global properties.

In this section, we perform left Bousfield localizations, see [6, 7, 16], and so our fixed homotopy cosmos $\mathcal{V}$ should be nice enough to ensure this process. Since left Bousfield localization at a set of maps produces a model structure under left properness and only a left semimodel structure without such property, we will distinguish two cases in the sequel.

The best possible scenario will be encoded by the following (see Theorem B.2):

**Hypothesis 4.11.** *$\mathcal{V}$ is a left proper combinatorial (resp. cellular) model category such that the projective model structure on $\mathtt{pFact}_X(\mathtt{V})$ is left proper.*

In particular, Hypothesis 4.11 holds when $\mathcal{V} = \mathcal{C}h(R)$ is the projective model structure on chain complexes for $R$ a commutative ring such that $\mathbb{Q} \subseteq R$; e.g. for $R$ a field with $\mathrm{char}(R) = 0$.

For more general homotopy cosmoi, in order to ensure that left Bousfield localization at a set of maps exists, we will have to consider:

**Hypothesis 4.12.** *$\mathcal{V}$ is a tractable[i] combinatorial (resp. cellular) model structure.*

Note that this requirement is already part of Hypothesis 1.2. Despite this fact, we will mention it in the statements involving left Bousfield localizations, in order to clarify its role.

**Local constancy on $D_X$-algebras.** First, we will explain that a convenient choice of unary operations in $D_X$ (or $E_X$) yields a left Bousfield localization of $\mathcal{A}lg_{D_X}(\mathcal{V})$ (resp. $\mathcal{A}lg_{E_X}(\mathcal{V})$) which is Quillen equivalent to the projective model structure on $\mathcal{D}_X$-algebras (resp. $\mathcal{E}_X$). The Quillen equivalence is induced by

$$D_X \longrightarrow \mathcal{D}_X \quad (\text{resp. } E_X \longrightarrow \mathcal{E}_X).$$

We will focus on $D_X$, because the other case is similar.

Recall from [2, Proposition 2.19] or [37, §5.2] that $\ell\colon D_X \to \mathcal{D}_X$ exhibits $\mathcal{D}_X$ as the $\infty$-localization of $D_X$ at the set of isotopy equivalences $\mathcal{J}_X$ in $D_X$ (as an $\infty$-operad).

**Definition 4.13.** A $D_X$-algebra $B$ is *locally constant* if it sends $\mathcal{J}_X$ to equivalences.

The idea is to represent the set of unary operations $\mathcal{J}_X$ as a set of morphisms $\mathcal{S}_X$ in $\mathtt{Alg}_{D_X}$ to perform a left Bousfield localization whose local objects are locally constant $D_X$-algebras.

In order to represent unary operations of $\mathcal{J}_X$ in $\mathtt{Alg}_{D_X}(\mathtt{V})$, we use the composite

$$\begin{array}{ccccc} \overline{D}_X^{op} & \xrightarrow{\text{Yoneda}} & \mathtt{Alg}_{\overline{D}_X}(\mathtt{V}) & \xrightarrow{\text{extension}} & \mathtt{Alg}_{D_X}(\mathtt{V}), \\ u & \longmapsto & y(u) & \longmapsto & \mathrm{ext}\, y(u) \end{array}$$

where $\overline{\mathcal{O}}$ denotes the underlying category of $\mathcal{O}$ and extension is the left adjoint of the canonical restriction functor along $\overline{\mathcal{O}} \to \mathcal{O}$. It is possible to show that there is a natural weak homotopy equivalence

$$\mathbb{R}\mathrm{Map}_{D_X}(\mathbb{L}\mathrm{ext}\, y(u), A) \simeq \mathbb{R}\mathrm{Map}_{\mathcal{V}}\left(\mathbb{I}, A(u)\right),$$

where $\mathbb{R}\mathrm{Map}$ denotes the homotopy mapping space (see the proof of Proposition 4.16). This natural identification almost serves to identify locally constant $D_X$-algebras as local objects for a set of maps. What is missing is that the homotopy mapping space $\mathbb{R}\mathrm{Map}_{\mathcal{V}}(\mathbb{I}, \star)$ may not detect equivalences in $\mathcal{V}$. One may fix this issue by the introduction of a sort of homotopy generators, whose existence is known in our cases of interest (see Proposition 4.17).

**Hypothesis 4.14.** *There is a set of (cofibrant) objects $\mathcal{G}$ in $\mathcal{V}$ which jointly detect equivalences through $\mathbb{R}\mathrm{Map}_{\mathcal{V}}$. That is, a map $f$ in $\mathcal{V}$ is an equivalence if $\mathbb{R}\mathrm{Map}_{\mathcal{V}}(x, f)$ is an equivalence for every $x \in \mathcal{G}$.*

[i] In both cases, tractability means that the generating sets of (trivial) cofibrations have cofibrant domains.

Assuming Hypothesis 4.14, one can make the following definition.

**Definition 4.15.** The set of maps $\mathcal{S}_X$ in $\mathcal{A}\mathrm{lg}_{D_X}(\mathcal{V})$ that *represents* $\mathcal{J}_X$ is

$$\mathbb{L}\mathrm{ext}(\mathcal{G} \overset{\mathbb{L}}{\odot} y(\mathcal{J}_X)) = \left\{ \begin{matrix} \mathbb{L}\mathrm{ext}(x \overset{\mathbb{L}}{\odot} y(j)) \colon \mathbb{L}\mathrm{ext}(x \overset{\mathbb{L}}{\odot} y(u)) \to \mathbb{L}\mathrm{ext}(x \overset{\mathbb{L}}{\odot} y(v)) \\ \text{for } x \in \mathcal{G} \text{ and } j \in \mathcal{J}_X \end{matrix} \right\},$$

where $\odot$ means $\mathcal{V}$-tensoring in $\overline{D}_X$-algebras, ext is the left adjoint of the restriction functor $\mathcal{A}\mathrm{lg}_{D_X}(\mathcal{V}) \to \mathcal{A}\mathrm{lg}_{\overline{D}_X}(\mathcal{V})$ and $y \colon \overline{D}_X^{\mathrm{op}} \to \mathcal{A}\mathrm{lg}_{\overline{D}_X}(\mathcal{V})$ denotes the enriched Yoneda functor.

**Proposition 4.16.** *Assume that Hypothesis 4.14 holds. Then,* $A \in \mathcal{A}\mathrm{lg}_{D_X}(\mathcal{V})$ *is* $\mathcal{S}_X$*-local if and only if it is locally constant.*

*Proof.* There is a chain of natural equivalences for a proj-fibrant $A \in \mathtt{Alg}_{D_X}(\mathtt{V})$

$$\begin{aligned} \mathbb{R}\mathrm{Map}_{D_X}\left(\mathbb{L}\mathrm{ext}(x \overset{\mathbb{L}}{\odot} y(u)), A\right) &\overset{(i)}{\simeq} \mathbb{R}\mathrm{Map}_{\overline{D}_X}\left(x \overset{\mathbb{L}}{\odot} y(u), \overline{A}\right) \\ &\overset{(ii)}{\simeq} \mathrm{Hom}_{\overline{D}_X}\left(x \odot Q\, y(u), \overline{A}^{\Delta^\bullet}\right) \\ &\overset{(iii)}{\cong} \mathrm{Hom}_{\mathtt{V}}\left(x, \underline{\mathrm{Hom}}_{\overline{D}_X}\left(Q\, y(u), \overline{A}^{\Delta^\bullet}\right)\right) \\ &\overset{(iv)}{\simeq} \mathrm{Hom}_{\mathtt{V}}\left(x, \overline{A}^{\Delta^\bullet}(u)\right) \\ &\overset{(v)}{\simeq} \mathbb{R}\mathrm{Map}_{\mathcal{V}}\left(x, A(u)\right), \end{aligned}$$

where we have applied: $(i)$ ext $\dashv \overline{(\star)}$ is a Quillen pair; $(ii)$ use that $x$ is cofibrant and that the homotopy mapping space can be computed by a choice of simplicial frame $\overline{A} \to \overline{A}^{\Delta^\bullet}$ (see [35, §17]); $(iii)$ follows from the adjunction given by $\mathcal{V}$-tensor and $\mathcal{V}$-enrichment on $\mathtt{Alg}_{\overline{D}_X}(\mathtt{V})$, denoted by $\underline{\mathrm{Hom}}_{\overline{D}_X}$; $(iv)$ is a consequence of enriched Yoneda lemma; $(v)$ comes from the fact that evaluation at $u$ sends the simplicial frame of $\overline{A}$ to a simplicial frame of $A(u)$, roughly speaking $\overline{A}^{\Delta^\bullet}(u) \simeq A(u)^{\Delta^\bullet}$, together with the assumption that $x$ is cofibrant. Combining Hypothesis 4.14 with this chain of natural equivalences, the claim holds. ☐

Before stating the main result of this subsection, let us point out that Hypothesis 4.14 already holds in the cases of interest.

**Proposition 4.17.** *If* $\mathcal{V}$ *is combinatorial (resp. left proper or tractable cofibrantly generated) model category, there is a set of cofibrant objects* $\mathcal{G}$ *in* $\mathcal{V}$ *which jointly detect equivalences through* $\mathbb{R}\mathrm{Map}_{\mathcal{V}}$*.*

*Proof.* In the combinatorial case, it is well known that it suffices to take the collection of cofibrant $\kappa$-small objects for a sufficiently large regular cardinal $\kappa$. In the cofibrantly generated cases, the idea is to take (replacements) of domains and codomains of the generating cofibrations (see [24, Proposition A.5]). ☐

We are finally ready for our main result in this subsection.

**Theorem 4.18.** *Assume that Hypothesis 4.12 holds. Then, the projective model structure* $\mathcal{A}\mathsf{lg}_{\mathsf{D}_\mathsf{X}}(\mathcal{V})$ *admits a left Bousfield localization at* $\mathcal{S}_\mathsf{X}$ *among left semimodel categories whose local objects are the locally constant* $\mathsf{D}_\mathsf{X}$*-algebras. We denote such a semimodel category* $\mathcal{A}\mathsf{lg}^{\mathsf{lc}}_{\mathsf{D}_\mathsf{X}}(\mathcal{V})$*. Moreover,*

$$\ell_\sharp\colon \mathcal{A}\mathsf{lg}^{\mathsf{lc}}_{\mathsf{D}_\mathsf{X}}(\mathcal{V}) \rightleftarrows \mathcal{A}\mathsf{lg}_{\mathcal{D}_\mathsf{X}}(\mathcal{V})\colon \ell^*$$

*establishes a Quillen equivalence.*

*Proof.* The existence of the left Bousfield localization follows from [7, Theorem A] and [16, Theorem B] under Hypothesis 4.12. For a full model structure use [6, Theorem 4.7] or [35, Theorem 4.1.1]. We have applied that the projective model structure $\mathcal{A}\mathsf{lg}_{\mathsf{D}_\mathsf{X}}$ is tractable (because $\mathcal{V}$ is tractable and the projective model structure is obtained by transfer) and additionally in the cellular case that cofibrations with cofibrant domains are effective monomorphisms (since those cofibrations are colorwise cofibrations). The characterization of local objects is given in Proposition 4.16.

The Quillen adjunction is a Quillen equivalence because $\mathsf{D}_\mathsf{X} \to \mathcal{D}_\mathsf{X}$ exhibits $\mathcal{D}_\mathsf{X}$ as the $\infty$-localization of $\mathsf{D}_\mathsf{X}$ at the set of isotopy equivalences $\mathcal{J}_\mathsf{X}$ due to [2, Proposition 2.19].

It remains to check that the adjunction $\ell_\sharp \dashv \ell^*$ descends to a Quillen equivalence between $\mathcal{A}\mathsf{lg}^{\mathsf{lc}}_{\mathsf{D}_\mathsf{X}}(\mathcal{V})$ and $\mathcal{A}\mathsf{lg}_{\mathcal{D}_\mathsf{X}}(\mathcal{V})$. To see that it gives a Quillen pair, we use [18, Lemma 3.5], which just requires Proposition 4.16. Let us show that it is moreover a Quillen equivalence. It suffices to see that $\ell^*$ induces an equivalence between homotopy categories. One gets such an equivalence

$$\ell^*\colon \mathsf{Ho}\left(\mathcal{A}\mathsf{lg}_{\mathcal{D}_\mathsf{X}}(\mathcal{V})\right) \xrightarrow{\sim} \mathsf{Ho}\left(\mathcal{A}\mathsf{lg}^{\mathsf{lc}}_{\mathsf{D}_\mathsf{X}}(\mathcal{V})\right),$$

since the fibrant objects of $\mathcal{A}\mathsf{lg}^{\mathsf{lc}}_{\mathsf{D}_\mathsf{X}}(\mathcal{V})$ are known to be locally constant proj-fibrant $\mathsf{D}_\mathsf{X}$-algebras by Proposition 4.16 and since $\mathsf{D}_\mathsf{X} \to \mathcal{D}_\mathsf{X}$ exhibits $\mathcal{D}_\mathsf{X}$ as the $\infty$-localization of $\mathsf{D}_\mathsf{X}$ at the set of isotopy equivalences $\mathcal{J}_\mathsf{X}$ due to [2, Proposition 2.19]. Notice that $\mathsf{D}_\mathsf{X}$ and $\mathcal{D}_\mathsf{X}$ are admissible and their associated homotopy categories of algebras obtained via model structures coincide with those obtained from $\infty$-categorical techniques by Hypothesis 1.2 (see [43, Theorems 5.11 and 7.11]). □

Proposition 2.19 in [2] (or equivalently, [40, Lemma 5.4.5.11], [37, §5.2]) also proves that $\mathcal{E}_\mathsf{X}$ is the $\infty$-localization of $\mathsf{E}_\mathsf{X}$ at all unary operations. Therefore, the analogous left Bousfield localization for $\mathsf{E}_\mathsf{X}$-algebras encodes local constancy in Definition 2.10.

**Proposition 4.19.** *Assume Hypothesis 4.12 holds. Then, the projective model structure* $\mathcal{A}\mathsf{lg}_{\mathsf{E}_\mathsf{X}}(\mathcal{V})$ *admits a left Bousfield localization among left semimodel categories whose local objects are the locally constant* $\mathsf{E}_\mathsf{X}$*-algebras. We denote such a semimodel category by* $\mathcal{A}\mathsf{lg}^{\mathsf{lc}}_{\mathsf{E}_\mathsf{X}}(\mathcal{V})$*. Moreover,*

$$\ell_\sharp\colon \mathcal{A}\mathsf{lg}^{\mathsf{lc}}_{\mathsf{E}_\mathsf{X}}(\mathcal{V}) \rightleftarrows \mathcal{A}\mathsf{lg}_{\mathcal{E}_\mathsf{X}}(\mathcal{V})\colon \ell^*$$

*establishes a Quillen equivalence.*

**Aside: Extension model structure in the presence of local constancy.** Let us discuss a simple variation of the extension-construction concerning the discrete setting $\mathsf{E}_\mathsf{X} \hookrightarrow \mathsf{D}_\mathsf{X} \hookrightarrow \mathsf{M}_\mathsf{X}$. To handle both cases uniformly, we consider a map $\iota\colon \mathtt{B} \to \mathtt{N}$ of operads in sets (for simplicity) and we select a collection of unary maps $\mathcal{J}$ in $\mathtt{B}$. Then, Hypothesis 4.3 on the Quillen pair $\iota_\sharp\colon \mathcal{A}\mathsf{lg}_{\mathtt{B}}(\mathcal{V}) \rightleftarrows \mathcal{A}\mathsf{lg}_{\mathtt{N}}(\mathcal{V}) : \iota^*$ imply the same requirements on the descended Quillen adjunction $\iota_\sharp\colon \mathcal{A}\mathrm{lg}^{\mathrm{lc}}_{\mathtt{B}}(\mathcal{V}) \rightleftarrows \mathcal{A}\mathrm{lg}^{\mathrm{lc}}_{\mathtt{N}}(\mathcal{V}) :\iota^*$, where we have localized both sides with respect to the corresponding set of maps representing $\mathcal{J}$ as in Definition 4.15. Simply observe that the Quillen adjunction descends since $\iota^*$ preserves (acyclic) proj-fibrations and $\mathcal{J}$-locally constant algebras, plus the fact that cofibrations are unchanged under left Bousfield localization.

The same arguments proving the existence and properties of the extension model structure $\mathcal{A}\mathrm{lg}_{\mathtt{N}}(\mathcal{V})_{\mathtt{B}}$ yield that $\mathcal{A}\mathrm{lg}^{\mathrm{lc}}_{\mathtt{N}}(\mathcal{V})$ can be right Bousfield localized using the descended pair $\iota_\sharp\colon \mathcal{A}\mathrm{lg}^{\mathrm{lc}}_{\mathtt{B}}(\mathcal{V}) \rightleftarrows \mathcal{A}\mathrm{lg}^{\mathrm{lc}}_{\mathtt{N}}(\mathcal{V}) :\iota^*$. The result, which we denote $\mathcal{A}\mathrm{lg}^{\mathrm{lc}}_{\mathtt{N}}(\mathcal{V})_{\mathtt{B}}$, is a left semimodel structure satisfying:

- A map of $\mathcal{J}$-locally constant algebras is a weak equivalence if and only if it is an equivalence when evaluated on colors of $\mathtt{B}$.
- It is cofibrantly generated; the generators are the image of those in $\mathcal{A}\mathrm{lg}^{\mathrm{lc}}_{\mathtt{B}}(\mathcal{V})$ along $\iota_\sharp$.
- It has the same cofibrations as $\mathcal{A}\mathrm{lg}_{\mathtt{N}}(\mathcal{V})_{\mathtt{B}}$.

Therefore, we obtain the same semimodel structure as if we first perform the extension model construction and then left Bousfield localize to impose local constancy.

Moreover, we only know the existence of $\mathcal{A}\mathrm{lg}^{\mathrm{lc}}_{\mathtt{N}}(\mathcal{V})_{\mathtt{B}}$ as a semimodel structure since, without stronger conditions, we need to use left Bousfield localization without left properness. However, if we happen to know that $\mathcal{A}\mathrm{lg}^{\mathrm{lc}}_{\mathtt{N}}(\mathcal{V})$ is a full model category (e.g. if $\mathcal{A}\mathrm{lg}_{\mathtt{N}}(\mathcal{V})$ is left proper) $\mathcal{A}\mathrm{lg}^{\mathrm{lc}}_{\mathtt{N}}(\mathcal{V})$ will be so. Notice that such right Bousfield localization only modifies the cofibrations of $\mathcal{A}\mathrm{lg}^{\mathrm{lc}}_{\mathtt{N}}(\mathcal{V})$ (as well as the weak equivalences), but the factorization and lifting axioms for the new cofibrations are unobstructed since so are those of $\mathcal{A}\mathrm{lg}^{\mathrm{lc}}_{\mathtt{B}}(\mathcal{V})$.

**Local constancy on prefactorization algebras.** The content of this subsection subsumes a combination of two Bousfield localizations which finally produce the lc-factorization model structure on prefactorization algebras over $X$.

Since we need to make references to different model structures and Quillen pairs between them, we adopt the following:

**Notation 4.20.** Arrows in diagrams of model categories represent left Quillen functors. The symbol $\sim_Q$ will denote a Quillen equivalence within these diagrams. LBL and RBL are acronyms for left and right Bousfield localization respectively.

Let us denote

$$\begin{array}{ccc} \mathrm{M}_X & \xleftarrow{\ \iota\ } & \mathrm{D}_X \\ {\scriptstyle\tilde{\ell}}\downarrow & & \downarrow{\scriptstyle\ell} \\ \mathcal{M}_X & \xleftarrow[\ i\ ]{} & \mathcal{D}_X \end{array}$$

the involved operads and maps between them on the subsequent construction.

First, we perform a right Bousfield localization using Theorem 4.5 for the inclusion of operads $\mathrm{M}_X \overset{\iota}{\hookleftarrow} \mathrm{D}_X$. This way we get the extension model structure, $\mathrm{p}\mathcal{F}\mathrm{act}_X(\mathcal{V})_{\mathrm{D}_X}$.

Define the set of maps $\widetilde{\mathcal{S}}_X$ in $\mathtt{pFact}_X(\mathtt{V})$ as the image of $\mathcal{S}_X$ along $\mathbb{L}\iota_\sharp$. Using the existence theorems for left Bousfield localizations as how they were employed in Theorem 4.18, we further produce the lax lc-factorization (semi)model structure on $\mathtt{pFact}_X(\mathtt{V})$.

**Theorem 4.21.** *Assume that Hypothesis 4.11 holds. Then,* $\mathtt{pFact}_X(\mathtt{V})$ *supports the* lax lc-factorization model structure $\mathrm{lax}\,\mathcal{F}\mathrm{act}^{\mathrm{lc}}_X(\mathcal{V})$*, which is the left Bousfield localization at* $\widetilde{\mathcal{S}}_X$ *of the extension model structure* $\mathrm{p}\mathcal{F}\mathrm{act}_X(\mathcal{V})_{\mathrm{D}_X}$*. Moreover, this model structure sits into a commutative diagram of (semi)model categories*

$$
\begin{array}{ccccc}
\mathrm{p}\mathcal{F}\mathrm{act}_X(\mathcal{V}) & \xleftarrow{\mathrm{RBL}} & \mathrm{p}\mathcal{F}\mathrm{act}_X(\mathcal{V})_{\mathrm{D}_X} & \xleftarrow{\sim_Q} & \mathcal{A}\mathrm{lg}_{\mathrm{D}_X}(\mathcal{V}) \\
\Big\downarrow \widetilde{\ell}_\sharp & & \Big\downarrow \mathrm{LBL} & & \Big\downarrow \mathrm{LBL} \\
 & & \mathrm{lax}\,\mathcal{F}\mathrm{act}^{\mathrm{lc}}_X(\mathcal{V}) & \xleftarrow{\sim_Q} & \mathcal{A}\mathrm{lg}^{\mathrm{lc}}_{\mathrm{D}_X}(\mathcal{V}) \\
 & & \Big\downarrow \mathrm{Q}\wr & & \Big\downarrow \wr\mathrm{Q} \\
\mathcal{A}\mathrm{lg}_{\mathcal{M}_X}(\mathcal{V}) & \xleftarrow[\mathrm{RBL}]{} & \mathrm{lax}\,\mathcal{F}\mathrm{act}^{\mathrm{enr}}_X(\mathcal{V}) & \xleftarrow[\sim_Q]{} & \mathcal{A}\mathrm{lg}_{\mathcal{D}_X}(\mathcal{V})
\end{array}
\qquad (\ell_\sharp : \mathcal{A}\mathrm{lg}_{\mathrm{D}_X}(\mathcal{V}) \to \mathcal{A}\mathrm{lg}_{\mathcal{D}_X}(\mathcal{V})) .
$$

*If Hypothesis 4.12 holds instead, the left Bousfield localization* $\mathrm{lax}\,\mathcal{F}\mathrm{act}^{\mathrm{lc}}_X(\mathcal{V})$ *exists as a left semimodel category.*

*Proof.* Under Hypothesis 4.12, it is clear that we have the existence of such a semimodel structure and that it fits in a diagram as depicted in the statement. To see that the marked arrows are Quillen equivalences, we proceed as follows. By Theorem 4.5, both upper and lower horizontal arrows on the right are Quillen equivalences. The middle horizontal arrow is a Quillen equivalence since it comes from left Bousfield localization of equivalent semimodel structures at essentially the same set of maps (see the proof of [9, Theorem 4.16]). Theorem 4.18 says that the lower vertical map on the right is a Quillen equivalence. The remaining Quillen equivalence is so by 2-out of-3.

The existence as a model structure under Hypothesis 4.11 comes from the aside after Proposition 4.19.

□

An important consequence of the Quillen equivalence $\mathrm{lax}\,\mathcal{F}\mathrm{act}^{\mathrm{enr}}_X(\mathcal{V}) \simeq_Q \mathrm{lax}\,\mathcal{F}\mathrm{act}^{\mathrm{lc}}_X(\mathcal{V})$ is a recognition of bifibrant objects in the latter, Proposition 4.24. Before stating it, we present a pair of lemmas which are essential in our proof.

**Lemma 4.22.** *A prefactorization algebra* $\mathcal{F} \in \mathtt{pFact}_X(\mathtt{V})$ *is* $\mathrm{D}_X$*-colocal if, and only if, for any* $Z \in \mathrm{Open}(X)$*, the canonical map* $\operatorname{hocolim}_{\overline{\mathrm{D}}_X \downarrow Z} \iota^*\mathcal{F} \to \mathcal{F}(Z)$ *is an equivalence.*

*Proof.* Note that the endofunctor $\mathfrak{Q} = \mathbb{L}\iota_\sharp\iota^*$,

$$\mathbb{L}\iota_\sharp\iota^* \colon \mathrm{p}\mathcal{F}\mathrm{act}_X(\mathcal{V}) \xrightarrow{\ \iota^*\ } \mathcal{A}\mathrm{lg}_{\mathrm{D}_X}(\mathcal{V}) \xrightarrow[\text{cof.rep.}]{Q} \mathcal{A}\mathrm{lg}_{\mathrm{D}_X}(\mathcal{V}) \xrightarrow{\ \iota_\sharp\ } \mathrm{p}\mathcal{F}\mathrm{act}_X(\mathcal{V}),$$

can be described alternatively as:

$$\mathfrak{Q}\mathcal{F}(Z) \;=\; \mathbb{L}\iota_\sharp\iota^*\mathcal{F}(Z) \;\simeq\; \mathbb{L}\iota_!\iota^*\mathcal{F}(Z) \;=\; \operatorname*{hocolim}_{\overline{\mathrm{D}}_X \downarrow Z} \iota^*\mathcal{F},$$

where the middle equivalence holds because of Proposition A.5. □

With such a characterization, one can show the following result.

**Lemma 4.23.** *Any locally constant lax factorization algebra over* $X$ *is* $D_X$*-colocal.*

*Proof.* Let $\mathcal{F}$ be a locally constant lax factorization algebra over $X$. By Lemma 4.22, we must show that

$$\operatorname*{hocolim}_{\overline{D}_Z} \mathcal{F}|_{\overline{D}_Z} \longrightarrow \mathcal{F}(Z)$$

is an equivalence for any $Z \in \mathrm{Open}(X)$ (since $\overline{D}_X \downarrow Z \simeq \overline{D}_Z$). The idea is that we can compute the homotopy colimit on the left using a nice Weiss cover of $Z$ made by discs $(U_i)_{i\in I}$ (whose existence is given by Lemma 3.14). To do so, form the commutative square

$$\begin{array}{ccc} \operatorname*{hocolim}_{\overline{D}_Z} \mathcal{F}|_{\overline{D}_Z} & \longrightarrow & \mathcal{F}(Z) \\ \uparrow & & \uparrow \\ \operatorname*{hocolim}_{S\subseteq I} \operatorname*{hocolim}_{\overline{D}_{U_S}} \mathcal{F}|_{\overline{D}_{U_S}} & \longrightarrow & \operatorname*{hocolim}_{S\subseteq I} \mathcal{F}(U_S) \end{array}$$

and notice that the lower horizontal arrow is an equivalence since $U_S \in \overline{D}_{U_S}$ is final and that the right vertical map is an equivalence since $\mathcal{F}$ is a lax factorization algebra. Therefore, it suffices to check that the left vertical map is an equivalence. We now use that $\bar{\ell}\colon \overline{D}_Z \to \overline{\mathcal{D}}_Z$ exhibits $\overline{\mathcal{D}}_Z$ as the $\infty$-localization of $\overline{D}_Z$ at the isotopy equivalences [37, §5.2], and analogously for $U_S$ in place of $Z$, to consider the commutative square

$$\begin{array}{ccc} \operatorname*{hocolim}_{\overline{\mathcal{D}}_Z} \left(\mathcal{F}|_{\overline{D}_Z}\right)^{\mathrm{lc}} & \longrightarrow & \operatorname*{hocolim}_{\overline{D}_Z} \mathcal{F}|_{\overline{D}_Z} \\ \uparrow & & \uparrow \\ \operatorname*{hocolim}_{S\subseteq I} \operatorname*{hocolim}_{\overline{\mathcal{D}}_{U_S}} \left(\mathcal{F}|_{\overline{D}_{U_S}}\right)^{\mathrm{lc}} & \longrightarrow & \operatorname*{hocolim}_{S\subseteq I} \operatorname*{hocolim}_{\overline{D}_{U_S}} \mathcal{F}|_{\overline{D}_{U_S}} \end{array},$$

where $\left(\mathcal{F}|_{\overline{D}_Z}\right)^{\mathrm{lc}}$ is the essentially unique factorization of $\mathcal{F}|_{\overline{D}_Z}$ through $\bar{\ell}\colon \overline{D}_Z \to \overline{\mathcal{D}}_Z$ (resp. for $U_S$). Notice that we are implicitly using [43, Theorem 7.11]. Since both horizontal arrows are equivalences by homotopy cofinality of $\infty$-localizations [19, Proposition 7.1.10], it remains to see that the left vertical map is an equivalence. This follows from the following facts[ii]:

- The nested homotopy colimit on the lower left corner can be rewritten as

$$\operatorname*{hocolim}_{S\subseteq I} \operatorname*{hocolim}_{\overline{\mathcal{D}}_{U_S}} \left(\mathcal{F}|_{\overline{D}_{U_S}}\right)^{\mathrm{lc}} \simeq \operatorname*{hocolim}_{S\subseteq I} \mathbb{L}p_!\left(\mathcal{F}|_{\overline{D}_\bullet}\right)^{\mathrm{lc}} \simeq \operatorname*{hocolim}_{\int_S \overline{\mathcal{D}}_{U_S}} \left(\mathcal{F}|_{\overline{D}_\bullet}\right)^{\mathrm{lc}},$$

  where $p\colon \int_S \overline{\mathcal{D}}_{U_S} \to \{S \subseteq I :\ \text{finite and non-empty}\}^{\mathrm{op}}$ is the ($\infty$-)cocartesian fibration associated to the functor $S \mapsto \overline{\mathcal{D}}_{U_S}$, $\mathbb{L}p_!$ denotes derived left Kan extension along $p$ and $\left(\mathcal{F}|_{\overline{D}_\bullet}\right)^{\mathrm{lc}}$ is the restriction along $p$ of

$$\operatorname*{hocolim}_{S\subseteq I} \overline{\mathcal{D}}_{U_S} \xrightarrow{\ \nu\ } \overline{\mathcal{D}}_Z \xrightarrow{\left(\mathcal{F}|_{\overline{D}_Z}\right)^{\mathrm{lc}}} \mathcal{V}\,.$$

  Notice that we use that $p$ is cocartesian to obtain the first equivalence, i.e. to ensure that $\overline{\mathcal{D}}_{U_S} \to \left(\int_S \overline{\mathcal{D}}_{U_S}\right)_{/S}$ is homotopy cofinal.

[ii] The first two are borrowed from Haugseng's answer to MO:question370788.

- By [26, Theorem 7.4], the functor $\int_S \overline{\mathcal{D}}\mathfrak{U}_S \longrightarrow \text{hocolim}_{S\subseteq I} \overline{\mathcal{D}}\mathfrak{U}_S$ exhibits its target as the $\infty$-localization at the cocartesian arrows of its domain. Thus, it is homotopy cofinal.
- Due to Lemma 3.10, one shows that the functor $\nu$ above is an equivalence. One possibility is to see both $\infty$-categories as domains of right fibrations over $\mathcal{D}\text{iscs}_n$ and use that lemma to compare fibers. Combined with the previous facts, this implies that the left vertical map in the previous square is an equivalence, concluding the proof. □

Due to Lemmas 4.22 and 4.23, it is possible to produce a cosheafification machine for the Weiss topology that preserves the appropriate algebraic structure, i.e. a functorial construction that turns precosheaves into cosheaves, hence this answers Gwilliam-Rejzner problem [29, Remark 2.33]. Such cosheafification articulates the recognition of bifibrant objects for $\text{lax}\,\mathcal{F}\text{act}_X^{\text{lc}}(\mathcal{V})$ and its existence is important on its own. One reason is that in usual categories as sets or abelian groups, cosheafification is only known to exist by means of general adjoint functor theorems and so, we have little control on its properties (see [45]).

**Proposition 4.24.** *The bifibrant objects in the lax lc-factorization (left semi)model structure are the proj-bifibrant locally constant lax factorization algebras.*

*Proof.* Recall that the lax lc-factorization model structure is constructed using two Bousfield localizations,

$$\text{p}\mathcal{F}\text{act}_X(\mathcal{V}) \xleftarrow{\text{RBL}} \text{p}\mathcal{F}\text{act}_X(\mathcal{V})_{D_X} \xrightarrow{\text{LBL}} \text{lax}\,\mathcal{F}\text{act}_X^{\text{lc}}(\mathcal{V}).$$

Right (resp. left) Bousfield localization does not affect fibrant (resp. cofibrant) objects. Hence, the cofibrant objects in $\text{lax}\,\mathcal{F}\text{act}_X^{\text{lc}}(\mathcal{V})$ are $D_X$-colocal proj-cofibrant $M_X$-algebras while the fibrant objects in $\text{lax}\,\mathcal{F}\text{act}_X^{\text{lc}}(\mathcal{V})$ are locally constant proj-fibrant $M_X$-algebras by Theorem 4.18. In the left semimodel case, we do not know a characterization of fibrant objects, but at least we know that the bifibrant objects are the locally constant bifibrant objects in $\text{p}\mathcal{F}\text{act}_X(\mathcal{V})_{D_X}$, and that suffices for our purposes.

Choose a bifibrant object $\mathcal{F}$ in $\text{lax}\,\mathcal{F}\text{act}_X^{\text{lc}}(\mathcal{V})$. We will show that $\mathcal{F}$ is a locally constant lax factorization algebra by constructing one, $\widetilde{\mathcal{F}}$, which is colorwise equivalent to $\mathcal{F}$. The main idea is to make use of the Quillen equivalence

$$\widetilde{\ell}_\sharp\colon \text{lax}\,\mathcal{F}\text{act}_X^{\text{lc}}(\mathcal{V}) \rightleftarrows \text{lax}\,\mathcal{F}\text{act}_X^{\text{enr}}(\mathcal{V})\colon \widetilde{\ell}^*$$

appearing in Theorem 4.21. Let us denote $\mathcal{F}_{\text{loc}} = \widetilde{\ell}^*\widetilde{\ell}_\sharp\mathcal{F}$ for simplicity.

Since $\mathcal{F}$ is cofibrant for the lax factorization model structure, $\widetilde{\ell}_\sharp\mathcal{F}$ is cofibrant in the enriched lax factorization model category and therefore it is an enriched lax factorization algebra by Proposition 4.6. In particular, its restriction $\mathcal{F}_{\text{loc}}$ along $\widetilde{\ell}$ is a locally constant lax factorization algebra. Moreover, the unit of the Quillen equivalence evaluated at $\mathcal{F}$ yields an equivalence $\mathcal{F} \to \mathcal{F}_{\text{loc}}$ in $\text{lax}\,\mathcal{F}\text{act}_X^{\text{lc}}(\mathcal{V})$, which is not, a priori, a colorwise equivalence. However, the class of equivalences between bifibrant objects in $\text{lax}\,\mathcal{F}\text{act}_X^{\text{lc}}(\mathcal{V})$ is the class of colorwise equivalences, since (co)local equivalences between (co)local objects are ordinary equivalences in a (right) left Bousfield localization.

Thus, it will suffice to show that $\mathcal{F}_{\text{loc}}$ is bifibrant in $\text{lax}\,\mathcal{F}\text{act}_X^{\text{lc}}(\mathcal{V})$. This is not always the case, but it will be after slightly modifying this algebra.

Recall that bifibrant in $\mathrm{lax}\mathcal{F}\mathrm{act}^{\mathrm{lc}}_{\mathrm{X}}(\mathcal{V})$ means proj-bifibrant locally constant $\mathrm{D}_{\mathrm{X}}$-colocal algebra as explained at the beginning of this proof. Let us focus on $\mathrm{D}_{\mathrm{X}}$-colocality. By Lemma 4.22, this property is equivalent to the following codescent condition: for any $\mathrm{Z} \in \mathrm{col}(\mathrm{M}_{\mathrm{X}})$, the canonical map

$$\operatorname*{hocolim}_{\overline{\mathrm{D}}_{\mathrm{X}}\downarrow \mathrm{Z}} \mathcal{F}_{\mathrm{loc}}(\star) \longrightarrow \mathcal{F}_{\mathrm{loc}}(\mathrm{Z})$$

is an equivalence. Taking into account that $\overline{\mathrm{D}}_{\mathrm{X}}\downarrow \mathrm{Z} \simeq \overline{\mathrm{D}}_{\mathrm{Z}}$, this codescent condition is deduced from $\mathcal{D}_{\mathrm{X}}$-colocality in exactly the same way as Weiss codescent in Proposition 4.6 replacing Lemma 3.10 by Lemma 3.12.

It remains to add the projective conditions on $\mathcal{F}_{\mathrm{loc}}$. Denoting R (resp. Q) a proj-fibrant (resp. proj-cofibrant) replacement in $\mathrm{p}\mathcal{F}\mathrm{act}_{\mathrm{X}}(\mathcal{V})$, we get a diagram

$$\begin{array}{ccccc} & & \mathrm{Q}\,\mathrm{R}\,\mathcal{F}_{\mathrm{loc}} & = & \widetilde{\mathcal{F}} \\ & \nearrow & \big\downarrow \wr & & \\ \mathcal{F} & \longrightarrow & \mathrm{R}\,\mathcal{F}_{\mathrm{loc}} & & \end{array} ,$$

where the lift exists by lifting properties in the projective model structure and it is an equivalence in $\mathrm{lax}\mathcal{F}\mathrm{act}^{\mathrm{lc}}_{\mathrm{X}}(\mathcal{V})$ by 2-out of-3. We conclude the proof since $\widetilde{\mathcal{F}}$ is a bifibrant object in $\mathrm{lax}\mathcal{F}\mathrm{act}^{\mathrm{lc}}_{\mathrm{X}}(\mathcal{V})$ and, by construction, colorwise equivalent to $\mathcal{F}$.

For the converse, pick a proj-bifibrant locally constant lax factorization algebra $\mathcal{F}'$. By the discussion at the beginning of the proof, we are reduced to check that $\mathcal{F}'$ is $\mathrm{D}_{\mathrm{X}}$-colocal and this was proven in Lemma 4.23. □

A trivial adaptation of the above arguments replacing $\mathrm{D}_{\mathrm{X}}$ by $\mathrm{E}_{\mathrm{X}}$, leads to:

**Theorem 4.25.** *Assume that Hypothesis 4.11 holds. Then,* $\mathtt{pFact}_{\mathrm{X}}(\mathcal{V})$ *supports the* lc-factorization model structure $\mathcal{F}\mathrm{act}^{\mathrm{lc}}_{\mathrm{X}}(\mathcal{V})$, *which is the left Bousfield localization at all unary maps in* $\mathrm{E}_{\mathrm{X}}$ *of the extension model structure* $\mathrm{p}\mathcal{F}\mathrm{act}_{\mathrm{X}}(\mathcal{V})_{\mathrm{E}_{\mathrm{X}}}$. *The lc-factorization model structure sits in a commutative square of (semi)model categories and Quillen equivalences*

$$\begin{array}{ccc} \mathcal{F}\mathrm{act}^{\mathrm{lc}}_{\mathrm{X}}(\mathcal{V}) & \xleftarrow{\sim\mathrm{Q}} & \mathcal{A}\mathrm{lg}^{\mathrm{lc}}_{\mathrm{E}_{\mathrm{X}}}(\mathcal{V}) \\ {\scriptstyle \mathrm{Q}\wr}\big\downarrow & & \big\downarrow{\scriptstyle \wr\mathrm{Q}} \\ \mathcal{F}\mathrm{act}^{\mathrm{enr}}_{\mathrm{X}}(\mathcal{V}) & \xleftarrow[\sim\mathrm{Q}]{} & \mathcal{A}\mathrm{lg}_{\mathcal{E}_{\mathrm{X}}}(\mathcal{V}) \end{array} .$$

*Moreover, the bifibrant objects in the lc-factorization model structure are the proj-bifibrant locally constant factorization algebras over* X.

*If Hypothesis 4.12 holds instead, the left Bousfield localization* $\mathcal{F}\mathrm{act}^{\mathrm{lc}}_{\mathrm{X}}(\mathcal{V})$ *exists as a left semimodel category. The characterization of bifibrant objects also holds.*

# 5 Summary, variations and generalizations

We have developed model structures that present different kinds of factorization algebras for smooth manifolds, but our methods can be applied in different situations. The choice of ordinary smooth manifolds is due to its simplicity.

Let us summarize our results and their requirements in each setting.

- **Enriched (lax) factorization model structure:** they are applications of the extension model machine (Theorem 4.5). We perform the pertinent analysis of the endofunctor $\mathfrak{Q} = \mathbb{L}\iota_\sharp\iota^*$, which requires a computation of factorization homology with context (Theorem 3.19) which follows from the existence of good Weiss covers (Lemma 3.14) and a cosheaf condition on embedding spaces (Lemma 3.10). The description of $\mathfrak{Q}$ and gluing properties of (lax) factorization algebras are combined to recognize cofibrant objects in Propositions 4.6 and 4.8.

- **(Lax) lc-factorization (left semi)model structure:** these model structures come from a combination of the right Bousfield localization given by the extension model structure together with a left Bousfield localization imposing local constancy. Left Bousfield localization produces model structures under left properness and left semimodel structures without such a hypothesis, so we should distinguish between two classes of homotopy cosmoi depending on this difference. In addition, we need to identify $\mathsf{D}_\mathsf{X} \to \mathcal{D}_\mathsf{X}$ (resp. $\mathsf{E}_\mathsf{X} \to \mathcal{E}_\mathsf{X}$) as an $\infty$-localizations of $\infty$-operads ([2, Proposition 2.19], [40, Lemma 5.4.5.11] or [37, §5.2]) to find the set of maps at which localize.

  On the other hand, one has to recognize the colocality condition for the extension model structure associated to $\mathsf{D}_\mathsf{X} \to \mathsf{M}_\mathsf{X}$ (resp. $\mathsf{E}_\mathsf{X} \to \mathsf{M}_\mathsf{X}$). This condition boils down to a codescent property (Lemma 4.22). Since embedding spaces satisfy that property (Lemma 3.12), enriched (lax) factorization algebras satisfy this variant of codescent. Proposition 4.24 implies that this suffices to identify bifibrant objects as locally constant (lax) factorization algebras.

At the end, and just for the record, we get a commutative diagram of (left semi)model categories and left Quillen functors between them

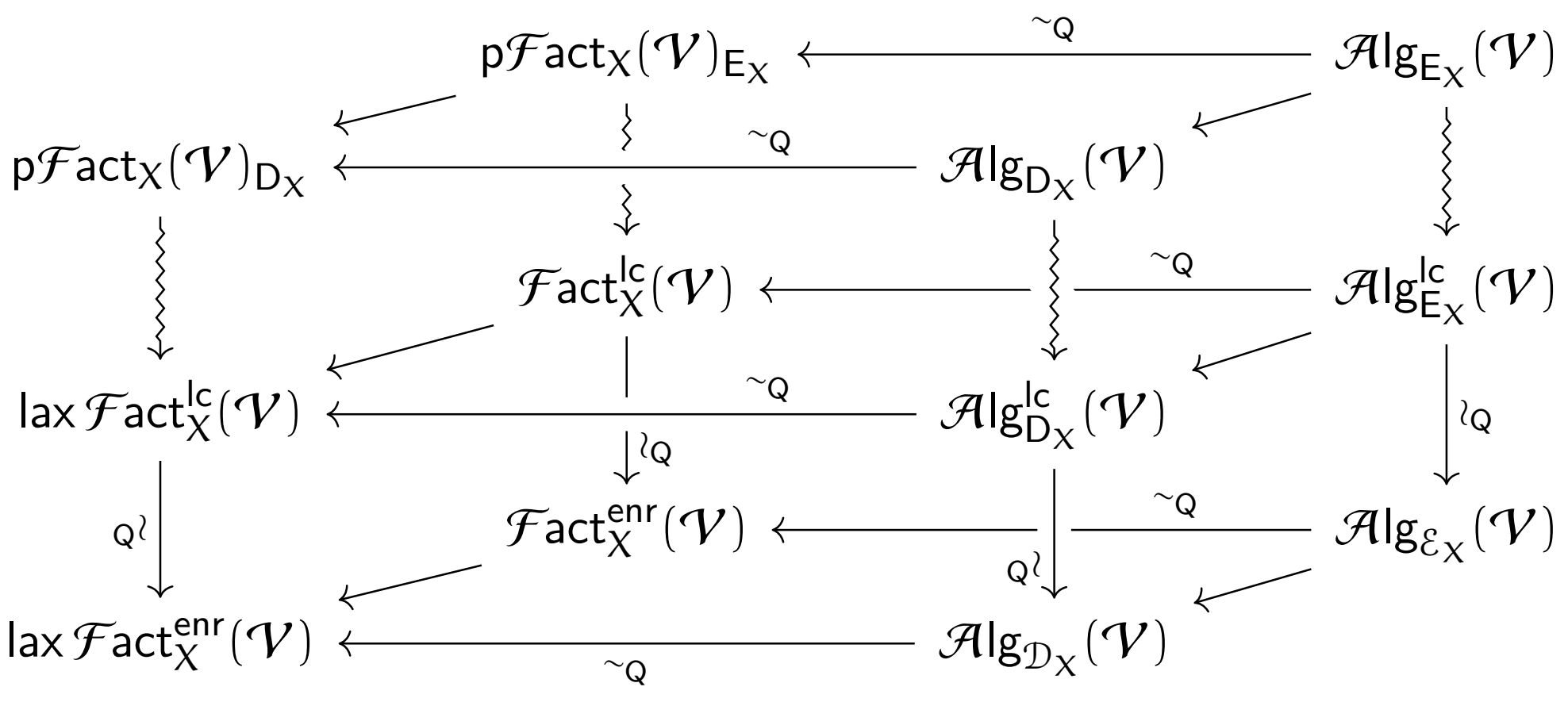

where the snaky arrows represent left Quillen functors of left Bousfield localizations.

Next, we collect some settings where one or both of these schemes are valid. They consist of varying the diagram of operads

$$\begin{array}{ccccc} \mathsf{M}_\mathsf{X} & \hookleftarrow & \mathsf{D}_\mathsf{X} & \hookleftarrow & \mathsf{E}_\mathsf{X} \\ \downarrow & & \downarrow & & \downarrow \\ \mathcal{M}_\mathsf{X} & \hookleftarrow & \mathcal{D}_\mathsf{X} & \hookleftarrow & \mathcal{E}_\mathsf{X} \end{array} .$$

**Tangent structures and context-free factorization algebras.** As customary in the field, one may add tangential structures on manifolds by fixing a continuous map $\tau\colon B \to BO(n)$. Then, a $n$-manifold with a B-structure is a manifold $Y \in \mathcal{M}\mathrm{fld}_n$ equipped with a factorization (up to homotopy) of their tangent bundle classifying map

$$[T\,Y]\colon Y \xrightarrow{\xi_Y} B \xrightarrow{\tau} BO(n).$$

Embedding spaces for manifolds with B-tangent structure, i.e. embeddings preserving the tangential structure up to homotopy, fit into homotopy pullback squares

$$\begin{array}{ccc} \mathrm{Emb}^B(Y,Z) & \longrightarrow & \mathrm{Emb}(Y,Z) \\ \downarrow & {}^{h}\lrcorner & \downarrow{\scriptstyle \text{tg-classifier}} \\ \mathrm{Map}_{/B}(Y,Z) & \xrightarrow[\tau]{} & \mathrm{Map}_{/BO(n)}(Y,Z) \end{array},$$

where $\mathrm{Map}_{/B}(\star,\bullet)$ denotes the mapping space in the homotopical slice category $\mathtt{Top}_{/B}$. The corresponding $\mathtt{Top}$-categories $\mathcal{M}\mathrm{fld}_n^B$, $\mathcal{D}\mathrm{iscs}_n^B$ inherit sm-structures. As expected, $\mathcal{E}_n^B$ is the full suboperad of $\mathcal{M}\mathrm{fld}_n^B$ on the color $\mathbb{D}$, which is equipped with the B-structure coming from its natural framing.

Applying the methods explained above, one gets a model structure

$$\mathcal{A}\mathrm{lg}_{\mathcal{M}\mathrm{fld}_n^B}(\mathcal{V})_{\mathcal{E}_n^B} \quad \left(\text{resp. } \mathcal{A}\mathrm{lg}_{\mathcal{M}\mathrm{fld}_n^B}(\mathcal{V})_{\mathcal{D}\mathrm{iscs}_n^B}\right)$$

presenting the homotopy theory of enriched (lax) factorization algebras with B-tangent structure. Of course, for $\tau = \mathrm{id}_{BO(n)}$, the resulting theories are what one obtains using

$$\mathcal{E}_n \hookrightarrow \mathcal{D}\mathrm{iscs}_n \hookrightarrow \mathcal{M}\mathrm{fld}_n\,.$$

*Remark* 5.1. The model structure $\mathcal{A}\mathrm{lg}_{\mathcal{M}\mathrm{fld}_n^B}(\mathcal{V})_{\mathcal{E}_n^B}$ presents the $\infty$-category

$$H(\mathcal{M}\mathrm{fld}_n^B, \mathcal{V}_\infty)$$

of factorization homology theories for B-structured $n$-manifolds, as defined in [2].

In these terms, it is quite natural to state that $\mathcal{E}_n$-algebras produce factorization algebras via factorization homology (see Theorem 2.14). Consider $X \in \mathcal{M}\mathrm{fld}_n^B$. Then, there is an induced diagram of operads

$$\begin{array}{ccccc} \mathrm{Disj}(X) & \hookleftarrow & \mathrm{Discs}(X) & \hookleftarrow & E_X \\ \downarrow & & \downarrow & & \downarrow \\ \mathcal{M}\mathrm{fld}_n^B & \hookleftarrow & \mathcal{D}\mathrm{iscs}_n^B & \hookleftarrow & \mathcal{E}_n^B \end{array}.$$

The claim is that

$$\mathcal{A}\mathrm{lg}_{\mathcal{E}_n^B}(\mathcal{V})_\infty \xrightarrow{\text{fact.hom.}} \mathcal{A}\mathrm{lg}_{\mathcal{M}\mathrm{fld}_n^B}(\mathcal{V})_\infty \xrightarrow{\text{rest.}} \mathrm{p}\mathcal{F}\mathrm{act}_X(\mathcal{V})_\infty$$

lands in factorization algebras. Using the model structures already introduced, this claim follows since the restriction functor associated to $\mathrm{Disj}(X) \to \mathcal{M}\mathrm{fld}_n^B$ induces a commutative square of model categories

$$\begin{array}{ccc} \mathcal{A}\mathrm{lg}_{\mathcal{M}\mathrm{fld}_n^B}(\mathcal{V})_{\mathcal{E}_n^B} & \xrightarrow{\text{rest.}} & \mathcal{F}\mathrm{act}_X^{\mathrm{lc}}(\mathcal{V}) \\ \downarrow & & \downarrow \\ \mathcal{A}\mathrm{lg}_{\mathcal{M}\mathrm{fld}_n^B}(\mathcal{V})_{\mathcal{D}\mathrm{iscs}_n^B} & \xrightarrow[\text{rest.}]{} & \mathrm{lax}\,\mathcal{F}\mathrm{act}_X^{\mathrm{lc}}(\mathcal{V}) \end{array}.$$

**Manifolds with boundary.** Define $\mathcal{M}\mathrm{fld}_n^\partial$ to be the $\mathrm{Top}$-category of $n$-manifolds possibly with boundary, $\partial X \hookrightarrow X$, and mapping spaces $\mathrm{Emb}^\partial(X, Y)$ given by those embeddings preserving the boundary. Of course, $\mathcal{M}\mathrm{fld}_n^\partial$ has a sm-structure given by disjoint union. Now, there are two building blocks, or local models, of manifolds with boundary

$$\mathbb{D} = \{x \in \mathbb{R}^n : \|x\| < 1\} \quad \text{and} \quad \mathbb{D}_+ = \{x \in \mathbb{D} : x_1 \geqslant 0\}.$$

Hence, $\mathcal{D}\mathrm{iscs}_n^\partial$ is the full subcategory of $\mathcal{M}\mathrm{fld}_n^\partial$ spanned by the objects $\{\mathbb{D}^{\sqcup p} \sqcup \mathbb{D}_+^{\sqcup q}\}_{0 \leqslant p,q < \infty}$ and $\mathcal{E}_n^\partial$ is the full suboperad of $\mathcal{M}\mathrm{fld}_n^\partial$ on the colors $\{\mathbb{D}, \mathbb{D}_+\}$.

For $X \in \mathcal{M}\mathrm{fld}_n^\partial$, one can consider variants with context of the previous objects (yielding operads as in §3):

$$\mathcal{M}_X^\partial \hookleftarrow \mathcal{D}_X^\partial \hookleftarrow \mathcal{E}_X^\partial,$$

and also discrete versions:

$$\mathrm{Disj}(X) \hookleftarrow \mathrm{D}_X^\partial \hookleftarrow \mathrm{E}_X^\partial,$$

where the colors of $\mathrm{D}_X$ (resp. $\mathrm{E}_X^\partial$) are those $Z \in \mathcal{M}\mathrm{fld}_{n/X}^\partial$ diffeomorphic/homeomorphic to $\mathbb{D}^{\sqcup p} \sqcup \mathbb{D}_+^{\sqcup q}$ with $0 \leqslant p, q < \infty$ (resp. $\mathbb{D}$ or $\mathbb{D}_+$).

As mentioned at the beginning of the current section, our model structures exist in this content if we check: $(i)$ sufficient codescent properties for embedding spaces, $(ii)$ the existence of good Weiss covers and $(iii)$ identification of $\infty$-localizations of operads.

$(i)$ The required properties of embedding spaces to make the theory work in this case are either deduced as those without boundary or additionally taking into account the following result ([36, Proposition 6.5]):

**Proposition 5.2.** *The derivative map*

$$\mathrm{Emb}^\partial(\mathbb{D}^{\sqcup p} \sqcup \mathbb{D}_+^{\sqcup q}, Z) \longrightarrow \mathrm{Fr}\big(\mathrm{Conf}_p(Z \backslash \partial Z)\big) \times \mathrm{Fr}\big(\mathrm{Conf}_q(\partial Z)\big)$$

*is a mixed acyclic fibration (i.e. a weak homotopy equivalence and a Hurewicz fibration), where* $\mathrm{Fr}(Y)$ *denotes the frame-bundle associated to a manifold* $Y$.

$(ii)$ Any manifold with boundary $X$ admits good Weiss covers by elements in $\mathrm{D}_X^\partial$. Use [11, Proposition 9.1] combined with a good Weiss cover of the (tubular neighborhood) of the boundary.

$(iii)$ The recognition of $\mathrm{E}_X^\partial \to \mathcal{E}_X^\partial$ (resp. $\mathrm{D}_X^\partial \to \mathcal{D}_X^\partial$) as an $\infty$-localization is also provided by [2]. It is possible to check this fact using weak approximations of $\infty$-operads (see [30]).

*Remark* 5.3. Tangent structures on manifolds with boundary can be considered as well. Defining the relevant embedding spaces is a little more involved, but in the end, it consists of using certain homotopy limits to modify $\mathrm{Emb}^\partial(\star, \bullet)$ (see [3, 30]).

**Hyperfactorization algebras and topological manifolds.** Much of the theory (and variations exposed above) is adaptable to work for topological manifolds as well, although the absence of differential techniques complicates a bit the discussion. For instance, for topological $n$-manifolds, the space of embeddings $\mathsf{Emb}(X, Y)$ is the space of injective open continuous maps equipped with the compact-open topology; to identify the homotopy type of, let us say, $\mathsf{Emb}(\mathbb{D}^{\sqcup k}, Z)$, instead of using the derivative map, one requires a bit of machinery dealing with germs of embeddings (see [30] or the original reference [40]).

The main difference with respect to the rest of the mentioned generalizations concerns the existence of good Weiss covers for topological manifolds. We are not currently aware of such an existence and this obstructs the previous line of discussion. However, what is easier to prove is that any topological manifold admits a *good Weiss hypercover*, i.e. an open hypercover whose levels are disjoint unions of open subsets homeomorphic to $\mathbb{D}^{\sqcup k}$ with $0 \leqslant k < \infty$. Consequently, our methods can be applied to find model structures to study locally constant (lax) *hyperfactorization algebras* over topological manifolds. These objects are prefactorization algebras which instead of *homotopical codescent for Weiss covers* satisfy *homotopical codescent with respect to Weiss hypercovers*. For concreteness: let $U_\bullet \colon \Delta^{\mathrm{op}} \to \mathsf{Open}(X)$ be a Weiss hypercover of $U \subseteq X$. Then, $\mathcal{F} \in \mathtt{pFact}_X(\mathtt{V})$ satisfies homotopy codescent for $U_\bullet$ if

$$\operatorname*{hocolim}_{r\in\Delta^{\mathrm{op}}} \mathcal{F}(U_r) \longrightarrow \mathcal{F}(U)$$

is an equivalence.

*Remark* 5.4. It might be the case that even without a result providing good Weiss covers for topological manifolds, the homotopy theory of hyperfactorization algebras and that of factorization algebras coincide. We have not explored this problem yet.

**Conically smooth stratified manifolds.** Let us very concisely collect the results that would make the construction of model structures for constructible (lax) hyperfactorization algebras on conically smooth stratified manifolds ([4, 5]) work[i]. By considering finite disjoint unions of basics ([5, Definition 2.2.1]) one can construct Weiss hypercovers due to [5, Proposition 3.2.23]. Codescent with respect to those hypercovers can be deduced from the proof of [5, Lemma 6.1.1] together with a configuration space argument, simply by evaluation on finite disjoint unions of basics. The recognition of the relevant $\infty$-localizations of $\infty$-operads in this setting is discussed at [4, Proposition 2.22] (or [37, §5.2] for a corrected version).

# A Partial monoidal categories

In this appendix, we define a mild generalization of symmetric monoidal (sm) categories, called *partial monoidal categories* or simply *pm-categories*, which share important properties with them. Our main result is that the derived counterparts of categorical Kan extension and operadic Kan extension along a *strong monoidal functor* between pm-categories coincide, see Proposition A.5. This result is fundamental to perform factorization homology computations.

Let us start by exploring a bit more the relation between sm-categories and operads.

[i] As in [37], we can work with constructible factorization algebras (not only *hyper*) if we consider conically smooth stratified manifolds admitting a *disk-factorizing basis*.

Given a sm-sSet-category $\mathtt{U}$, there is an associated $\mathtt{sSet}$-operad $\mathtt{U}_{\mathsf{opd}}$ whose $\mathtt{sSet}$-objects of multimorphisms are given by

$$\mathtt{U}_{\mathsf{opd}}\begin{bmatrix}\{X_i\}_i\\ Y\end{bmatrix} = \underline{\mathrm{Hom}}_{\mathtt{U}}\Big(\bigotimes_i X_i, Y\Big).$$

The composition products are defined combining the enriched symmetric monoidal structure and the $\mathtt{sSet}$-category structure. If we consider the category of symmetric monoidal $\mathtt{sSet}$-categories equipped with lax monoidal functors as morphisms, $\mathtt{smCat}(\mathtt{sSet})$, the assignment $(\star)_{\mathsf{opd}}\colon \mathtt{smCat}(\mathtt{sSet}) \to \mathtt{Opd}(\mathtt{sSet})$ is functorial. One also observes:

- The essential image of $(\star)_{\mathsf{opd}}$ can be characterized. More concretely, an operad $\mathcal{O}$ comes from a sm-sSet-category if (i) the functor $d \mapsto \mathcal{O}\left[\begin{smallmatrix}\underline{c}\\ d\end{smallmatrix}\right]$ is representable for any finite list of colors $\underline{c} = \{c_i\}_{i\in I}$ in $\mathrm{col}(\mathcal{O})$, and (ii) any canonical map associated to a function $f\colon I \to J$ of finite sets by representability

  $$\bigotimes_{i\in I} c_i \longrightarrow \bigotimes_{j\in J}\Big(\bigotimes_{i\in f^{-1}(j)} c_i\Big)$$

  is an isomorphism (see [32] or [34] for the set-valued case).

- The functor $(\star)_{\mathsf{opd}}$ is fully-faithful. That is,

  $$\mathtt{LaxMon}\left[\mathtt{U}', \mathtt{U}''\right] \cong \mathtt{Opd}(\mathtt{sSet})(\mathtt{U}'_{\mathsf{opd}}, \mathtt{U}''_{\mathsf{opd}})$$

  is a bijection for any pair of sm-sSet-categories[i]. This fact follows by naturality and representability of the $\mathtt{sSet}$-objects of multioperations.

- The functor $(\star)_{\mathsf{opd}}$ admits a left-adjoint called *enveloping sm-category* and denoted $\mathsf{Env}_{\otimes}\colon \mathtt{Opd}(\mathtt{sSet}) \to \mathtt{smCat}(\mathtt{sSet})$, where the target category is now considered with strong monoidal functors as morphisms. Given an operad $\mathcal{O}$, the enveloping sm-category $\mathsf{Env}_{\otimes}(\mathcal{O})$ is defined as follows:

  - $\mathrm{ob}(\mathsf{Env}_{\otimes}(\mathcal{O})) = \mathrm{ob}(\mathsf{Fin} \downarrow \mathrm{col}(\mathrm{O}))$, where Fin denotes the category of finite sets;
  - for $\underline{c} = \{c_i\}_{i\in I}$ and $\underline{d} = \{d_j\}_{j\in J}$,

    $$\underline{\mathrm{Hom}}_{\mathsf{Env}_{\otimes}(\mathcal{O})}(\underline{c}, \underline{d}) = \coprod_{f\colon I\to J}\prod_{j\in J}\mathcal{O}\begin{bmatrix}\{c_i\}_{i\in f^{-1}(j)}\\ d_j\end{bmatrix},$$

    where the coproduct runs over maps $f\colon I \to J$ in Fin;
  - The monoidal structure comes from the functor defined on objects by

    $$\begin{array}{ccc}\mathsf{Env}_{\otimes}(\mathcal{O})^{\times 2} & \longrightarrow & \mathsf{Env}_{\otimes}(\mathcal{O})\\ (\underline{c}, \underline{d}) & \longmapsto & \underline{c} \boxplus \underline{d} := \{c_i\}_{i\in I} \amalg \{d_j\}_{j\in J}\end{array}.$$

[i] In fact, $(\star)_{\mathsf{opd}}$ can be upgraded to a fully-faithful 2-functor, yielding equivalences of categories between lax monoidal functors and operadic algebras.

- The unit $\mathcal{O} \to \mathsf{Env}_\otimes(\mathcal{O})_{\mathsf{opd}}$ identifies the category of $\mathcal{O}$-algebras with the category of *strong monoidal functors* $\mathtt{StMon}[\mathsf{Env}_\otimes(\mathcal{O}), \mathtt{V}]$. That is, it yields an equivalence of categories
$$\mathtt{Alg}_{\mathcal{O}}(\mathtt{V}) \simeq \mathtt{StMon}\left[\mathsf{Env}_\otimes(\mathcal{O}), \mathtt{V}\right],$$
via the identification $\mathtt{Alg}_{\mathsf{Env}_\otimes(\mathcal{O})_{\mathsf{opd}}}(\mathtt{V}) = \mathtt{LaxMon}\left[\mathsf{Env}_\otimes(\mathcal{O}), \mathtt{V}\right]$. See Remark A.3.

Observing that ordinary left Kan extension preserves (strong) monoidality, Horel in [36] applied the former ideas to compute operadic Kan extension along $\mathcal{O} \to \mathcal{P}$ via ordinary Kan extension along $\mathsf{Env}_\otimes(\mathcal{O}) \to \mathsf{Env}_\otimes(\mathcal{P})$. He also showed how to compute derived operadic Kan extensions in these terms. A related result in the $\infty$-categorical framework was given in [4, Lemma 2.16]. We are going to present a generalization of those conclusions to cover partial monoidal categories. This broader generality is an essential tool to perform *factorization homology* computations. See [40] for a thorough treatment of sm-envelopes in the higher categorical context.

It is time to define partial monoidal categories. Roughly speaking, they are symmetric monoidal categories whose objects cannot be always multiplied. See [38, 48] for definitions of partial algebraic structures. We have chosen to define partial monoidal categories using the fully-faithful functor $(\star)_{\mathsf{opd}}\colon \mathtt{pmCat}(\mathtt{sSet}) \to \mathtt{Opd}(\mathtt{sSet})$ for brevity.

**Definition A.1.** Let $\mathtt{B}$ be a $\mathtt{sSet}$-category.

- A *partial (symmetric) monoidal structure*, or *pm-structure*, on $\mathtt{B}$ is given by an operad $\mathcal{B} \equiv \mathtt{B}_{\mathsf{opd}}$ with underlying category $\overline{\mathcal{B}} = \mathtt{B}$ which satisfies: $(i)$ the functor $b \mapsto \mathcal{B}\left[\begin{smallmatrix}\{a_i\}_i\\ b\end{smallmatrix}\right]$ is either representable by $\boxtimes_i a_i \in \mathtt{B}$, i.e. it is naturally isomorphic to $b \mapsto \mathtt{B}\left[\begin{smallmatrix}\boxtimes_i a_i\\ b\end{smallmatrix}\right]$, or constant with value $\mathbb{0} \in \mathtt{sSet}$, and $(ii)$ the canonical maps associated to a function $f\colon I \to J$ in Fin by representability
$$\boxtimes_{i\in I} b_i \longrightarrow \boxtimes_{j\in J}\Big(\boxtimes_{i\in f^{-1}(j)} b_i\Big)$$
are isomorphisms (whenever representatives exist, and the existence of the right-hand side implies the existence of the left-hand side).[ii]

- A *lax monoidal functor* $g\colon \mathtt{B} \to \mathtt{B}'$ between pm-$\mathtt{sSet}$-categories is just a map between the associated operads $\mathcal{B} \to \mathcal{B}'$. We say that the lax monoidal functor $g$ is *strong* if $(i)$ the morphisms in $\mathtt{B}'$,
$$\boxtimes_i g(b_i) \longrightarrow g\left(\boxtimes_i b_i\right),$$
associated by representability to the image through $g$ of the universal multimorphisms $\Phi_{\{b_i\}_i}\colon \{b_i\}_i \to \boxtimes_i b_i$ are isomorphisms, and additionally $(ii)$ $\mathcal{B}'\left[\begin{smallmatrix}g(\underline{a})\\ \star\end{smallmatrix}\right]$ is representable iff $\mathcal{B}\left[\begin{smallmatrix}\underline{a}\\ \star\end{smallmatrix}\right]$ is representable.

*(Example)* A.2. Let $X$ be a topological space and $\mathsf{Open}(X)$ its poset of open subsets seen as a category. Then, the operad $\mathsf{Disj}(X)$ from Definition 2.3 is a pm-structure on $\mathsf{Open}(X)$. The functor $Z \mapsto \mathsf{Disj}(X)\left[\begin{smallmatrix}\{U_i\}_{i\in I}\\ Z\end{smallmatrix}\right]$ is representable if the open sets in $\{U_i\}_{i\in I}$ are pairwise disjoint, and in this case a representative is $\bigsqcup_i U_i \subseteq X$. In particular, when $I$ is empty, the representative is the empty subset $\varnothing \subseteq X$.

[ii] It is possible to define this structure by providing partially defined functors $\boxtimes^{(m)}\colon \mathtt{B}^{\times m} \dashrightarrow \mathtt{B}$, a monoidal unit (if it has one) and compatibilities. In other words, pm-categories are partial $\mathbb{E}_\infty$-algebras in $\mathtt{Cat}(\mathtt{sSet})$.

*Remark* A.3. Let $\mathtt{B}$ be a pm-$\mathtt{sSet}$-category. Then, the category of $\mathtt{B}_{\mathsf{opd}}$-algebras in $\mathtt{V}$ is equivalent to that of functors in $[\mathtt{B},\mathtt{V}]$ equipped with a lax monoidal structure, i.e.

$$\mathtt{Alg}_{\mathtt{B}_{\mathsf{opd}}}(\mathtt{V}) = \mathtt{LaxMon}\,[\mathtt{B},\mathtt{V}]\,.$$

From our previous definitions, this is a tautology. However, the idea is that, by means of the universal multimorphisms $(\Phi_{\underline{o}})_{\underline{o}}$, a lax monoidal structure on a functor $A\colon \mathtt{B}\to\mathtt{V}$ consists of providing operations $A(\Phi_{\underline{o}})\colon \bigotimes_i A(o_i)\to A(\boxtimes_i o_i)$ subject to natural conditions.

With this in mind, the following non-homotopical result is easy.

**Lemma A.4.** *Let* $\iota\colon \mathtt{B}\to\mathtt{N}$ *be a strong monoidal functor between pm-*$\mathtt{sSet}$*-categories. Then, the following square of decorated functors commutes*

$$\begin{array}{ccccccc}
[\mathtt{B},\mathtt{V}] & = & \mathtt{Alg}_{\overline{\mathtt{B}}}(\mathtt{V}) & \underset{\longleftarrow}{\overset{\iota_!}{\longrightarrow}} & \mathtt{Alg}_{\overline{\mathtt{N}}}(\mathtt{V}) & = & [\mathtt{N},\mathtt{V}] \\
 & & \uparrow{\scriptstyle \overline{(\star)}} & & \uparrow{\scriptstyle \overline{(\star)}} & & \\
\mathtt{LaxMon}[\mathtt{B},\mathtt{V}] & = & \mathtt{Alg}_{\mathtt{B}_{\mathsf{opd}}}(\mathtt{V}) & \underset{\longleftarrow}{\overset{\iota_\sharp}{\longrightarrow}} & \mathtt{Alg}_{\mathtt{N}_{\mathsf{opd}}}(\mathtt{V}) & = & \mathtt{LaxMon}[\mathtt{N},\mathtt{V}]
\end{array}\,.$$

*In other words, there is a natural isomorphism* $\iota_!\overline{A}\cong\overline{\iota_\sharp A}$.

*Proof.* By the universal property satisfied by $\iota_\sharp A$, one could verify the claim by equipping $\iota_!\overline{A}$ with a lax monoidal structure, see Remark A.3, and showing that it is universal. Recall that we can compute categorical left Kan extensions using coends, i.e. $\iota_!\overline{A}(o)=\int^{b\in\overline{\mathtt{B}}}\mathtt{N}\left[\begin{smallmatrix}\iota b\\ o\end{smallmatrix}\right]\otimes A(b)$. Then, the structure map $\mathtt{N}\left[\begin{smallmatrix}\{o_i\}_{i\in I}\\ u\end{smallmatrix}\right]\otimes\bigotimes_{i\in I}\iota_!\overline{A}(o_i)\longrightarrow\iota_!\overline{A}(u)$ we need to construct is given by the following composition

$$\begin{array}{c}
\mathtt{N}_{\mathsf{opd}}\left[\begin{smallmatrix}\{o_i\}_{i\in I}\\ u\end{smallmatrix}\right]\otimes\bigotimes_{i\in I}\left(\int^{b_i}\mathtt{N}_{\mathsf{opd}}\left[\begin{smallmatrix}\iota b_i\\ o_i\end{smallmatrix}\right]\otimes A(b_i)\right)\\
\cong\big\downarrow\ \text{Fubini, }\otimes\text{ is cocontinuous, symmetry}\\
\int^{(b_i)_i}\mathtt{N}_{\mathsf{opd}}\left[\begin{smallmatrix}\{o_i\}_{i\in I}\\ u\end{smallmatrix}\right]\otimes\bigotimes_{i\in I}\mathtt{N}\left[\begin{smallmatrix}\iota b_i\\ o_i\end{smallmatrix}\right]\otimes\bigotimes_{i\in I}A(b_i)\\
\big\downarrow\ \text{compose in }\mathtt{N}_{\mathsf{opd}}\text{, universal operation*}\\
\int^{(b_i)_i}\mathtt{N}_{\mathsf{opd}}\left[\begin{smallmatrix}\{\iota b_i\}_{i\in I}\\ u\end{smallmatrix}\right]\otimes A(\boxtimes_{i\in I}b_i)\\
\cong\big\downarrow\ \mathtt{N}\text{ is a pm-category*, coYoneda lemma}\\
\int^{(b_i)_i}\mathtt{N}\left[\begin{smallmatrix}\boxtimes_i\iota b_i\\ u\end{smallmatrix}\right]\otimes\left(\int^{b}\mathtt{B}\left[\begin{smallmatrix}b\\ \boxtimes_i b_i\end{smallmatrix}\right]\otimes A(b)\right)\\
\cong\big\downarrow\ \text{Fubini, }\otimes\text{ is cocontinuous, }\iota\text{ strong monoidal}\\
\int^{b}\left(\int^{(b_i)_i}\mathtt{N}\left[\begin{smallmatrix}\iota(\boxtimes_i b_i)\\ u\end{smallmatrix}\right]\otimes\mathtt{B}\left[\begin{smallmatrix}b\\ \boxtimes_i b_i\end{smallmatrix}\right]\right)\otimes A(b)\\
\big\downarrow\ \text{composition in }\mathtt{N}\circ(\mathrm{id}\otimes\iota)\\
\int^{b}\mathtt{N}\left[\begin{smallmatrix}\iota b\\ u\end{smallmatrix}\right]\otimes A(b).
\end{array}$$

Notice that the only case where this makes sense is when $\{o_i\}_i$ can be tensored in $\mathtt{N}$, i.e. $\mathtt{N}_{\mathsf{opd}}\left[{}^{\{o_i\}_i}_{\star}\right]$ is representable (otherwise, $\mathtt{N}_{\mathsf{opd}}\left[{}^{\{o_i\}_i}_{o}\right] = \mathbb{0}$ and the previous composition is the unique map from the initial object). Then, any family of objects $(b_i)_i$ in $\mathtt{B}$ such that $\mathtt{N}\left[{}^{\iota\, b_i}_{o_i}\right]$ is non-empty yields a collection $\{\iota\, b_i\}_i$ that can be tensored in $\mathtt{N}$; simply observe that in other case we would obtain a map of simplicial sets

$$\mathbb{0} \neq \left(\mathtt{N}_{\mathsf{opd}}\begin{bmatrix}\{o_i\}_i\\ \boxtimes_i\, o_i\end{bmatrix} \times \prod_i \mathtt{N}_{\mathsf{opd}}\begin{bmatrix}\iota\, b_i\\ o_i\end{bmatrix}\right) \xrightarrow{\text{composition}} \mathtt{N}_{\mathsf{opd}}\begin{bmatrix}\{\iota\, b_i\}_i\\ \boxtimes_i\, o_i\end{bmatrix} = \mathbb{0}.$$

By Condition (ii) in the definition of strong monoidal functor, this implies that $\{b_i\}_i$ can be tensored in $\mathtt{B}$. These considerations ensure that the steps marked with $^*$ in the previous construction are sound.

Checking that these choices fulfill the requirements is a lengthy but easy computation. See [39] for nomenclature and results about coends. □

Now we shall prove that the natural isomorphism $\iota_!\overline{A} \cong \overline{\iota_\sharp\, A}$ holds, possibly as an equivalence, when replacing $\iota_!$ and $\iota_\sharp$ by their homotopical analogues. Equivalently, we can compare derived Kan extensions (Proposition A.5). To achieve this goal, we adapt Horel's strategy from [36, §3]. Recall that we assume Hypothesis 1.2 on $\mathcal{V}$.

**Proposition A.5.** *Let $\iota\colon \mathtt{B} \to \mathtt{N}$ be a strong monoidal functor between pm-$\mathtt{sSet}$-categories. Then, the following square of decorated functors commutes*

$$\begin{array}{ccc}
\mathrm{Ho}\,\mathcal{A}\mathrm{lg}_{\overline{\mathtt{B}}}(\mathcal{V}) & \overset{\mathbb{L}\iota_!}{\rightleftarrows} & \mathrm{Ho}\,\mathcal{A}\mathrm{lg}_{\overline{\mathtt{N}}}(\mathcal{V}) \\
\uparrow \overline{(\star)} & & \uparrow \overline{(\star)} \\
\mathrm{Ho}\,\mathcal{A}\mathrm{lg}_{\mathtt{B}_{\mathsf{opd}}}(\mathcal{V}) & \overset{\mathbb{L}\iota_\sharp}{\rightleftarrows} & \mathrm{Ho}\,\mathcal{A}\mathrm{lg}_{\mathtt{N}_{\mathsf{opd}}}(\mathcal{V})
\end{array}.$$

*In other words, there is a natural equivalence $\mathbb{L}\iota_!\overline{A} \simeq \overline{\mathbb{L}\iota_\sharp\, A}$.*

*Proof.* It suffices to check that, for a proj-cofibrant $\mathtt{B}_{\mathsf{opd}}$-algebra QA, the canonical comparison map $\mathbb{L}\iota_!\overline{QA} \to \overline{\iota_\sharp\, QA}$ is an equivalence. Since $\mathcal{V}$ comes with a lax symmetric monoidal left Quillen functor from simplicial sets with the Kan-Quillen model structure, we know that QA is objectwise cofibrant by [43, Theorem 6.7] (note that we are using that $\mathcal{V}$ is symmetroidal). Thus, $\mathbb{L}\iota_!\overline{QA}(o)$ can be computed as the geometric realization of the Bar construction

$$\mathrm{Bar}_\bullet\left(\overline{\mathtt{N}}\left[{}^{\iota(\star)}_{o}\right], \overline{\mathtt{B}}, \overline{QA}\right) \equiv \left(\cdots \textstyle\coprod_{u,v} \overline{\mathtt{N}}\left[{}^{\iota(v)}_{o}\right] \otimes \overline{\mathtt{B}}\left[{}^{u}_{v}\right] \otimes QA(u) \rightrightarrows \coprod_u \overline{\mathtt{N}}\left[{}^{\iota(u)}_{o}\right] \otimes QA(u)\right),$$

as in the proof of [36, Proposition 3.9]. Following loc.cit., such geometric realization is weakly equivalent to $\iota_!\overline{QA}(o)$ (and hence to $\overline{\iota_\sharp\, QA}(o)$ by the previous lemma) if the functor $\star \otimes_{\overline{\mathtt{B}}} \overline{QA}\colon [\overline{\mathtt{B}}^{\mathsf{op}}, \mathtt{sSet}] \to \mathtt{V}$ preserves weak equivalences. We deduce this homotopical property making use of the same proof strategy as the third item in [36, Proposition 3.8], i.e. by rephrasing this tensor product over $\overline{\mathrm{B}}$ in terms of an operadic extension. More precisely, given a morphism $P \to P'$ in $[\overline{\mathrm{B}}^{op}, \mathrm{Set}]$, one builds a commutative square

$$\begin{array}{ccc} P \underset{\overline{\mathrm{B}}}{\otimes} \overline{QA} = \big((\alpha_P)_! \overline{QA}\big)(\infty) & \longrightarrow & \big((\alpha_{P'})_! \overline{QA}\big)(\infty) = P' \underset{\overline{\mathrm{B}}}{\otimes} \overline{QA} \\ \cong\big\downarrow & & \big\downarrow\cong \\ \big((\alpha_P)_\sharp\, QA\big)(\infty) & \longrightarrow & \big((\alpha_{P'})_\sharp\, QA\big)(\infty) \end{array},$$

where $\alpha_P$ denotes a canonical operad map $B_{opd} \to B_P$ whose target $(\mathrm{col}(B) \amalg \{\infty\})$-colored sSet-operad is defined as in loc.cit. with the only exception when the input list $\{b_i\}_{i\in I}$ satisfies $b_i \in \mathrm{col}(B)$ for all $i \in I$,

$$B_P \begin{bmatrix} \{b_i\}_{i\in I} \\ \infty \end{bmatrix} := \begin{cases} P(\boxtimes_i b_i) & \text{if } \{b_i\}_i \text{ can be tensored in } B, \\ \mathbb{0} & \text{otherwise.} \end{cases}$$

Notice that $B_P$ is a pm-category; indeed, $B_P \left[\begin{smallmatrix} \{c_i\}_{i\in I} \\ \star \end{smallmatrix}\right]$ is either $\mathbb{0}$, the representable functor $\overline{B}_P \left[\begin{smallmatrix} \boxtimes_i c_i \\ \star \end{smallmatrix}\right]$ if the list $\{c_i\}_i$ is of elements in $\mathrm{col}(B)$ and they can be tensored or $\overline{B}_P \left[\begin{smallmatrix} \infty \\ \star \end{smallmatrix}\right]$. Therefore, the vertical arrows in the square above are isomorphisms by Lemma A.4.

It remains to see that the lower horizontal map in the square is an equivalence when so is $P \to P'$. For this purpose, we think about it as the evaluation of the unit transformation

$$(\alpha_P)_\sharp\, QA \longrightarrow \beta^* \beta_\sharp (\alpha_P)_\sharp\, QA \cong \beta^* (\alpha_{P'})_\sharp\, QA$$

at $\infty \in B_P$, where $\beta\colon B_P \to B_{P'}$ is a canonical equivalence of operads induced by $P \to P'$ and $\mathrm{id}_B$ which satisfies $\beta\alpha_P = \alpha_{P'}$. We can conclude the result by using [43, Theorem 7.5] together with the fact that $(\alpha_P)_\sharp\, QA$ is proj-cofibrant since so is QA (note that, among other things, it has been applied that $\mathcal{V}$ is symmetric flat and symmetric h-monoidal). □

We close this appendix with a variation of [17, Theorem 4.52] that is required in §3.2.

**Proposition A.6.** *Let $\mathcal{V}$ be a closed symmetric monoidal model category satisfying Hypothesis 1.2. Let $\phi\colon \mathcal{O} \to \mathcal{P}$ be a map of weakly well-pointed or locally cofibrant operads*[iii]*. Then, the derived mate transformation $\mathbb{L}\phi_! \overline{A} \longrightarrow \overline{\mathbb{L}\phi_\sharp A}$ is an equivalence for any $\mathcal{O}$-algebra A if the canonical map $\mathbb{L}\phi_!\, \mathcal{O} \left[\begin{smallmatrix} \underline{c} \\ \star \end{smallmatrix}\right] \longrightarrow \mathcal{P} \left[\begin{smallmatrix} \phi(\underline{c}) \\ \star \end{smallmatrix}\right]$ is an equivalence in $[\overline{\mathcal{P}}, \mathcal{V}]$ for all $\underline{c}$.*

*Proof.* Without loss of generality, we can assume that A is cellular proj-cofibrant. In this case, the derived mate transformation is nothing but the nullary component $\underline{c} = \mathbb{0}$ of the canonical map $\mathbb{L}\phi_!\, \mathcal{O}_A \left[\begin{smallmatrix} \underline{c} \\ \star \end{smallmatrix}\right] \to \mathcal{P}_{\phi_\sharp A} \left[\begin{smallmatrix} \phi(\underline{c}) \\ \star \end{smallmatrix}\right]$ in $[\overline{\mathcal{P}}, \mathcal{V}]$. Hence, the claim will follow if we show more generally that this map is an equivalence for any $\underline{c}$. By the usual inductive argument for cell-attachments, it suffices to prove the following: given a pushout in $\mathcal{A}\mathrm{lg}_{\mathcal{O}}(\mathcal{V})$

$$\begin{array}{ccc} \bullet & \xrightarrow{\mathcal{O} \circ j} & \bullet \\ \big\downarrow & \ulcorner & \big\downarrow \\ A & \longrightarrow & A[j] \end{array},$$

[iii] Notice that under our simplifying assumption, the monoidal unit $\mathbb{I} \in \mathcal{V}$ is cofibrant, weakly well-pointed implies locally cofibrant.

where $j$ is a cofibration with cofibrant domain in $\mathcal{V}^{\mathrm{col}(\mathcal{O})}$ concentrated on one color (say $b$), the map $\mathbb{L}\phi_! \mathcal{O}_{A[j]} \left[\begin{smallmatrix}\underline{c}\\ \star\end{smallmatrix}\right] \to \mathcal{P}_{\phi_\sharp(A[j])} \left[\begin{smallmatrix}\phi(\underline{c})\\ \star\end{smallmatrix}\right]$ in $[\overline{\mathcal{P}}, \mathcal{V}]$ is an equivalence for any $\underline{c}$, provided the analogous maps for $A$ are so. Now, use [17, Lemma 3.22] to see the map we want to analyze as the colimit of an $\omega$-ladder in $[\overline{\mathcal{P}}, \mathcal{V}]$

$$\begin{array}{ccccccccc}
\mathbb{L}\phi_! \mathcal{O}_{A[j],0} \left[\begin{smallmatrix}\underline{c}\\ \star\end{smallmatrix}\right] & \longrightarrow & \cdots & \longrightarrow & \mathbb{L}\phi_! \mathcal{O}_{A[j],t-1} \left[\begin{smallmatrix}\underline{c}\\ \star\end{smallmatrix}\right] & \longrightarrow & \mathbb{L}\phi_! \mathcal{O}_{A[j],t} \left[\begin{smallmatrix}\underline{c}\\ \star\end{smallmatrix}\right] & \longrightarrow & \cdots \\
\downarrow & & & & \downarrow & & \downarrow & & \\
\mathcal{P}_{\phi_\sharp(A[j]),0} \left[\begin{smallmatrix}\phi(\underline{c})\\ \star\end{smallmatrix}\right] & \longrightarrow & \cdots & \longrightarrow & \mathcal{P}_{\phi_\sharp(A[j]),t-1} \left[\begin{smallmatrix}\phi(\underline{c})\\ \star\end{smallmatrix}\right] & \longrightarrow & \mathcal{P}_{\phi_\sharp(A[j]),t} \left[\begin{smallmatrix}\phi(\underline{c})\\ \star\end{smallmatrix}\right] & \longrightarrow & \cdots
\end{array}$$

whose first term is $\mathbb{L}\phi_! \mathcal{O}_A \left[\begin{smallmatrix}\underline{c}\\ \star\end{smallmatrix}\right] \to \mathcal{P}_{\phi_\sharp A} \left[\begin{smallmatrix}\phi(\underline{c})\\ \star\end{smallmatrix}\right]$ and whose successive steps $(t-1) \Rightarrow (t)$ fit into commutative cubes

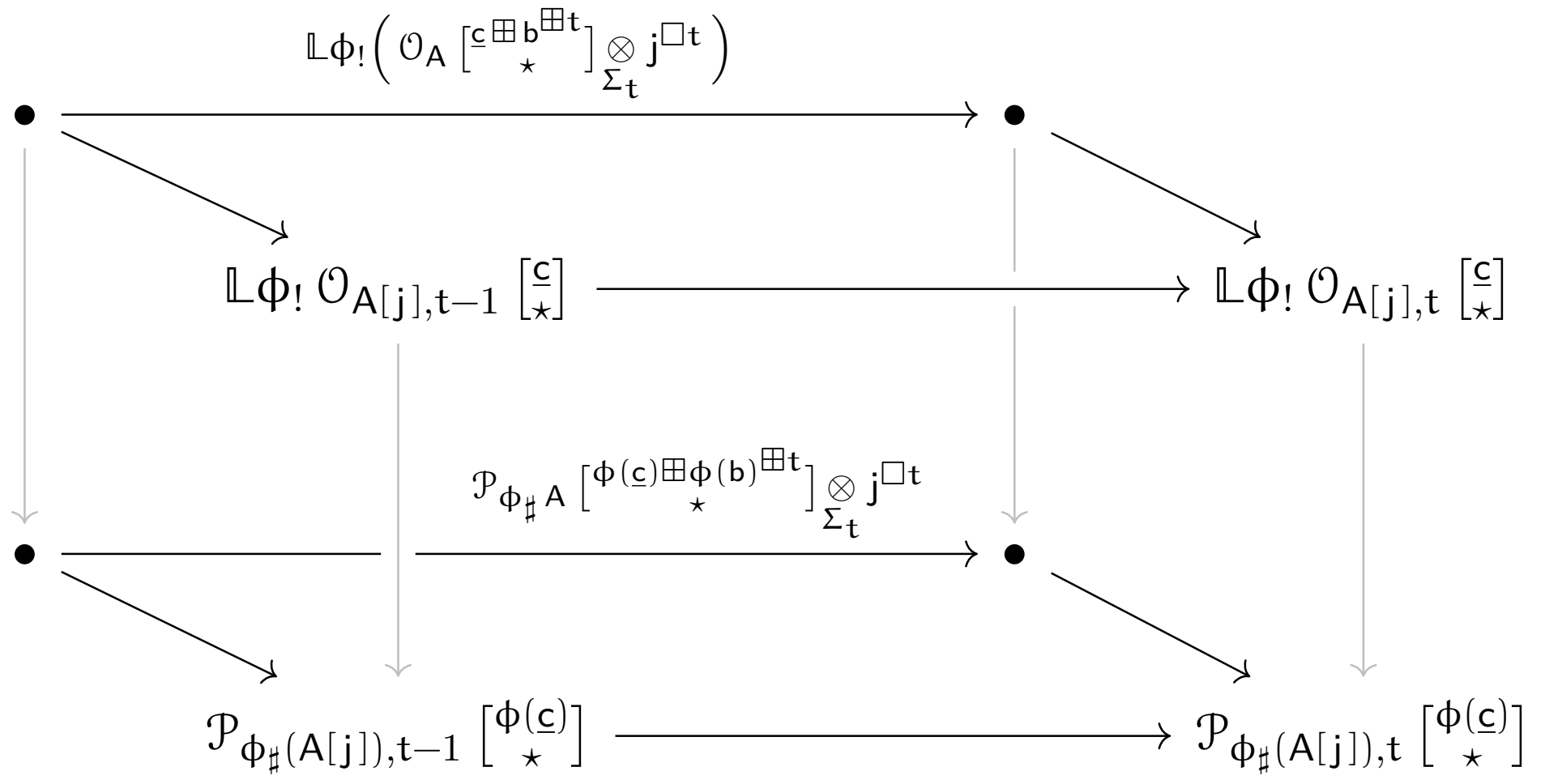

Thus, if we prove that (i) the black squares above (bottom and top faces of the cube) are homotopy pushouts, and (ii) the gray vertical arrows in the back face of the cube are equivalences, it would follow by induction over $t$ that $\mathbb{L}\phi_! \mathcal{O}_{A[j]} \left[\begin{smallmatrix}\underline{c}\\ \star\end{smallmatrix}\right] \to \mathcal{P}_{\phi_\sharp(A[j])} \left[\begin{smallmatrix}\phi(\underline{c})\\ \star\end{smallmatrix}\right]$ is an equivalence (notice that Hypotheses 1.2 imply that $\omega$-colimits in $\mathcal{V}$ preserve equivalences).

Let us begin with (i). Since $\mathbb{L}\phi_!$ preserves homotopy colimits, we are simply reduced to check, for any $d \in \mathrm{col}(\mathcal{O})$, that the pushout square in $\mathcal{V}$ (see [17, Lemma 3.20])

$$\begin{array}{ccc}
\bullet & \xrightarrow{\ \mathcal{O}_A \left[\begin{smallmatrix}\underline{c} \boxplus b^{\boxplus t}\\ d\end{smallmatrix}\right] \underset{\Sigma_t}{\otimes} j^{\square t}\ } & \bullet \\
\downarrow & \ulcorner & \downarrow \\
\mathcal{O}_{A[j],t-1} \left[\begin{smallmatrix}\underline{c}\\ d\end{smallmatrix}\right] & \longrightarrow & \mathcal{O}_{A[j],t} \left[\begin{smallmatrix}\underline{c}\\ d\end{smallmatrix}\right]
\end{array}$$

has cofibrant vertices and its top horizontal map is a cofibration (the same argument applies to $\mathcal{P}_{\phi_\sharp A}$ in place of $\mathcal{O}_A$). The symmetroidal condition on $\mathcal{V}$ from Hypotheses 1.2, the cofibrant assumptions on $\mathcal{O}$ (weakly well-pointed or locally cofibrant operad) and $j$ being a cofibration with cofibrant domain imply these claims (see the proof of [17, Lemma 3.35] and [43, Theorem 6.7]).

Regarding (ii), observe that symmetric flatness on $\mathcal{V}$ implies that $\star \otimes_{\Sigma_t} (\mathrm{co})\mathrm{dom}(j^{\square t})$ sends proj-equivalences in the category $[\Sigma_t^{\mathrm{op}}, \mathcal{V}]$ to equivalences in $\mathcal{V}$, i.e. it is a derived functor. Since $\mathrm{codom}(j^{\square t}) = \mathrm{codom}(j)^{\otimes t}$, the codomain case follows directly from the definition of symmetric flatness. For $\mathrm{dom}(j^{\square t})$, a simple induction using symmetric flatness and Harper's filtration [31, Definition 7.10], recalling that $\mathrm{dom}(j^{\square t}) = Q^t_{t-1}$ in his notation, does the job. To form $\mathbb{L}\phi_!\big(\mathcal{O}_A\left[\begin{smallmatrix}\underline{c}\boxplus \underline{b}^{\boxplus t}\\ \star\end{smallmatrix}\right]\otimes_{\Sigma_t}(\mathrm{co})\mathrm{dom}(j^{\square t})\big)$ note that we see $\mathcal{O}_A\left[\begin{smallmatrix}\underline{c}\boxplus \underline{b}^{\boxplus t}\\ \star\end{smallmatrix}\right]$ as an object in $[\Sigma_t^{\mathrm{op}}\times\overline{\mathcal{O}},\mathcal{V}]$. Thus, the derived left Kan extension $\mathbb{L}\phi_!$ and the derived functor $\star\otimes_{\Sigma_t}(\mathrm{co})\mathrm{dom}(j^{\square t})$ affect different variables. Therefore, we obtain a commutative triangle

$$\begin{array}{ccc}
\mathbb{L}\phi_!\Big(\mathcal{O}_A\left[\begin{smallmatrix}\underline{c}\boxplus \underline{b}^{\boxplus t}\\ \star\end{smallmatrix}\right]\otimes_{\Sigma_t}(\mathrm{co})\mathrm{dom}(j^{\square t})\Big) & \xrightarrow{\sim} & \Big(\mathbb{L}\phi_!\,\mathcal{O}_A\left[\begin{smallmatrix}\underline{c}\boxplus \underline{b}^{\boxplus t}\\ \star\end{smallmatrix}\right]\Big)\otimes_{\Sigma_t}(\mathrm{co})\mathrm{dom}(j^{\square t}) \\
\downarrow & \swarrow & \\
\mathcal{P}_{\phi_\sharp A}\left[\begin{smallmatrix}\phi(\underline{c})\boxplus\phi(\underline{b})^{\boxplus t}\\ \star\end{smallmatrix}\right]\otimes_{\Sigma_t}(\mathrm{co})\mathrm{dom}(j^{\square t}) & &
\end{array}$$

where the maps we need to analyze for $(ii)$ appear as the gray vertical arrows, and the diagonal maps are equivalences since, by assumption, $\mathbb{L}\phi_!\,\mathcal{O}_A\left[\begin{smallmatrix}\underline{c}'\\ \star\end{smallmatrix}\right]\to\mathcal{P}_{\phi_\sharp A}\left[\begin{smallmatrix}\phi(\underline{c}')\\ \star\end{smallmatrix}\right]$ is an equivalence for any $\underline{c}'$.

□

# B On the left properness of prefactorization algebras

This appendix is devoted to prove that Hypothesis 4.11 holds true in relevant contexts for practitioners of factorization algebras, e.g. when $\mathcal{V} = \mathcal{C}h(R)$ is the projective model structure on chain complexes over a commutative ring $R$ containing $\mathbb{Q}$. This result might seem surprising since it was claimed in [43, §9.7], without proof, that the projective model structure $\mathrm{p}\mathcal{F}\mathrm{act}_X(\mathcal{V}) = \mathcal{A}\mathrm{lg}_{M_X}(\mathcal{V})$ is not left proper. Also notice that projective model structures on operadic algebras are rarely left-proper.

Our argument relies on the explicit enveloping operad computation, Proposition B.1, for prefactorization algebras. See, e.g., [17, §3.1] for background and definitions around enveloping operads. Here, we simply recall that $M_{X,\mathcal{F}}\left[\begin{smallmatrix}\underline{U}\\ V\end{smallmatrix}\right]$ sits in a reflexive coequalizer

$$M^1_{X,\mathcal{F}}\left[\begin{smallmatrix}\underline{U}\\ V\end{smallmatrix}\right] \underset{d_1}{\overset{d_0}{\rightrightarrows}} M^0_{X,\mathcal{F}}\left[\begin{smallmatrix}\underline{U}\\ V\end{smallmatrix}\right] \xrightarrow{\ \mathrm{coeq}\ } M_{X,\mathcal{F}}\left[\begin{smallmatrix}\underline{U}\\ V\end{smallmatrix}\right]\ ,$$

where: (1) $M^i_{X,\mathcal{F}}\left[\begin{smallmatrix}\underline{U}\\ V\end{smallmatrix}\right]$, for $i = 0, 1$, is given by a colimit, over a particular groupoid of decorated trees, of a natural tagging functor built from the underlying $\Sigma_{\mathrm{Open}(X)}$-collection of $M_X$ and the prefactorization algebra $\mathcal{F}$; (2) $d_0$ is induced from the action of $M_X$ on $\mathcal{F}$; and (3) $d_1$ comes from the operadic structure on $M_X$. To provide some intuition, let us describe $M^0_{X,\mathcal{F}}$ concretely; $M^1_{X,\mathcal{F}}$ is just a slightly more complicated version of this which requires trees of height $\leqslant 1$ (see [17, Definition 3.8]).

Let $O = \mathrm{Open}(X)$ denote the set of open subsets of $X$, and consider the groupoid $G := (\mathrm{Fin}\downarrow O\times\{\downarrow, \wr\})^{\simeq}$ of finite sets equipped with a function into $O\times\{\downarrow,\wr\}$ and compatible bijections between them. We want to think about the elements of $G$ as the leaves in an $O$-colored corolla with root labeled by $V\in O$. The function into $\{\downarrow,\wr\}$ just splits these leaves into two subsets, let us say snaky and straight leaves to match the notation of [17]. By selecting the snaky arrows only, we get a functor $G \longrightarrow (\mathsf{Fin}_O)^{\simeq}$ explicitly defined as follows: $\big((Z,g)\colon I \to O\times\{\downarrow,\wr\}\big) \longmapsto \big(Z|\colon g^{-1}\{\wr\} \to O\big)$. Looking at $\underline{U}$ as an element of $(\mathsf{Fin}_O)^{\simeq}$, one can finally define the groupoid used to construct $M^0_{X,\mathcal{F}}\left[\begin{smallmatrix}\underline{U}\\ V\end{smallmatrix}\right]$. That is, we take the slice of the previous functor over $\underline{U}$, $G\left[\begin{smallmatrix}\underline{U}\\ V\end{smallmatrix}\right] := G \downarrow \underline{U}$. Then, we set $M^0_{X,\mathcal{F}}\left[\begin{smallmatrix}\underline{U}\\ V\end{smallmatrix}\right]$ to be the colimit of

$$G\left[\begin{smallmatrix}\underline{U}\\ V\end{smallmatrix}\right]^{op} \xrightarrow{(i)} G^{op} \times O = \left\{\begin{smallmatrix}O\text{-corollas with}\\ \text{two types of leaves}\end{smallmatrix}\right\}^{op} \xrightarrow{(ii)} \mathtt{V},$$

where $(i)$ is the slice projection in the first factor and the constant functor with value $V$ on the second, and $(ii)$ is the natural tagging functor

$$(M_X;\mathcal{F})\colon \Big(\big((Z,g)\colon I \to O\times\{\downarrow,\wr\}\big), V\Big) \longmapsto M_X\begin{bmatrix}\{Z_i\}_{i\in I}\\ V\end{bmatrix} \cdot \bigotimes_{j\in g^{-1}\{\downarrow\}} \mathcal{F}(Z_j)$$

mentioned above (see [17, Definition 3.2], recalling that $M_X$ is an operad in sets).

**Proposition B.1.** *Let $\mathcal{F} \in \mathtt{pFact}_X(\mathtt{V})$ be a prefactorization algebra over $X$. Then, the enveloping operad $M_{X,\mathcal{F}}$ has $\mathtt{V}$-objects of multimorphisms*

$$M_{X,\mathcal{F}}\begin{bmatrix}\underline{U}\\ V\end{bmatrix} = \begin{cases}\mathcal{F}(V_{\underline{U}}) & \textit{if } \bigsqcup_i U_i \subseteq V,\\[1ex] \mathbb{0} & \textit{otherwise},\end{cases}$$

*where $V_{\underline{U}} := V \setminus \bigcup_i \overline{U_i}$ is biggest open subset of $V$ disjoint to all the $U_i$'s. Furthermore, its operadic composition is canonically induced from the prefactorization products of $\mathcal{F}$.*

*Proof.* We show that this is the case by extending the (reflexive) coequalizer which defines $M_{X,\mathcal{F}}\left[\begin{smallmatrix}\underline{U}\\ V\end{smallmatrix}\right]$ to a split coequalizer

$$M^1_{X,\mathcal{F}}\left[\begin{smallmatrix}\underline{U}\\ V\end{smallmatrix}\right] \underset{d_1}{\overset{d_0}{\rightrightarrows}} M^0_{X,\mathcal{F}}\left[\begin{smallmatrix}\underline{U}\\ V\end{smallmatrix}\right] \xrightarrow{\ d\ } M_{X,\mathcal{F}}\left[\begin{smallmatrix}\underline{U}\\ V\end{smallmatrix}\right] = \mathcal{F}(V_{\underline{U}}) \ ,$$

$t\colon M^0_{X,\mathcal{F}}\left[\begin{smallmatrix}\underline{U}\\ V\end{smallmatrix}\right] \dashrightarrow M^1_{X,\mathcal{F}}\left[\begin{smallmatrix}\underline{U}\\ V\end{smallmatrix}\right]$ $s\colon M_{X,\mathcal{F}}\left[\begin{smallmatrix}\underline{U}\\ V\end{smallmatrix}\right] \dashrightarrow M^0_{X,\mathcal{F}}\left[\begin{smallmatrix}\underline{U}\\ V\end{smallmatrix}\right]$

i.e. we need to construct maps $s, t$ as depicted above satisfying $ds = \mathrm{id}$, $d_1 t = \mathrm{id}$, $sd = d_0 t$. Even though a categorical proof is possible, we perform the necessary checks using elements for the sake of conciseness. For this purpose, let us denote $\mu\left[\begin{smallmatrix}\underline{Z}\\ V\end{smallmatrix}\right]$ the unique element of $M_X\left[\begin{smallmatrix}\underline{Z}\\ V\end{smallmatrix}\right]$, which transforms under an $O$-colored permutation $\sigma\colon \underline{Z} \to \underline{Z}\cdot\sigma$ as $\mu\left[\begin{smallmatrix}\underline{Z}\\ V\end{smallmatrix}\right]\cdot\sigma = \mu\left[\begin{smallmatrix}\underline{Z}\cdot\sigma\\ V\end{smallmatrix}\right]$, and let us use square brackets to denote classes in the quotients $M^i_{X,\mathcal{F}}\left[\begin{smallmatrix}\underline{U}\\ V\end{smallmatrix}\right]$, for $i = 0, 1$. Then,

- The map $d\colon M^0_{X,\mathcal{F}}\left[\begin{smallmatrix}\underline{U}\\ V\end{smallmatrix}\right] \longrightarrow \mathcal{F}(V_{\underline{U}})$ is given by

$$\Big[\mu\left[\begin{smallmatrix}\underline{Z}\\ V\end{smallmatrix}\right], \textstyle\bigotimes_i f_{W_i}\Big] = \Big[\big(\mu\left[\begin{smallmatrix}(V_{\underline{U}};\underline{U})\\ V\end{smallmatrix}\right] \circ_1 \mu\left[\begin{smallmatrix}\underline{W}\\ V_{\underline{U}}\end{smallmatrix}\right]\big), \textstyle\bigotimes_i f_{W_i}\Big] \longmapsto \mu\left[\begin{smallmatrix}\underline{W}\\ V_{\underline{U}}\end{smallmatrix}\right]\big(\textstyle\bigotimes_i f_{W_i}\big),$$

  which uses the identities $\mu\left[\begin{smallmatrix}\underline{Z}\\ V\end{smallmatrix}\right] = \mu\left[\begin{smallmatrix}(\underline{W};\underline{U})\\ V\end{smallmatrix}\right]\cdot \mathsf{sh} = \big(\mu\left[\begin{smallmatrix}(V_{\underline{U}};\underline{U})\\ V\end{smallmatrix}\right] \circ_1 \mu\left[\begin{smallmatrix}\underline{W}\\ V_{\underline{U}}\end{smallmatrix}\right]\big)\cdot \mathsf{sh}$, where $\mathsf{sh}\colon (\underline{W};\underline{U}) \to \underline{Z}$ is the unique $O$-colored shuffle yielding $\underline{Z}$ by shuffling the straight leaves (labeled by $\underline{W}$) and the snaky leaves (labeled by $\underline{U}$).

- The map $s\colon \mathcal{F}(V_{\underline{U}}) \longrightarrow M^0_{X,\mathcal{F}}\left[\begin{smallmatrix}\underline{U}\\ V\end{smallmatrix}\right]$ is just the inclusion of a factor in the colimit; explicitly,

$$f_{V_{\underline{U}}} \longmapsto \left[\mu\left[\begin{smallmatrix}(V_{\underline{U}};\underline{U})\\ V\end{smallmatrix}\right], f_{V_{\underline{U}}}\right].$$

- The map $t\colon M^0_{X,\mathcal{F}}\left[\begin{smallmatrix}\underline{U}\\ V\end{smallmatrix}\right] \longrightarrow M^1_{X,\mathcal{F}}\left[\begin{smallmatrix}\underline{U}\\ V\end{smallmatrix}\right]$, similarly to $d$, is given by

$$\left[\mu\left[\begin{smallmatrix}Z\\ V\end{smallmatrix}\right], \bigotimes_i f_{W_i}\right] = \left[\left(\mu\left[\begin{smallmatrix}(V_{\underline{U}};\underline{U})\\ V\end{smallmatrix}\right] \circ_1 \mu\left[\begin{smallmatrix}W\\ V_{\underline{U}}\end{smallmatrix}\right]\right), \bigotimes_i f_{W_i}\right] \longmapsto \left[\mu\left[\begin{smallmatrix}(V_{\underline{U}};\underline{U})\\ V\end{smallmatrix}\right], \mu\left[\begin{smallmatrix}W\\ V_{\underline{U}}\end{smallmatrix}\right], \bigotimes_i f_{W_i}\right].$$

Regarding the equations that must be satisfied, one readily observes $ds = \mathrm{id}$ and $sd = d_0 t$, since

$$sd\left[\mu\left[\begin{smallmatrix}Z\\ V\end{smallmatrix}\right], \bigotimes_i f_{W_i}\right] = \left[\mu\left[\begin{smallmatrix}(V_{\underline{U}};\underline{U})\\ V\end{smallmatrix}\right], \mu\left[\begin{smallmatrix}W\\ V_{\underline{U}}\end{smallmatrix}\right]\left(\bigotimes_i f_{W_i}\right)\right] = d_0 t\left[\mu\left[\begin{smallmatrix}Z\\ V\end{smallmatrix}\right], \bigotimes_i f_{W_i}\right].$$

The remaining identity $d_1 t = \mathrm{id}$ holds because

$$\left[\mu\left[\begin{smallmatrix}(W;\underline{U})\\ V\end{smallmatrix}\right], \bigotimes_i f_{W_i}\right] = \left[\mu\left[\begin{smallmatrix}(W;\underline{U})\\ V\end{smallmatrix}\right]\cdot \mathrm{sh}, \bigotimes_i f_{W_i}\right] = \left[\mu\left[\begin{smallmatrix}Z\\ V\end{smallmatrix}\right], \bigotimes_i f_{W_i}\right],$$

for sh the O-colored shuffle described above.

Identifying the operadic structure on $M_{X,\mathcal{F}}$ is just a routinary task. □

**Theorem B.2.** *Assume that $\mathcal{V}$ satisfies, apart from Hypothesis 1.2, the conditions: (i) the model category $\mathcal{V}$ is left proper, (ii) finite coproducts in $\mathcal{V}$ are homotopical, and (iii) taking coinvariants in $\mathcal{V}$ for any finite group action is homotopical. Then, the following model structures on $\mathtt{pFact}_X(V)$ are left proper:*

- *the projective model structure $\mathrm{p}\mathcal{F}\mathrm{act}_X(\mathcal{V})$;*
- *the model structure $\mathrm{p}\mathcal{F}\mathrm{act}^{\mathrm{lc}}_X(\mathcal{V})$ obtained by localizing $\mathrm{p}\mathcal{F}\mathrm{act}_X(\mathcal{V})$ at the set $\widetilde{\mathcal{S}}_X$ defined right before Theorem 4.21, or its analogue for $\mathrm{E}_X$ in place of $\mathrm{D}_X$.*

*In particular, this holds when $\mathcal{V} = \mathcal{C}\mathrm{h}(R)$ for $R$ a commutative ring such that $\mathbb{Q} \subseteq R$.*

*Proof.* For the left properness of $\mathrm{p}\mathcal{F}\mathrm{act}_X(\mathcal{V})$, we apply [17, Proposition 4.36 and Remark 4.40]. The hypotheses in loc.cit. which are not directly included in the previous list follow easily from our assumptions: (a) pretty smallness can be replaced by the fact that, in $\mathcal{V}$, filtered colimits are homotopical; (b) cofibrant objects being flat is the simplest case of symmetric flatness (which holds in $\mathcal{V}$). Furthermore, the explicit computation in Proposition B.1 implies that any proj-equivalence of prefactorization algebras $\mathcal{F} \xrightarrow{\sim} \mathcal{F}'$ induces an equivalence of enveloping operads $M_{X,\mathcal{F}} \xrightarrow{\sim} M_{X,\mathcal{F}'}$. Hence, $\mathrm{p}\mathcal{F}\mathrm{act}_X(\mathcal{V})$ is left proper relative to all its objects (i.e. it is left proper).

For the second dot, recall that left Bousfield localization preserves left properness. □

*Remark* B.3. Using symmetric flatness and tractability for $\mathcal{V}$, it is even possible to drop condition (iii) in Theorem B.2 using that $\star \otimes_{\Sigma_t} (\mathrm{co})\mathrm{dom}(j^{\square t})$ is a derived functor for any generating cofibration $j$ (see the proof of Proposition A.6). We do not demonstrate this fact in order to keep the appendix relatively short.